\documentclass[12pt,a4paper]{article}
\usepackage{amsmath}
\usepackage{amssymb}
\usepackage{amsxtra}
\usepackage{amscd}
\usepackage{amsthm}
\usepackage{hyperref}
\usepackage[mathscr]{eucal} 
\usepackage{graphicx}
\usepackage{mathtools}
\DeclareMathOperator{\MyProd}{\scalebox{1.4}{$\mathrm{I\kern-0.2ex I}$}}

\theoremstyle{plain}

\newtheorem{sbthm}[subsubsection]{Theorem}
\newtheorem{sbprop}[subsubsection]{Proposition}

\newtheorem{sblem}[subsubsection]{Lemma}

\theoremstyle{definition}

\newtheorem{sbrem}[subsubsection]{Remark}

\newtheorem{sbpara}[subsubsection]{}

\newenvironment{pf}{\proof[\proofname]}{\endproof}

\begin{document}

\title{Height functions for motives, II}

\author
{Kazuya Kato}

\maketitle

\begin{abstract}

This is the Part II of our paper ''Height functions for motives''. We consider more general period domains and the height functions on more general sets of motives. We also consider the corresponding Hodge theoretic variant of Nevanlinna theory.

\end{abstract}

\newcommand{\lr}[1]{\langle#1\rangle}
\newcommand{\ul}[1]{\underline{#1}}
\newcommand{\eq}[2]{\begin{equation}\label{#1}#2 \end{equation}}
\newcommand{\ml}[2]{\begin{multline}\label{#1}#2 \end{multline}}
\newcommand{\ga}[2]{\begin{gather}\label{#1}#2 \end{gather}}
\newcommand{\mc}{\mathcal}
\newcommand{\mb}{\mathbb}
\newcommand{\surj}{\twoheadrightarrow}
\newcommand{\inj}{\hookrightarrow}
\newcommand{\red}{{\rm red}}
\newcommand{\codim}{{\rm codim}}
\newcommand{\rank}{{\rm rank}}
\newcommand{\Pic}{{\rm Pic}}
\newcommand{\Div}{{\rm Div}}
\newcommand{\divi}{{\rm div}}
\newcommand{\Hom}{{\rm Hom}}
\newcommand{\Ext}{{\rm Ext}}
\newcommand{\im}{{\rm im}}
\newcommand{\fil}{{\rm fil}}
\newcommand{\gp}{{\rm gp}}
\newcommand{\Spec}{{\rm Spec}}
\newcommand{\Sing}{{\rm Sing}}
\newcommand{\Char}{{\rm char}}
\newcommand{\Tr}{{\rm Tr}}
\newcommand{\Gal}{{\rm Gal}}
\newcommand{\Min}{{\rm Min}}
\newcommand{\Upsilonlt}{{\rm mult}}
\newcommand{\Max}{{\rm Max}}
\newcommand{\Alb}{{\rm Alb}}
\newcommand{\gr}{{\rm gr}}
\newcommand{\Ker}{{\rm Ker}}
\newcommand{\Lie}{{\rm Lie}}
\newcommand{\infi}{{\rm inf}}
\newcommand{\et}{{\rm \acute{e}t}}
\newcommand{\pole}{{\rm pole}}
\newcommand{\ti}{\times }
\newcommand{\modu}{{\rm mod}}
\newcommand{\Ab}[1]{{\mathcal A} {\mathit b}/#1}
\newcommand{\sA}{{\mathcal A}}
\newcommand{\sB}{{\mathcal B}}
\newcommand{\sC}{{\mathcal C}}
\newcommand{\sD}{{\mathcal D}}
\newcommand{\sE}{{\mathcal E}}
\newcommand{\sF}{{\mathcal F}}
\newcommand{\sG}{{\mathcal G}}
\newcommand{\sH}{{\mathcal H}}
\newcommand{\sI}{{\mathcal I}}
\newcommand{\sJ}{{\mathcal J}}
\newcommand{\sK}{{\mathcal K}}
\newcommand{\sL}{{\mathcal L}}
\newcommand{\sM}{{\mathcal M}}
\newcommand{\sN}{{\mathcal N}}
\newcommand{\sO}{{\mathcal O}}
\newcommand{\sP}{{\mathcal P}}
\newcommand{\sQ}{{\mathcal Q}}
\newcommand{\sR}{{\mathcal R}}
\newcommand{\sS}{{\mathcal S}}
\newcommand{\sT}{{\mathcal T}}
\newcommand{\sU}{{\mathcal U}}
\newcommand{\sV}{{\mathcal V}}
\newcommand{\sW}{{\mathcal W}}
\newcommand{\sX}{{\mathcal X}}
\newcommand{\sY}{{\mathcal Y}}
\newcommand{\sZ}{{\mathcal Z}}
\newcommand{\A}{{\mathbb A}}
\newcommand{\B}{{\mathbb B}}
\newcommand{\C}{{\mathbb C}}
\newcommand{\D}{{\mathbb D}}
\newcommand{\E}{{\mathbb E}}
\newcommand{\F}{{\mathbb F}}
\newcommand{\G}{{\mathbb G}}
\renewcommand{\H}{{\mathbb H}}
\newcommand{\I}{{\mathbb I}}
\newcommand{\J}{{\mathbb J}}
\newcommand{\M}{{\mathbb M}}
\newcommand{\N}{{\mathbb N}}
\renewcommand{\P}{{\mathbb P}}
\newcommand{\Q}{{\mathbb Q}}
\newcommand{\R}{{\mathbb R}}
\newcommand{\T}{{\mathbb T}}
\newcommand{\U}{{\mathbb U}}
\newcommand{\V}{{\mathbb V}}
\newcommand{\W}{{\mathbb W}}
\newcommand{\X}{{\mathbb X}}
\newcommand{\Y}{{\mathbb Y}}
\newcommand{\Z}{{\mathbb Z}}
\newcommand{\pic}{{\text{Pic}(C,\sD)[E,\nabla]}}
\newcommand{\ocd}{{\Omega^1_C\{\sD\}}}
\newcommand{\oc}{{\Omega^1_C}}
\newcommand{\al}{{\alpha}}
\newcommand{\be}{{\beta}}
\newcommand{\ta}{{\theta}}
\newcommand{\ve}{{\varepsilon}}
\newcommand{\phe}{{\varphi}}
\newcommand{\om}{{\overline M}}
\newcommand{\sym}{{\text{Sym}(\om)}}
\newcommand{\an}{{{\rm an}}}
\newcommand{\Ad}{{{\rm Ad}}}
\newcommand{\bs}{{\backslash}}
\newcommand{\lra}{\longrightarrow}
\newcommand{\Sig}{{\Sigma}}
\newcommand{\sig}{{\sigma}}
\newcommand{\dR}{{{\rm dR}}}
\newcommand{\fg}{{\frak {g}}}
\newcommand{\Rep}{{\rm {Rep}}}
\newcommand{\La}{{{\Lambda}}}
\newcommand{\Mh}{{{\rm Mor}_{\rm hor}}}
\newcommand{\la}{{{\lambda}}}
\newcommand{\hor}{{{\rm hor}}}
\newcommand{\pst}{{{\rm pst}}}

{\bf Contents.}

\medskip

\S0. Introduction

\S1. Period domains and motives

\S2. Curvature forms and  Hodge theory

\S3. Height functions 

\S4. Speculations

\medskip

\setcounter{section}{-1}
\section{Introduction}
This is a sequel of our paper \cite{Ka3} concerning height functions for mixed motives over number fields, which we call Part I. The new subjects in this Part II are as follows.

\medskip

1. We consider more general period domains $X(\C)$ and more general sets $X(F)$ of motives over number fields $F$ than Part I. See Section 1.2 for the definition of $X(\C)$ and Section 1.3 for the definition of $X(F)$ of this paper. 

\medskip

2. We start a Hodge theoretic variant of Nevanlinna theory, which we would call {\it Hodge-Nevanlinna theory}. 

\medskip

Vojta compares 
number theory and Nevanlinna theory (see \cite{Vo}, for example) . He compares height functions

\medskip

(1)  $V(F) \to \R_{>0}\;$ in number theory,

(2) $\sM(B, V(\C))\to \R\;$ in Nevanlinna theory. 

\medskip

Here $F$ is a number field, $V$ is an algebraic variety over $F$, $V(F)$ and $V(\C)$ denote the sets of $F$-points and $\C$-points of $V$, respectively, $B$ is a connected one-dimensional complex analytic manifold (that is, a connected Riemann surface) endowed with a finite flat morphism $B\to \C$, and we denote by $\sM(B, V(\C))$ the set of meromorphic functions from $B$ to $V(\C)$. If $V$ is a dense open subvariety of a projective variety $\bar V$, a meromorphic function from $B$ to $V(\C)$ is nothing but a holomorphic function $f: B\to \bar V(\C)$ such that $f(B)\not\subset \bar V(\C)\smallsetminus V(C)$. For $f\in \sM(B, V(\C))$, Nevanlinna theory asks how often $f$ has singularities on $B$, that is, 
how often we have $f(z)\in \bar V(\C)\smallsetminus V(\C)$ ($z\in B$), relating this question to height functions of $f$. 
(In the classical Nevanlinna theory, $B=\C$ endowed with the identity morphism $B\to \C$, $\bar V={\bf P}^1(\C)$, and Nevanlinna theory asks how often  a meromorphic function on $\C$ has values in a given finite subset of ${\bf P}^1(\C)$.)
\medskip

In this paper, we compare the height functions in the above (1) and (2) and our height functions

\medskip

(I) $X(F)\to \R_{>0}\;$ in number theory,

(II) $\sM_{\text{hor}}(B, X(\C))\to \R\;$ in Hodge-Nevanlinna theory,

(III) $ \sM_{\text{hor}}(C, X(\C))\to \R\;$ in Hodge theory.

\medskip

Here $F$ and $B$ are as above, $C$ in (III) is a connected compact Riemann surface  (that is, a smooth projective curve over $\C$), $X(F)$ is a set of motives over $F$, $X(\C)$ is a period domain which classifies Hodge structures, and for $Y=B$ or $C$, $\sM_{\text{hor}}(Y, X(\C))$ denotes the set of horizontal meromorphic functions $f$ from $Y$ to $X(\C)$. See  \ref{hormer} for the precise definition of $\sM_{\text{hor}}(Y, X(\C))$. 
 If $\bar X(\C)\supset X(\C)$ denotes the toroidal partial compactification of $X(\C)$ consisting of the classes of nilpotent orbits of rank $\leq 1$ (\cite{KU}, \cite{KNU} Part III, \cite{KP}, \cite{KNU2}), $\sM_{\hor}(Y, X(\C))$ is identified with the set of horizontal morphisms $f:Y\to \bar X(\C)$ such that $f(Y)\not\subset \bar X(\C)\smallsetminus X(\C)$. In our Hodge-Nevanlinna theory (II), for $Y=B$, we ask how often $f$ has singularities as a meromorphic function from $B$ to $X(\C)$, that is,  how often we have $f(z)\in \bar X(\C)\smallsetminus X(\C)$ ($z\in B$), relating this question to Hodge theoretic height functions of $f$.

 (I) and (III) were already considered in Part I,
 but as is said above, we consider more general situation in this Part II.  (In Part I, we used the notation $\Mh(C, \bar X(\C))$ for $\sM_{\hor}(C, X(\C))$ in (III), where $\text{Mor}$ stands for morphisms whereas $\sM$ stands for meromorphic maps, regarding elements of $\sM_{\hor}(C, X(\C))$ as horizontal morphisms $C\to \bar X(\C)$.) In this Part II, $X(F)$ is the set of $G$-mixed motives over $F$ and $X(\C)$ is the set of $G$-mixed Hodge structures, where $G$ is a linear algebraic group over $\Q$. Here a $G$-mixed motive (resp. $G$-mixed Hodge structure) means an exact $\otimes$-functor from the category of linear representations of $G$ to the category of mixed motives (resp. mixed Hodge structures). Mumford-Tate domains (\cite{GGK}) appear as standard examples of $X(\C)$, and higher Albanese manifolds  \cite{Ha} also appear as examples of $X(\C)$. For various Shimura varieties, the set of their $\C$-points give special cases of $X(\C)$ and their $F$-points give points of  $X(F)$. 
  In (II) (resp. (III)), for $Y=B$ (resp. $C$), $\sM_{\text{hor}}(Y, X(\C))$ is understood as set of isomorphism classes of  variations of $G$-mixed Hodge structure  with logarithmic degeneration on $Y$. 
 
   We expect that the comparisons of (1), (2), (I), (II), (III) shed new lights to the arithmetic of the world of motives.
 Since $X(F)$ is usually not the set of $F$-points of an algebraic variety, the study of (I) is not reduced to the study of (1). Since $\bar X(\C)$ is usually not a complex analytic space (it is a ''logarithmic manifold'' in the sense of \cite{KU}, 3.5.7), the study of (II) is not reduced to the study of (2).  
As in  Section \ref{s2.6},  the height functions in (I), (II), (III) are related not only philosophically, but also actually via asymptotic formulas.

The organization of this paper is as follows. In Section 1, we define $X(\C)$ and $X(F)$.  Section 2 is a preparation for Hodge-Nevanlinna theory. There we study some Hodge theoretic curvature forms (theorems \ref{spadepos}, \ref{htp3}). In Section 3, we define height functions in (I), (II), (III). In Section 4, we present questions. 

 Acknowledgements.   The author thanks Teruhisa Koshikawa for helpful discussions and valuable advice. He thanks Minhyong Kim and Paul Vojta for teaching him various aspects in Diophantine geometry. 
  He also thanks Spencer Bloch, Chikara Nakayama, Sampei Usui, and Takeshi Saito for valuable advice.

 The author is partially supported by NFS grants DMS 1303421 and DMS 1601861.

\section{Period domains and motives}\label{s1}

 In Section 1.2 (resp. Section 1.3), we define the set $X(\C)$ (resp. $X(F)$) which appears in Section 0. Sections 1.1 is a preparation for Section 1.2. In Section 1.4, we give some examples.  The contents of Section 1.1 and 1.2 are contained in the theory of Mumford-Tate domains in \cite{GGK} if the linear algebraic group $G$ is reductive and if $\sG=G$ (\ref{HAB}).  

\subsection{The period domain $D$}\label{s1.1}

In this Section \ref{s1.1}, we consider period domains $D(G, \Upsilon)$ (see \ref{defD1} for the definition) and their generalizations $D(G, \sG, H_b)$ (see \ref{defD2} for the definition). In the case $G$ is a reductive algebraic group, $D(G, \Upsilon)$ is the Mumford-Tate domain studied in \cite{GGK}. For a linear algebraic group $G$ in general, $D(G,\Upsilon)$ is the period domain considered in \cite{KNU2}.

\begin{sbpara}\label{1.1a} 
In \ref{1.1a}--\ref{modff}, we prepare notation.

For a linear algebraic group $G$ over $\Q$, let $G_u$ be the unipotent radical of $G$ and let $G_{\text{red}}$ be the reductive group $G/G_u$.

  Let $\Rep(G)$ be the category of 
finite-dimensional linear representations of $G$ over $\Q$.

\end{sbpara}

\begin{sbpara}\label{1.1b} Let $\Q$MHS (resp. $\R$MHS) be the category of mixed $\Q$ (resp. $\R$)-Hodge structures. 

Let $\Q$HS (resp. $\R$HS) be the full subcategory of $\Q$MHS (resp. $\R$MHS) consisting of objects which are direct sums of pure objects.

For $H\in \Q$MHS (resp. $\R$MHS), let $H_\Q$ (resp. $H_\R$) be the underlying $\Q$ (resp. $\R$)-vector space of $H$.

On the other hand, for $H\in \Q$MHS, let $\R \otimes_\Q H$ be the associated mixed $\R$-Hodge structure. For such $H$, $H_\R$ denotes the  $\R$-vector space $\R\otimes_\Q H_\Q$ underlying the mixed $\R$-Hodge structure $\R\otimes_\Q H$ though this may be a bit confusing.  

\end{sbpara}

\begin{sbpara}\label{modff}
For a commutative ring $R$, let $\text{Mod}_{ff}(R)$ be the category of finitely generated free $A$-modules.
In particular, for a field
$k$, $\text{Mod}_{ff}(k)$ is the category of finite-dimensional $k$-vector spaces. 
\end{sbpara}

\begin{sbpara}\label{Gandw} In the rest of this paper,  let $G$ be a linear algebraic group over $\Q$.  
We assume that we are given  a homomorphism $w: {\bf G}_m\to G_{\text{red}}=G/G_u$ whose image is in the center of $G_{\text{red}}$ such that for some (equivalently,  for any) lifting ${\bf G}_m\to G$ of $w$, the adjoint action of ${\bf G}_m$ on $\Lie(G_u)$ is of weight $\leq -1$.

Note that any $V\in \Rep(G)$ is endowed with a canonical $G$-stable increasing filtration, which we denote by $W_{\bullet}V$ and call the weight filtration, such that for any lifting ${\bf G}_m\to G$ of $w$ and for $i\in \Z$, $W_iV$ is the part of weight $\leq i$ of $V$ for the action of ${\bf G}_m$. 
\end{sbpara}

 \begin{sblem}\label{Hred} Let $H:\Rep(G)\to \Q$MHS be an exact $\otimes$-functor. Then 
for any $V\in \text{Rep}(G_{\text{red}})\subset \text{Rep}(G)$, we have $H(V)\in \Q\text{HS}\subset \Q\text{MHS}$.

\end{sblem}

\begin{pf} Since any representation of the reductive group $G_{\text{red}}$ is semi-simple, the weight filtration of $V\in \Rep(G)$ splits. 
\end{pf}

\begin{sbpara}\label{SCR} 

  Let $S_{\C/\R}$ be the Weil restriction of the multiplicative group ${\mathbb G}_m$ from $\C$ to $\R$, which represents the functor $R\mapsto (\C\otimes_\R R)^\times$ for commutative rings $R$ over $\R$.
 Let $w: {\mathbb G}_{m,\R}\to S_{\C/\R}$ be the homomorphism which represents the natural maps $R^\times\to (\C\otimes_\R R)^\times$ 
 for commutative rings $R$ over $\R$.
 
 As in \cite{De71}, the category $\R\text{HS}$ (\ref{1.1b}) is equivalent to the category of 
 finite-dimensional linear representations  of $S_{\C/\R}$ over $\R$. For a finite-dimensional representation $V$ of $S_{\C/\R}$ over $\R$, the corresponding object of $\R$HS  has $V$ as the underlying $\R$-structure and has the Hodge decomposition $$V_\C:=\C\otimes_\R V=\bigoplus_{p,q\in \Z} \; V_\C^{p,q},$$ where
$$V_\C^{p,q}=\{v\in V_\C\;|\; [z]v= z^p{\bar z}^qv\;\text{for}\;z\in \C^\times\}.$$
Here $[z]$ denotes $z$ regarded as an element of $S_{\C/\R}(\R)=\C^\times$. 
For a finite-dimensional representation $V$ of $S_{\C/\R}$ over $\R$, the part of $V$ of weight $w$ of the corresponding Hodge structure coincides with the part of $V$ of weight $w$ for the action of ${\bf G}_{m,\R}$ via ${\bf G}_{m,\R}\overset{w}\to S_{\C/\R}$.
\end{sbpara}

\begin{sbpara}\label{hb}
Consider a $G_{\text{red}}(\R)$-conjugacy class $\Upsilon$ of homomorphisms $S_{\C/\R}\overset{h}\to G_{\text{red},\R}$ of algebraic groups over $\R$ 
such that the 
composition ${\bf G}_{m,\R}\overset{w}\to S_{\C/\R}\overset{h}\to G_{\red,\R}$ coincides with the homomorphism induced by $w: {\bf G}_m\to G_{\text{red}}$ (\ref{Gandw}).  

\end{sbpara}

\begin{sbpara}\label{defD1} For $\Upsilon$ as in \ref{hb}, we define the period domain
$D(G, \Upsilon)$ to be the set of isomorphism classes of exact $\otimes$-functors $H: \Rep(G) \to \Q$MHS preserving the underlying $\Q$-vector spaces with weight filtrations satisfying the following condition (i).

\medskip

(i) The homomorphism $S_{\C/\R}\to G_{\red,\R}$ corresponding to $\Rep(G_{\red})\overset{H}\to \R$HS via the theory of Tannakian categories
belongs to $\Upsilon$.   

\medskip

In the case $G$ is reductive, this set $D(G, \Upsilon)$ is identified with $\Upsilon$ itself, and it is a Mumford-Tate domain studied in \cite{GGK}.
In general, this set $D(G,\Upsilon)$ coincides with the set which is denoted by $D$ in \cite{KNU2} associated to $\Upsilon$.

\end{sbpara}

 \begin{sbpara}\label{nonemp}
 The set $D(G,\Upsilon)$ is not empty. In fact, there is a homomorphism $S_{\C/\R}\to G_\R$ such that the composition $S_{\C/\R}\to G_\R\to G_{\red,\R}$ belongs to $\Upsilon$. This homomorphism defines an exact $\otimes$-functor $\Rep(G)\to \Q$MHS which is an element of $D(G, \Upsilon)$.

  \end{sbpara}

\begin{sbpara}\label{HAB}

 Consider an algebraic normal subgroup $\sG$ of $G$ (defined over $\Q$). Let $\sQ=G/\sG$. 

 Consider an exact $\otimes$-functor $$H_b:\Rep(G) \to \Q\text{MHS}$$ which keeps the underlying $\Q$-vector spaces with weight filtrations.  
 (The subscript $b$ in  $H_b$ presents the role of $H_b$ as a base point.)

 \end{sbpara}

  \begin{sbpara}\label{defD2}   For $(\sG, H_b)$ as in \ref{HAB}, we define the period domain $D(G, \sG, H_b)$ as follows.
 
 Let $\Upsilon$ be the $G_{\red}(\R)$-conjugacy class of the homomorphism $S_{\C/\R}\to G_{\text{red},\R}$  corresponding to the exact $\otimes$-functor $\Rep(G_{\text{red}})\to \R$HS induced by $H_b$ (\ref{Hred}). Let $H_b|_{\sQ}$ be the composite exact $\otimes$-functor $\Rep(\sQ)\to \Rep(G)\overset{H_b}\to \Q$MHS and let 
  $\Upsilon_{\sQ}$ be the $\sQ_{\text{red}}(\R)$-conjugacy class of homomorphisms $S_{\C/\R}\to \sQ_\R$ induced by $H_b|_{\sQ}$.
 Define $D=D(G, \sG, H_b)\subset D(G, \Upsilon)$ to be the inverse image of $\text{class}(H_b|_{\sQ})\in D(\sQ, \Upsilon_{\sQ})$ under the canonical map $D(G, \Upsilon)\to D(\sQ, \Upsilon_{\sQ})\;;\;\text{class}(H)\mapsto \text{class}(H|_{\sQ})$. 

Note that we have $\text{class}(H_b)\in D$.

\end{sbpara}

\begin{sbpara}\label{sG=G1} In fact, $D(G, \Upsilon)$ is regarded as the case $\sG=G$ of $D(G,\sG, H_b)$. For $\Upsilon$ as in \ref{hb}, by fixing any $H_b\in D(G, \Upsilon)$, we have $D(G,\Upsilon)=D(G,G, H_b)$. On the other hand, for $H_b$ as in \ref{HAB}, we have $D(G,G, H_b)=D(G,\Upsilon)$ where $\Upsilon$ is determined by $H_b$ as in \ref{defD2}.  

\end{sbpara}

 We consider the complex analytic structures of these period domains. We first review the complex analytic structure of $D(G, \Upsilon)$ which is given in \cite{GGK}, \cite{KNU2}.

\begin{sblem}\label{11lem}
$G(\R)G_u(\C)$ acts on $D(G, \Upsilon)$ transitively.

\end{sblem}

\begin{pf}  \cite{KNU2} Prop. 2.2.5. 
\end{pf}

\begin{sbpara}\label{DGan} 
The set  $D(G, \Upsilon)$ is regarded as a complex analytic manifold as follows.

Note that $H\in D(G, \Upsilon)$ is regarded as an exact $\otimes$-functor from $\Rep(G)$ to the category $\sC\supset \Q$MHS of finite-dimensional $\Q$-vector spaces $V$ such that $V$ is endowed with an increasing filtration (called the weight filtration) and $V_\C$ is endowed with a decreasing filtration (called the Hodge filtration). By \ref{11lem}, in the set of isomorphism classes of exact $\otimes$-functors from $\Rep(G)$ to $\sC$, we have a unique
 $G(\C)$-orbit $\check{D}(G, \Upsilon)$  containing $D(G, \Upsilon)$ 
 Then $G(\C)$ acts on $\check{D}(G,\Upsilon)$ transitively and the isotropy group in $G(\C)$ of each point of $\check{D}(G, \Upsilon)$ is an algebraic subgroup of $G(\C)$. 
Hence $\check{D}(G, \Upsilon)$ is a complex analytic manifold. 

By the following lemma, $D(G, \Upsilon)$ is also a complex analytic manifold.

\end{sbpara}

\begin{sblem}\label{open}
 $D(G, \Upsilon)$ is open in $\check{D}(G, \Upsilon)$ 
 
 \end{sblem} 
 
 \begin{pf} (\cite{GGK}, \cite{KNU2}) Prop. 2.2.7.
 \end{pf}

\begin{sbpara}\label{Drep} Let $\sA$ be the category of complex analytic spaces in the sense of Grothendieck (the structure sheaf can have non-zero nilpotent sections). For a commutative ring $R$ and for $Y\in \sA$, let $\text{Mod}_{ff}(R,Y)$ be the category of locally constant sheaves of finitely generated free  $R$-modules on $Y$. Let $\Q$MHS$(Y)$ be the category of pairs $(\sH_\Q, \fil)$, where  $\sH_\Q\in \text{Mod}_{ff}(\Q, Y)$ endowed with an increasing filtration by locally constant $\Q$-subsheaves (called the weight filtration) and $\fil$ is a decreasing filtration on $\sH_{\sO}:=\sO_Y\otimes_\Q \sH_\Q$ by subbundles (called the Hodge filtration) such that for each $y\in Y$, the fiber $(\sH_{\sQ,y}, \fil(y))$ at $y$ is a $\Q$-mixed Hodge structure.

Then the complex analytic manifold $D(G, \Upsilon)$  represents the functor $\sA\to$ (Sets) which sends $Y\in \sA$ to the set of isomorphism classes of exact $\otimes$-functors $\sH:\Rep(G)\to \Q$MHS$(Y)$ such that $\sH(V)_\Q$ with the weight filtration for $V\in \Rep(G)$ coincides with the constant sheaf $V$ with the weight filtration and such that the following (i) is satisfied.

\medskip

(i) For any $y\in Y$, the homomorphism $S_{\C/\R}\to G_{\text{red},\R}$ induced by the $\otimes$-functor $\Rep(G_{\text{red}})\to \R\text{HS}\;;\;V\mapsto \R\otimes_\Q \sH(V)(y)$ belongs to $\Upsilon$. 

\medskip

For $D=D(G, \Upsilon)$, consider the universal object $\sH_D\in \Q$MHS$(D)$. Note that this object need not satisfy the Griffiths transversality, that is, it need not be a variation of mixed Hodge structure. For $V\in\Rep(G)$,  we have the universal Hodge filtration on $\sO_D\otimes_\Q \sH_D(V)_\Q= \sO_D\otimes_\Q V$.
\end{sbpara}

\begin{sbprop}\label{tanD}  Let $D=D(G,\Upsilon)$.  Concerning the tangent bundle $T_D$ of $D$, we have a canonical isomorphism 
$$T_D\cong (\sO_D\otimes_\Q \Lie(G))/\fil^0.$$
Here $\fil$ denotes the universal Hodge filtration on $\sO_D\otimes_\Q \sH_D(\Lie(G))_\Q = \sO_D \otimes_\Q \Lie(G)$ where  $\Lie(G)$ is endowed with the adjoint action of $G$. 
\end{sbprop}

\begin{pf} Let $U$ be an open set of $D$ and let $\tilde U=U[\epsilon]/(\epsilon^2)$ 
be the complex analytic space whose underlying topological space is the same as that of $U$  and whose structure sheaf is $\sO_U[\epsilon]/(\epsilon^2)$ (with $\epsilon$ an indeterminate). By the usual infinitesimal understanding of the tangent bundle, $T_D(U)$ is identified with the set of all morphisms $f: \tilde U \to D$ such that the composition $U \to \tilde U \overset{f}\to D$ is the inclusion morphism. 
 Since $D$ represents the functor described in  \ref{Drep}, this set is  identified with the set $\Gamma(U, (\sO_U\otimes_{\Q} \Lie(G))/\fil^0)$ as follows. Let $h$ be an element of the last set. Then $h$ corresponds to 
 $f$ which sends  $V\in \Rep(G)$ to $(V, \fil)\in \Q$MHS$(\tilde U)$  where $\fil$ is as follows. Let $\fil'$ be the universal filtration on $\sO_U\otimes_{\Q} V$. 
  Locally, lift $h$ to $\sO_{\tilde U}\otimes_\Q \Lie(G)$. Then 
  $\fil^r= \{x+ \epsilon h(x)\;|\; x\in (\fil')^r\}+ \epsilon (\fil')^r$, where $h(x)$ is defined by the action of $\Lie(G)$ on $V$. 
\end{pf}

\begin{sbprop}\label{smooth}

The morphism $D(G, \Upsilon)\to D(\sQ, \Upsilon_{\sQ})$ is smooth.
\end{sbprop}
\begin{pf} This follows from the surjectivity of the induced maps of tangent spaces (\ref{tanD}), which follows from the surjectivity of $\Lie(G)\to \Lie(\sQ)$.  \end{pf}

\begin{sbpara}  
By definition, $D(G, \sG, H_b)$ is a closed complex analytic subspace of $D(G, \Upsilon)$.

\end{sbpara}

\begin{sbprop}\label{Dan} Let $D=D(G, \sG, H_b)$.

(1) $D$ is a complex analytic manifold.

(2) The tangent bundle of $D$ is canonically isomorphic to $(\sO_D \otimes_\Q \Lie(\sG))_{\sO}/\fil^0$. 
Here $\fil$ is the Hodge filtration on $\sO_D \otimes_\Q \sH_D(\Lie(\sG))_\Q= \sO_D\otimes_\Q \Lie(\sG)$ where $\Lie(\sG)$ is endowed with the adjoint action of $G$.

\end{sbprop}

\begin{pf} This follows from \ref{smooth}.
\end{pf}

\begin{sbprop}\label{forb}
 $D(G, \sG, H_b)$ is a finite disjoint union of $\sG(\R)\sG_u(\C)$-orbits. 
 
 \end{sbprop}
 
 \begin{pf} 
First we consider the case $G$ is reductive. For $x\in D(G, \sG, H_b)$, let $\varphi_x: S_{\C/\R}\to G(\R)$ be the homomorphism associated to $x$. Let $x\in D(G, \sG, H_b)$ be the class of $H_b$, and let $J$ (resp. $J_{\sQ}$) be the algebraic subgroup of $G$ (resp. $\sQ$) consisting of elements $a$ such that $a\varphi_x(z)a^{-1}= \varphi_x(z)$ (resp. $a\varphi_{x,\sQ}(z)a^{-1}= \varphi_{x,\sQ}(z)$ where 
$\varphi_{x,\sQ}$ denotes the composition $S_{\C/\R}\overset{\varphi_x}\to G_\R\to \sQ_\R$), and let $I$ be the inverse image of $J_{\sQ}$ in $G$. Then $J\subset I$ and $\sG\subset I$. 

\medskip

Claim 1. We have $I(\R)/J(\R)\overset{\cong}\to D(G, \sG, H_b)\;;\; g \mapsto gx.$

\medskip

In fact, if $y\in D(G,\sG, H_b)$, $y=gx$ for some $g\in G(\R)$ in $D(G, \Upsilon)$ where $\Upsilon$ is associated to $x$. Since $\varphi_x$ and $\varphi_y$ induce the same homomorphism $S_{\C/\R}\to \sQ_{\R}$, we see that the image of $g$ in $\sQ(\R)$ belongs to $J_{\sQ}(\R)$. This proves Claim 1.

 \medskip

Claim 2. The map $I(\R) \to D(G, \sG, H_b)$ induces a surjection $I(\R)/I(\R)^{\circ}\to \sG(\R)\bs D(G, \sG, H_b)$ where $I(\R)^{\circ}$ denotes the connected component of $I(\R)$ containing $1$. 
\medskip

Proposition \ref{forb}  in the case $G$ is reductive follows from Claim 2 and the finiteness of $I(\R)/I(\R)^{\circ}$. 

We prove Claim 2. Let $g_1, g_2\in I(\R)$ and assume $g_2=g_1g$ for some $g\in I(\R)^{\circ}$. We prove $g_2\in\sG(\R) g_1J(\R)$. The image $g_{\sQ}$ of $g$ in $J_{\sQ}(\R)$ belongs to the connected component of $J_{\sQ}(\R)$ containing $1$ and hence $g_{\sQ}=\exp(A)$ for some  $A\in \Lie(J_{\sQ})_\R$. 
Since $G$ is reductive, the surjection $\Lie(G)\to \Lie(\sQ)$ is regarded as the projection to a direct factor of the Lie algebra $\Lie(G)$, and hence the map $\Lie(J)\to \Lie(J_{\sQ})$ is surjective. 
Take an element $\tilde A\in \Lie(J)_\R$ which is sent to $A$ and let $a= \exp(\tilde A) \in J(\R)$. Then 
 $h:= ga^{-1} $ belongs to $\sG(\R)$. We have $g_2= g_1ha= g_1hg_1^{-1} g_1 a$. Since $\sG$ is normal in $G$, $g_1hg_1^{-1}\in \sG(\R)$.

Now we do not assume $G$ is reductive. By the case $G$ is reductive treated above, it is sufficient to prove that if $y_1, y_2\in D(G, \sG, H_b)$ and if the images of $y_1$ and $y_2$ in $D(G_{\text{red}}, \Upsilon)$ coincide, then $y_2=gy_1$ for some $g\in \sG_u(\C)$. Take a homomorphism $h: S_{\C/\R}\to G_\R$ which lifts the common homomorphism $S_{\C/\R}\to G_{\red,\R}$ induced by $y_1$ and by $y_2$ and let $H$ be the element of $D(G, \Upsilon)$ corresponding to $h$. Then $H$ induces $\Rep(G)\to \R$HS. By \cite{CKS} 2.20, $y_j=e^{i\delta(y_j)}u(y_j)H$ for unique $u(y_j)\in G_u(\R)$ and $\delta(y_j)\in \Lie(G_u)_\R$ such that for any $V\in \Rep(G)$ and for any $p,q\in \Z$, $\delta(y_i)$ sends the $(p,q)$-Hodge component  of $u(y_j)H(V)$  into the sum of $(p',q')$-Hodge components of $u(y_j)H(V)$ with $p'<p$, $q'<q$. Since $y_{1,\sQ}=y_{2,\sQ}$ in $D(\sQ, \Upsilon_{\sQ})$, we have $u(y_1)_{\sQ}=u(y_2)_{\sQ}$ in $\sQ_u(\R)$ and $\delta(y_1)_{\sQ}= \delta(y_2)_{\sQ}$ in $\Lie(\sQ_u)_\R$. Thus $y_2=gy_1$ with $g=e^{i\delta(y_2)}u(y_2)u(y_1)^{-1}e^{-i\delta(y_1)}$ and $g_{\sQ}=1$. Hence $g\in \sG_u(\C)$. 
\end{pf}

\begin{sbpara} The action of $\sG(\R)\sG_u(\C)$ on $D(G, \sG, H_b)$ need not be transitive. 

Example. Let $G= J \cup J\sig\subset GL_2$, where $$J=\{\begin{pmatrix} a &-b \\ b& a\end{pmatrix}\in GL_2\}, \quad \sig=\begin{pmatrix} 0& 1 \\ 1 & 0\end{pmatrix}.$$ Let $w: {\bf G}_m\to G$ be the canonical embedding as scaler matrices. Consider the $G(\R)$-conjugacy class $\Upsilon=\{h_{+}, h_{-}\}$ where $h_{\pm}: S_{\C/\R}\to J \subset G$ is defined by 
$$h_{\pm}(a+bi)=  \begin{pmatrix} a&\mp b\\ \pm b & a\end{pmatrix}.$$ 
We have $D(G, \Upsilon)=\{h_+,h_-\}$. Let $\sG= G\cap SL_2$ and let $H_b=h_+\in D(G, \sG)$.  Then $\text{det}: \sQ=G/\sG \overset{\cong}\to {\bf G}_m$ 
 and hence $D(\sQ,\Upsilon_{\sQ})$ consists of one element. Hence $D(G, \sG, H_b)$ consists of two elements $h_+$ and $h_-$.  But $\sG(\R)\cong SO(2,\R)$ is connected and $\sG_u=\{1\}$, and hence the action of $\sG(\R)\sG_u(\C)$ on $D(G,\sG, H_b)$ is trivial and is not transitive. 
\end{sbpara}

\subsection{The period domain  $X(\C)$}\label{s1.2}

In this Section \ref{s1.2}, we consider the period domain $X(\C)=X_{G, \sG, H_b, K}(\C)$ (see \ref{defXC2} for the definition). We also describe a special case  $X_{G, \Upsilon, K}(\C)$ (see \ref{defXC1}) of $X_{G, \sG, H_b, K}(\C)$ hoping that it works as a guide to the more complicated object $X_{G, \sG, H_b, K}(\C)$.

\begin{sbpara}

Let ${\bf A}^f_\Q$ be the non-archimedean component of the adele ring of $\Q$.

\end{sbpara}

\begin{sbpara}\label{defXC1}  Assume we are given $\Upsilon$ as in \ref{hb}. Let $K$ be an open compact subgroup of $G({\bf A}^f_\Q)$. 
Define $$X_{G,\Upsilon, K}(\C):= G(\Q) \bs  (D(G,\Upsilon)\times (G({\bf A}^f_\Q)/K)).$$
Here $G(\Q)$ acts on $D(G,\Upsilon)\times (G({\bf A}^f_\Q)/K)$ diagonally from the left.

\end{sbpara}

\begin{sbpara}\label{defXC2} Assume that we are given $(\sG, H_b)$ as in \ref{HAB}.  
Let $K$ be an open compact subgroup of $\sG({\bf A}^f_\Q)$.

Define
$$X(\C)=X_{G,\sG, H_b, K}(\C):=\sG(\Q)\bs (D(G, \sG, H_b) \times (\sG({\bf A}^f_\Q)/K)),$$
where $\sG(\Q)$ acts on $D(G, \sG,H_b)\times (\sG({\bf A}^f_\Q)/K)$ diagonally from the left. 

\end{sbpara}

\begin{sbpara}\label{sG=G2} $X_{G,\Upsilon,K}(\C)$ is the case $\sG=G$ of $X_{G,\sG, H_b, K}(\C)$ (\ref{sG=G1}). 
\end{sbpara}

\begin{sbpara}\label{defXCa} As is explained in \ref{defXCac} below, 
$X_{G,\Upsilon,K}(\C)$ is identified with 
the set of isomorphism classes of pairs $(H, \lambda)$, where 

\medskip

$H$ is an exact $\otimes$-functor from $\Rep(G)$ to $\Q$MHS 

\medskip
\noindent 
such that for some
isomorphism of $\otimes$-functors from $\Rep(G)$ to $\text{Mod}_{ff}(\Q)$
$$\theta: (V \mapsto H(V)_\Q)\overset{\cong}\to (V\mapsto V)\quad\text{preserving the weight filtrations},$$
the homomorphism $S_{\C/\R}\to G_{\red}(\R)$ induced by the restrictions of $H$ and $\theta$ to $\Rep(G_{\red})$  belongs to $\Upsilon$, and 

\medskip

$\lambda$ (called the $K$-level structure) is a mod $K$ class of an isomorphism of $\otimes$-functors from $\Rep(G)$ to $\text{Mod}_{ff}({\bf A}^f_\Q)$
$$\tilde \lambda: (V\mapsto {\bf A}^f_\Q \otimes_\Q V) \overset{\cong}\to (V\mapsto {\bf A}^f\otimes_\Q H(V)_\Q)\quad\text{preserving the weight filtrations}.$$

\end{sbpara}

\begin{sbpara}\label{defXCac} 
In \ref{defXCa}, we identify the class of a pair
 $(H,\la)$  as above with   $\text{class}(H', g)\in X_{G,\Upsilon,K}(\C)$, where $H'=H$ endowed with the identification of $H(V)_{\Q}$ with $V$ via $\theta$, and $g=({\bf A}^f_\Q\otimes \theta) \circ \tilde \la\in G({\bf A}^f_\Q)$. Conversely, we identify $\text{class}(H',g)\in X_{G,\Upsilon,K}(\C)$   ($\text{class}(H')\in D$, $g\in G({\bf A}^f_\Q)$) with the class of 
  the pair $(H, \la)$, where 
$H=H'$,  $\tilde \la=g$. 

\end{sbpara}

\begin{sbpara}\label{defXCb} As is explained in \ref{defXCc} below, the set $X_{G,\sG,H_b,K}(\C)$ is identified with the set of isomorphism classes of triples $(H, \xi, \lambda)$, where 

\medskip

$H$ is an exact $\otimes$-functor from $\Rep(G)$ to $\Q$MHS,  

\medskip

$\xi$ is an isomorphism of $\otimes$-functors from $\Rep(\sQ)$ to $\Q$MHS
$$H|_{\sQ}\overset{\cong}\to H_b|_{\sQ},$$
where $H|_{\sQ}$ and $H_b|_{\sQ}$ denote the restrictions of $H$ and $H_b$  to $\Rep(\sQ)\subset \Rep(G)$, respectively, 

\medskip
$\la$ is as in \ref{defXCa},

\medskip
\noindent
such that there are an isomorphism of $\otimes$-functors from $\Rep(G_{\text{red}})$ to $\R$HS
$$\nu: (V\mapsto \R\otimes_\Q H(V)) \overset{\cong}\to (V\mapsto \R \otimes_\Q H_b(V))$$
 and an isomorphism of $\otimes$-functors from $\Rep(G)$ to $\text{Mod}_{ff}(\Q)$ 
$$\theta: (V \mapsto H(V)_\Q)\overset{\cong}\to (V\mapsto V)\quad\text{preserving the weight filtrations}$$
which satisfy the following (i)--(iii).

(i) $\xi$ and $\nu$ induce the same isomorphism of functors $(V\mapsto \R \otimes_\Q H(V)) \overset{\cong}\to (V\mapsto \R \otimes_\Q H_b(V))$ from $\Rep(\sQ_{\text{red}})$ to $\R$HS.

(ii) $\xi$ and $\theta$ induce the same isomorphism of functors $(V\mapsto H(V)_\Q)\overset{\cong}\to (V\mapsto H_b(V)_\Q)$ from $\Rep(\sQ)$ to $\text{Mod}_{ff}(\Q)$.

(iii) The automorphism  $({\bf A}^f_\Q\otimes \theta) \circ \tilde \la$ of the $\otimes$-functor 
$$\Rep(G) \to \text{Mod}_{ff}({\bf A}^f_\Q)\;;\; V\mapsto {\bf A}^f_\Q \otimes_\Q V$$ 
(which is an element of $G({\bf A}^f_\Q)$ by the theory of Tannakian categories) belongs to $\sG({\bf A}^f_\Q)$.

\end{sbpara}

\begin{sbpara}\label{defXCc}  
In  \ref{defXCb}, we identify the class of a triple $(H, \xi, \la)$ as above with   $\text{class}(H', g)\in X_{G,\sG,H_b,K}(\C)$, where $H'=H$ endowed with the identification of $H(V)_{\Q}$ with $V$ via $\theta$, and $g=({\bf A}^f_\Q\otimes \theta) \circ \tilde \la\in \sG({\bf A}^f_\Q)$. Conversely, we identify $\text{class}(H',g)\in X_{G,\sG,H_b,K}(\C)$   ($\text{class}(H')\in D$, $g\in \sG({\bf A}^f_\Q)$) with the class of 
  the  triple $(H, \xi, \la)$, where 
$H=H'$,  $\tilde \la=g$, and $\xi$ is the evident isomorphism.

In the case $\sG=G$, the identification $X_{G,\sG, H_b,K}(\C)= X_{G, \Upsilon, K}(\C)$ corresponds to $(H, \xi, \la)\mapsto (H,\la)$ ($\xi$ is unique and is the evident isomorphism in this case). 

\end{sbpara}

Let $(\sG, H_b, K)$ be as in \ref{defXC2}. 
We define a structure of a complex manifold on $X(\C)=X_{G,\sG, H_b,K}(\C)$ assuming \ref{1.2.1}, \ref{1.2.1b} (concerning polarization) and assuming \ref{neat} (concerning level structure) below.

\begin{sbpara}\label{1.2.1} 

We assume  we are given 
$V_0\in \Rep(G)$, a non-degenerate bilinear form $$\langle\;,\;\rangle_{0,w}: \gr^W_wV_0\times \gr^W_wV_0\to \Q\cdot (2\pi i)^{-w}$$ for each $w\in \Z$ which is symmetric if $w$ is even and anti-symmetric if $w$ is odd, and a homomorphism $$\eta: G\to {\bf G}_m$$
having the following properties (i)--(iii).

\medskip

(i) $\langle gx,gy\rangle_{0,w}= \eta(g)^w\langle x, y\rangle_{0,w}$ for any $w\in \Z, \; g\in G, \;x,y\in \gr^W_wV_0$.

(ii) The composition ${\bf G}_m\overset{w}\to G_{\text{red}}\overset{\eta}\to {\bf G}_m$ is $x\mapsto x^2$. 

(iii)  The homomorphism $G\to \text{Aut}(V_0) \times {\bf G}_m$ is injective, where the part $G\to {\bf G}_m$ is $\eta$. 
\end{sbpara}

\begin{sbpara}\label{1.2.1b}

Consider the one-dimensional  $\Q$-vector space $\Q\cdot (2\pi i)^{-1}$ on which $G$ acts via $\eta$.  
By (ii),  $H_b(\Q\cdot (2\pi i)^{-1})$ is a one-dimensional Hodge structure of weight $2$ whose underlying $\Q$-vector space is $\Q \cdot (2\pi i)^{-1}$, and hence it is canonically identified with the Hodge structure $\Q(-1)$. 
The homomorphism $\gr^W_wV_0\otimes \gr^W_wV_0 \to \Q\cdot (2\pi i)^{-w}= (\Q\cdot (2\pi i)^{-1})^{\otimes w}$ defined by $\langle\;,\;\rangle_{0,w}$ is a $G$-homomorphism, and hence induces $\gr^W_wH_b(V_0)\otimes \gr^W_wH_b(V_0) \to \Q(-w)$.

We assume the following (*). 
\medskip

(*) This homomorphism $\gr^W_wH_b(V_0) \otimes \gr^W_wH_b(V_0)\to \Q(-w)$ is a polarization of $\gr^W_wH_b(V_0)$ for any $w\in \Z$. 

\end{sbpara}

\begin{sbpara}\label{neat} 
 From now on, we assume that $K$ satisfies the following condition (*) which we call the neat condition. 

\medskip

(*) If  $g\in \sG({\bf A}^f_\Q)$ and $\gamma\in \sG(\Q) \cap gKg^{-1}\subset \sG({\bf A}^f_\Q)$ and if the action of $\gamma$ on $\gr^W_wV_0$ preserves $\langle\;,\;\rangle_{0,w}$ for any $w\in \Z$, then for any $V\in \Rep(G)$, the subgroup of $\C^\times$ generated by all eigen values of $\gamma: \C\otimes_\Q V\to \C\otimes_\Q V$ is torsion free.

\end{sbpara}

\begin{sbrem}
 If  $L$ is a $\Z$-lattice of $V_0$ and if $n\geq 3$ is an integer, and if the action of any element $g$ of $K$ on ${\bf A}^f_\Q \otimes_\Q V_0$ satisfies the following (i) and (ii), then $K$ satisfies the neat condition (\ref{neat}).

(i) $g\hat L= \hat L$, where $\hat L= \hat \Z\otimes_\Z L$. 

(ii) The automorphism of $L/nL=\hat L/n\hat L$ induced by $g$ is the identify map.

\medskip

Hence a sufficiently small open compact subgroup $K$ of $\sG({\bf A}^f_\Q)$ satisfies the neat condition. 
\end{sbrem}

\begin{sbprop}\label{Aut1} (As is said, we assume \ref{1.2.1}, \ref{1.2.1b}, \ref{neat}.) The automorphism group of a triple $(H, \xi, \la)$ as in \ref{defXCb} is trivial.
\end{sbprop}

\begin{pf} 
The proof goes in the same way as in the corresponding part of Part 7 (the proof of Proposition 2.2.7 of Part I). 
\end{pf}

By \ref{Aut1}, we will identify a triple $(H, \xi, \la)$ as in  \ref{defXCb} with its class in $X(\C)$. 

We will sometimes denote $(H, \xi, \la)\in X(\C)$ also simply as $H\in X(\C)$ omitting $\xi$ and $\la$.
\begin{sbpara}\label{canpol1}
For $(H, \xi, \la)\in X(\C)$, we have a canonical polarization $\langle\;,\;\rangle_w$ on $\gr^W_wH(V_0)$ for each $w\in \Z$ defined as follows. 

Take $\theta$ in \ref{defXCb}. 
Let $g\in \sG({\bf A}^f_\Q)$ be the composition $({\bf A}^f_\Q\otimes \theta)\circ \tilde \la$, where $\tilde \la$ is a representative of $\la$ (\ref{defXCb}), and let $c= \prod_p |\eta(g)_p|_p\in \Q_{>0}$ (here $p$ ranges over all prime numbers, $\eta(g)_p$ denotes the $p$-component of $\eta(g)\in ({\bf A}^f_\Q)^\times$,   and $|\;|_p$ denotes the standard absolute value of $\Q_p$). Then $c$ is independent of the choice of the representative $\tilde \la$ of $\la$ because $|\eta(k)_p|_p=1$ for any $k\in K$ and any $p$. 

Let $\gr^W_wH(V_0)_\Q \overset{\cong}\to \gr^W_wH_b(V_0)_\Q= \gr^W_wV_0$ be the isomorphism defined by $\theta$ and let  $\langle\;,\;\rangle_w'$ be the bilinear form $\gr^W_wH(V_0)_\Q\times \gr^W_wH(V_0)_\Q\to \Q \cdot (2\pi i)^{-w}$ corresponding to $\langle\;,\;\rangle_{0,w}$ via this isomorphism.
By the existence of the isomorphism $\nu$ (\ref{defXCb}), we have either $\langle\;,\;\rangle_w'$ is a polarization for any $w$ or $\langle\;,\;\rangle_w^{'}:=(-1)^w \langle\;,\;\rangle_w'$ is a polarization for any $w$. 
We define the canonical polarization  $\langle\;,\;\rangle_w$ on $\gr^W_wH(V_0)$ to be $c^w\langle\; ,\;\rangle_w'$ in the former case and $(-1)^wc^w\langle\;,\;\rangle_w'$ in the latter case. Then $\langle\;,\;\rangle_w$ is independent of the choice of $\theta$.

\end{sbpara}

\begin{sbprop}\label{paction} (1) The action of $\sG(\Q)$ on $D(G, \sG, H_b)\times (\sG({\bf A}^f_\Q)/K)$ is proper. This action is fixed-point free. 

(2) The canonical surjection $D(G, \sG, H_b) \times (\sG({\bf A}^f_\Q)/K)\to X(\C)$ is a local homeomorphism (for the quotient topology on $X(\C)$). 

\end{sbprop}

\begin{pf} (2) follows from (1).  (1) is reduced to the fact that for any $g\in \sG({\bf A}^f_\Q)$, the action of $\sG(\Q)\cap gKg^{-1}$ on $D(G, \sG, H_b)$ is proper and fixed-point-free. The last thing is proved as follows.

Let $\Phi= ((h(w,r)_{w,r\in \Z}, V_0, W, (\langle\;,\;\rangle_{0,w})_{w\in \Z})$, where $h(w,r)$ is the dimension of $\gr^r \gr^W_w H_b(V_0)_\C$ as a $\C$-vector space, and let $D^{\pm}_{\Phi}$ be the corresponding space $D^{\pm}$ in Part I, Section 2.2. We have an embedding $D(G, \Upsilon)\subset D^{\pm}_{\Phi}$. 
Let $G_{\Phi}$ be the algebraic group over $\Q$ defined as the subgroup of $\text{Aut}_\Q(V_0, W) \times {\bf G}_m$ consisting of all elements $(g,t)$ such that $\langle gx,gy\rangle_{0,w}= t^w\langle x,y\rangle_{0,w}$ for all $w\in \Z$ and all $x,y\in \gr^W_w$. (This group $G_{\Phi}$ was denoted by $G$ in Part I.) We have an embedding $G\overset{\subset}\to G_{\Phi}$.

By  these embeddings $D(G, \sG, H_b)\subset D^{\pm}_{\Phi}$ and $G\subset G_{\Phi}$, we are  reduced to the fact that the action of $G_{\Phi}(\Q)\cap gKg^{-1}$ on $D^{\pm}_{\Phi}$ is proper and fixed-point-free and this follows from Part I, Section 2.2.  
\end{pf}

\begin{sbpara} By \ref{paction} (2), there is a unique structure of a complex analytic manifold on $X_{G,\sG, H_b,K}(\C)$ for which the canonical surjection $D(G, \sG, H_b) \times (\sG({\bf A}^f_\Q)/K) \to X_{G,\sG, H_b,K}(\C)$ is a morphism of  complex analytic manifolds and is locally an isomorphism. 

\end{sbpara}

Note that $\sG_u=\sG\cap G$ and hence $\sG_{\red}$ is the image of $\sG$ in $G_{\red}$.

\begin{sbprop}\label{Xred}

Let $K_{\red}$ be the image of $K$ in $\sG_{\red}({\bf A}^f_\Q)$. Then $(G_{\red}, \sG_{\red}, H_b|_{G_{\red}}, K_{\red})$ also satisfies the neat condition. 

\end{sbprop}

\begin{pf} We are reduced to the fact that $\sG(\Q)\cap gKg^{-1}\to \sG_{\red}(\Q) \cap g_{\red}K_{\red}g_{\red}^{-1}$ is surjective for any $g\in \sG({\bf A}^f_\Q)$ where $g_{\red}$ denotes the image of $g$ in $\sG_{\red}({\bf A}^f_\Q)$. This last fact is reduced to the following well known fact: For a unipotent algebraic group $\sU$ over $\Q$ and for an open subgroup $U$ of $\sU({\bf A}^f_\Q)$, $\sU({\bf A}^f_\Q)= \sU(\Q)U$. \end{pf}

\begin{sbpara}\label{Xred2}
By \ref{Xred},  we have the complex analytic manifold $X_{\red}(\C):= X_{G_{\red}, \sG_{\red}, H_b|_{G_{\red}}, K_{\red}}(\C)$. 
We have a canonical morphism 
$X_{G, \sG, H_b, K}(\C) \to X_{G_{\red}, \sG_{\red}, H_b|_{G_{\red}}, K_{\red}}(\C)$. 
\end{sbpara}

\begin{sbpara}\label{XCrep}

The complex analytic manifold  $X(\C)=X_{G,\sG, H_b, K}(\C)$ represents the functor $\sA\to $ (Sets) which sends $Y\in \sA$ to the set of isomorphism classes of triples $(\sH, \xi, \la)$, where:

$\sH$ is an exact $\otimes$-functors $\Rep(G)\to \Q$MHS$(Y)$ (\ref{Drep}),

$\xi$ is an isomorphism of $\otimes$-functors from $\Rep(\sQ)$ to $\Q$MHS$(Y)$
$$H|_{\sQ}\overset{\cong}\to H_b|_{\sQ},$$
where $H|_{\sQ}$ and $H_b|_{\sQ}$ denote the restrictions of $H$ and $H_b$  to $\Rep(\sQ)\subset \Rep(G)$ ($H_b(V)$ for $V\in \Rep(G)$ is defined to be a constant $\Q$MHS $H_b(V)$ on $Y$), respectively.

$\la$ (called the $K$-level structure) is a global section of the quotient sheaf $\sI/K$, where $\sI(U)$ for an open set $U$ of $Y$ is the set of isomorphisms of $\otimes$-functors from $\Rep(G)$ to $\text{Mod}_{ff}({\bf A}^f_\Q)(U)$
$$(V\mapsto {\bf A}^f_\Q \otimes_\Q V)\overset{\cong}\to (V \mapsto {\bf A}^f_\Q\otimes_\Q \sH_\Q|_U),$$

\medskip
\noindent 
such that there are an isomorphism of $\otimes$-functors from $\Rep(G_{\text{red}})$ to $\R$HS
$$\nu: (V\mapsto \R\otimes_\Q \sH(V)(y)) \overset{\cong}\to (V\mapsto \R \otimes_\Q H_b(V))$$
 at each $y\in Y$ and an isomorphism of $\otimes$-functors from $\Rep(G)$ to $\text{Mod}_{ff}(\Q,Y)$ 
$$\theta: (V \mapsto \sH(V)_\Q)\overset{\cong}\to (V\mapsto V)\quad\text{preserving the weight filtrations}$$
locally on $Y$, which satisfy the following (i)--(iii).

(i) At each $y\in Y$, $\xi$ and $\nu$ induce the same isomorphism of functors $(V\mapsto \R \otimes_\Q \sH(V)(y)) \overset{\cong}\to (V\mapsto \R \otimes_\Q H_b(V))$ from $\Rep(\sQ_{\text{red}})$ to $\R$HS.

(ii) Locally on $Y$, $\xi$ and $\theta$ induce the same isomorphism of functors $(V\mapsto \sH(V)_\Q)\overset{\cong}\to (V\mapsto H_b(V)_\Q)$ from $\Rep(\sQ)$ to $\text{Mod}_{ff}(\Q, Y)$.

(iii) Locally on $Y$, the automorphism  $({\bf A}^f_\Q\otimes \theta) \circ \tilde \la$ of the $\otimes$-functor 
$$\Rep(G) \to \text{Mod}_{ff}({\bf A}^f_\Q)\;;\; V\mapsto {\bf A}^f_\Q \otimes_\Q V$$ 
belongs to $\sG({\bf A}^f_\Q)$.

\end{sbpara}

\begin{sbpara}\label{faith} 
In this paper, we will often use the following fact: If $V_1\in \Rep(G)$ is a faithful representation of $G$, $V_1$ generates $\Rep(G)$ as a $\otimes$-category (\cite{De82} 3.1a).

\end{sbpara}

\begin{sbpara}\label{XCpol}  Taking $Y=X(\C)$ in the above, we have the universal object
$\sH_{X(\C)}: \Rep(G)\to \Q$MHS$(X(\C))$. 

The canonical polarizations on the fibers of $\sH_{X(\C)}(V_0)$ at each point of $X(\C)$ (\ref{canpol1}) give a canonical polarization on $\sH_{X(\C)}$. 

For any $V\in \Rep(G)$ and $w\in \Z$, we have a homomorphism $\gr^W_wV \otimes \gr^W_w V \to 
\Q\cdot (2\pi i)^{-w}$ of representations of $G$ such that $\sH(\gr^W_wV) \otimes \sH(\gr^W_wV)\to \Q \cdot (2\pi i)^{-w}$ is a polarization of $\sH(\gr^W_wV)$. This is true for $V=V_0$ (the canonical polarization) and is true for general $V$ because $\gr^W(V_0)$ and $\Q\cdot (2\pi i)^{-1}$ generate the $\otimes$-category $\Rep(G_{\red})$ (\ref{faith}).   
\end{sbpara}

\begin{sbpara}\label{TX} Concerning the tangent bundle $T_{X(\C)}$ of  $X(\C)$, we have a canonical isomorphism
$$T_{X(\C)}\cong \sH_{X(\C)}(\Lie(\sG))_{\sO}/\fil^0,$$
where $\Lie(\sG)$ is endowed with the adjoint action of $G$ and $\fil$ denote the Hodge filtration. 

This follows from \ref{XCrep} and \ref{Dan} (2). 
\end{sbpara} 

\begin{sbpara}\label{TXhor}  We define the horizontal tangent bundle $T_{X(\C),\hor}$ of $X(\C)$ as 
 $\fil^{-1}/\fil^0$, where $\fil$ is the Hodge filtration of  $\sH_{X(\C)}(\Lie(\sG))_{\sO}$.

It is a subbundle  of the tangent bundle $T_{X(\C)}=\sH_{X(\C)}(\Lie(\sG))_{\sO}/\fil^0$.

\end{sbpara}

\begin{sbpara} In \ref{Xred}, \ref{Xred2}, $T_{X_{\red}}$ and $T_{X_{\red} ,\hor}$ are identified with
$(\gr^W_0 \sH_{X(\C)}(\Lie(\sG)))/\fil^0$ and $\gr^{-1}\gr^W_0\sH_{X(\C)}(\Lie(\sG))$, respectively. 

\end{sbpara}

 \begin{sbpara}\label{varH} Let $Y$ be a complex analytic manifold. We discuss horizontal morphisms from $Y$ to a period domain.

Let $\Q\text{VMHS}(Y)$ be the category of variations of mixed $\Q$-Hodge structure on $Y$. It is a full subcategory of  $\Q$MHS$(Y)$. By definition, an object of $\Q$MHS$(Y)$ belongs to $\Q$VMHS$(Y)$ if and only if it satisfies Griffiths transversality.

A morphism $f: Y\to X(\C)=X_{G,\sG,H_b,K}(\C)$ is said to be {\it horizontal} if the morphism $T_Y \to T_{X(\C)}$ of tangent bundles associated to $f$ factors through the horizontal tangent bundle $T_{X(\C),\hor}\subset T_{X(\C)}$.  A morphism $Y\to X(\C)$ is horizontal if and only if the corresponding object $\sH: \Rep(G) \to \Q$MHS$(Y)$ satisfies that $\sH(V)\in \Q$VMHS$(Y)$ for any $V\in \Rep(G)$.

 \end{sbpara}

\begin{sbpara}\label{hormer} Let $Y$ be a one-dimensional complex analytic manifold (that is, a Riemann surface). 
 Let $\Q\text{VMHS}_{\log}(Y)$ be the full subcategory of $\cup_R \; \Q\text{VMHS}(Y\smallsetminus R)$, where $R$ ranges over all discrete subsets of $Y$, consisting of objects which satisfy the conditions in Part I, Section 1.6 at any point of $R$ when we replace $C$ there by $Y$. 
 
 We say an element of $\cup_R \Mh(Y\smallsetminus R, X(\C))$, where $R$ ranges over all discrete subsets of $Y$, is {\it meromorphic} on $Y$ if the corresponding functor $\Rep(G)\to \cup_R \; \Q\text{VMHS}(Y\smallsetminus R)$ is $\Rep(G)\to \Q$VMHS$_{\log}(Y)$.
 Let 
 $\sM_{\hor}(Y, X(\C))$ be the subset  of $\cup_R \Mh(Y\smallsetminus R, X(\C))$ consisting of all elements which are meromorphic on $Y$.

\end{sbpara}

\begin{sbpara}\label{Hgen} For $Y$ as in \ref{hormer}, we define the subset $\sM_{\hor, \text{gen}}(Y, X(\C))$ of $\sM_{\hor}(Y, X(\C))$. 

 Let 
$f\in \sM_{\hor}(Y, X(\C))$ and let $\sH: \Rep(G) \to \Q\text{VMHS}_{\log}(Y)$ be the corresponding exact $\otimes$-functor. Then by our definition, $f$ belongs to $\sM_{\hor,\text{gen}}(Y, X(\C))$ if and only if there is no algebraic subgroup $G'$ of $G$ such that $\dim(G')<\dim(G)$ and such that  $\sH$ is isomorphic to the composition $\Rep(G)\to \Rep(G') \overset{\sH'}\to \Q\text{VMHS}_{\log}(Y)$ for some exact $\otimes$-functor $\sH': \Rep(G')\to \Q\text{VMHS}_{\log}(Y)$

\end{sbpara}

\subsection{The set $X(F)$ of motives}\label{s1.3}

In this Section \ref{s1.3},
we consider a set  $X(F)=X_{G,\sG, M_b, ,K}(F)$ of $G$-mixed motives (see \ref{XFB} for the definition). We also describe a  special case  $X_{G, \Upsilon,K}(F)$ (see \ref{XFA}), which is simpler,  as a guide. 
We define $X_{G,\Upsilon,K}(F)$ (resp. $X_{G,\sG,M_b,K}(F)$)  imitating the presentation \ref{defXCa} of $X_{G,\Upsilon,K}(\C)$ (resp. \ref{defXCb} of $X_{G,\sG,H_b,K}(\C)$).

\begin{sbpara} As in Part I, we use the formulation of mixed motives due to Jannsen \cite{Ja}. See Part I, Section 1.1 for a summary of points in \cite{Ja} which are important in our study. 

For a number field $F$, we denote the category of mixed motives with $\Q$-coefficients over $F$  by $MM(F)$.

\end{sbpara}
 
\begin{sbpara}\label{XFsetting1} Let $G$ and $w: {\bf G}_m\to G_{\text{red}}$ be as in \ref{Gandw}. Let $V_0\in \Rep(G)$, $\langle\;,\;\rangle_{0,w}:\gr^W_wV_0\times \gr^W_wV_0\to \Q\cdot (2\pi i)^{-w}$ for $w\in \Z$, and $\eta: G\to {\bf G}_m$ be as in \ref{1.2.1}. 
Let $\sG$ be as in \ref{HAB}. 
In what follows, when we discuss $X_{G, \Upsilon, K}(F)$, $\sG$ denotes $G$. Let $K$ be an open subgroup of $\sG({\bf A}^f_\Q)$ satisfying the neat condition (\ref{neat}).

\end{sbpara}

 \begin{sbpara}\label{XFA} Let $\Upsilon$ be as in \ref{hb}. 
 
 For a number field $F\subset \C$, let $X_{G, \Upsilon, K}(F)$ be the set of  isomorphism classes of pairs $(M, \la)$, where
 
\medskip

$M$ is an exact $\otimes$-functor $\Rep(G)\to MM(F)$,

 \medskip

$\lambda$ (called the $K$-level structure) is a mod $K$-class of an isomorphism $\tilde \lambda$ of $\otimes$-functors from $\Rep(G)$ to $\text{Mod}_{ff}({\bf A}^f_\Q)$ $$(V\mapsto {\bf A}^f_\Q\otimes_\Q V) \overset{\cong}\to (V\mapsto  M(V)_{et})\quad \text{preserving the weight filtrations}$$ 
satisfying the condition (i) below, 

\medskip
\noindent
such that there is an isomorphism  of $\otimes$-functors from $\Rep(G)$ to $\text{Mod}_{ff}(\Q)$ 
$$\theta: (V \mapsto M(V)_B)\overset{\cong}\to (V\mapsto V)\quad\text{preserving the weight filtrations}$$
satisfying the conditions (ii) and (iii) below. Here $(-)_B$ denotes the Betti realization with respect to the embedding $F\subset \C$. 

\medskip

(i) There is a homomorphism $\Gal(\bar F/F)\to K$ such that $\tilde \la (k(\sig)x)= \sig \tilde \la(x)$ for any $V\in \Rep(G)$ and for any $x\in {\bf A}^f_\Q\otimes_\Q V$.

(iii) The homomorphism $S_{\C/\R}\to G_{\red,\R}$ induced by the restrictions of $M_H$ and $\theta$ to $\Rep(G_{\red})$ belongs to $\Upsilon$. Here $(-)_H$ denotes the associated $\Q$-mixed Hodge structure with respect to the embedding $F\subset \C$.

(iii) Consider the object $\Q\cdot (2\pi i)^{-1}$ of $\Rep(G)$ on which $G$ acts via $\eta$. Then the isomorphism $M(\Q\cdot (2\pi i)^{-1})_B\cong \Q\cdot (2\pi i)^{-1}$ in $\text{Mod}_{ff}(\Q)$ induced by $\theta$ comes from an isomorphism 

(*) $M(\Q \cdot (2\pi i)^{-1})\cong \Q(-1)$ in $MM(F_0)$, 

\noindent 
and concerning  the morphism $p_w:\gr^W_wM(V_0)\otimes \gr^W_wM(V_0)\to \Q(-w)$ in $MM(F_0)$ induced by the $G$-homomorphism $\gr^W_wV_0\otimes\gr^W_wV_0\to \Q\cdot (2\pi i)^{-w}$ and the above isomorphism (*), we have either $p_w$ is a polarization of $\gr^W_wM(V_0)$ for any  $w\in \Z$ or $(-1)^wp_w$ is a polarization of $\gr^W_wM(V_0)$ for any $w\in \Z$. 

 \end{sbpara}

\begin{sbpara}\label{XFsetting2}

To discuss $X_{G, \sG, M_b, K}(F)$, 
we assume that we are given a number field $F_0\subset \C$ and an exact $\otimes$-functor $$M_b: \Rep(G)\to MM(F_0)$$ 
and that we are given an isomorphism between $\otimes$-functors from $\Rep(G)$ to $\text{Mod}_{ff}(\Q)$
$$(V\mapsto M_b(V)_B)\cong (V\mapsto V) \quad \text{preserving weight filtrations}.$$ 
Here $(\;)_B$ denotes the Betti realization with respect to the inclusion map $F_0\overset{\subset}\to \C$. We regard this isomorphism as an identification. 

We further assume the following (i) and (ii).

Consider the action of $G$ on $\Q \cdot (2\pi i)^{-1}$ via $\eta$. By the condition (ii) in \ref{1.2.1}, we have a unique isomorphism $M_b(\Q \cdot (2\pi i)^{-1})_H\cong \Q(-1)$ of Hodge structures whose underlying isomorphism of $\Q$-vector spaces is the identity map. Here $(\;)_H$ denotes the associated Hodge structure with respect to the inclusion map $F_0\overset{\subset}\to \C$. We assume 

\medskip

(i) There is an isomorphism of motives $M_b(\Q\cdot (2\pi i)^{-1})\cong \Q(-1)$ whose underlying isomorphism of Hodge structures is the above one.

\medskip

Since $\langle\;,\;\rangle_{0,w}$ gives a $G$-homomorphism $\gr^W_wV_0\otimes \gr^W_wV_0\to \Q \cdot (2\pi i)^{-w}=(\Q \cdot (2\pi i)^{-1})^{\otimes w}$, we have a morphism
$\gr^W_wM_b(V_0) \otimes \gr^W_wM_b(V_0) \to \Q(-w)$ of motives over $F_0$. We assume

\medskip
(ii) The last morphism is a polarization on $\gr^W_wM_b(V_0)$ for any $w\in \Z$. 

\medskip 
We further assume

\medskip

(iii) The homomorphism  $\Gal(\bar F/F)\to G({\bf A}^f_\Q)$ defined by the action of $\Gal(\bar F/F)$ on $M_b(V)_{et}= {\bf A}^f_\Q\otimes_\Q V$ ($V\in \Rep(G)$) satisfies $\sig K\sig^{-1}=K$ for any $\sig\in \Gal(\bar F/F)$.

 \end{sbpara}

\begin{sbpara}\label{XFB} Let $\sG$ be as in \ref{HAB}, let $F_0$ and $M_b$ be as in \ref{XFsetting2}, and let $F$  be a finite extension of $F_0$ in $\C$.
Let
$X(F)=X_{G, \sG, M_b,K}(F)$ be the set of isomorphism classes of triples $(M, \xi, \lambda)$, where 

\medskip

$M$ is an exact $\otimes$-functor from $\Rep(G)$ to $MM(F)$, 

\medskip

 $\xi$ is an isomorphism $M|_{\sQ}\overset{\cong}\to M_b|_{\sQ}$, where $M|_{\sQ}$ (resp. $M_b|_{\sQ}$) denotes the restriction of $M$ (resp. $M_b$) to $\Rep(\sQ)$, and 
 \medskip
 
$\la$ (called the $K$-level structure) is a mod $K$ class of an isomorphism of $\otimes$-functors from $\Rep(G)$ to $\text{Mod}_{ff}({\bf A}^f_\Q)$
$$(V\mapsto M_b(V)_{et})\overset{\cong}\to (V\mapsto  M(V)_{et})\quad \text{preserving the weight filtrations}$$
 satisfying the condition (i) below,

\medskip
\noindent
such that there are an isomorphism of $\otimes$-functors from $\Rep(G_{\text{red}})$ to $\R$HS
$$\nu: (V\mapsto \R\otimes_\Q M(V)_H) \overset{\cong}\to (V\mapsto \R \otimes_\Q M_b(V)_H)$$
 and an isomorphism of $\otimes$-functors from $\Rep(G)$ to $\text{Mod}_{ff}(\Q)$ 
$$\theta: (V \mapsto M(V)_B)\overset{\cong}\to (V\mapsto V)\quad\text{preserving the weight filtrations}$$
which satisfy the following  conditions (ii)--(iv) and the condition (iii) in \ref{XFA}.

(i) For any $\sig\in \Gal(\bar F/F)$, there is $k(\sig)\in K$ satisfying $\sig\tilde \lambda (x) = \tilde \lambda (k(\sig)\sig x)$ for any $V\in \Rep(G)$ and $x\in M_b(V)_{et}$ ($\bar F$ denotes the algebraic closure of $F$ in $\C$; here $K$  acts on $M_b(V)_{et}$ because $M_b(V)_{et}={\bf A}^f_\Q\otimes_\Q M_b(V)_B= {\bf A}^f_\Q \otimes_\Q V$),

(ii) $\xi$ and $\nu$ induce the same isomorphism of functors $(V\mapsto \R \otimes_\Q M(V)_H) \overset{\cong}\to (V\mapsto \R \otimes_\Q M_b(V)_H)$ from $\Rep(\sQ_{\text{red}})$ to $\R$HS.

(iii) $\xi$ and $\theta$ induce the same isomorphism of functors $(V\mapsto H(V)_B)\overset{\cong}\to (V\mapsto M_b(V)_B=V)$ from $\Rep(\sQ)$ to $\text{Mod}_{ff}(\Q)$.

(iv) The automorphism  $({\bf A}^f_\Q\otimes \theta) \circ \tilde \la$ of the $\otimes$-functor 
$$\Rep(G) \to \text{Mod}_{ff}({\bf A}^f_\Q)\;;\; V\mapsto {\bf A}^f_\Q \otimes_\Q V$$ 
(which is an element of $G({\bf A}^f_\Q)$ by the theory of Tannakian categories) belongs to $\sG({\bf A}^f_\Q)$.

\end{sbpara}

 \begin{sbpara}\label{FtoC}

 We have a canonical map $X_{G, \Upsilon, K}(F)\to X_{G,\Upsilon,K}(\C)$ (resp. $X(F)=X_{G,\sG,M_b,K}(F)\to X(\C)=X_{G,\sG,M_{b,H},K}(\C)$). Using  \ref{defXCa} (resp. \ref{defXCb}), it is given by $\text{class}(M,\la)\mapsto \text{class}(M_H, \la)$ (resp. 
 $\text{class}(M, \xi, \la)\mapsto \text{class}(M_H, \xi_H, \la)$). 
 
 \end{sbpara}

\begin{sbpara}\label{XFtoXF} (1) Let $(M_1,\xi_1,\la_1)\in X_{G,\sG,M_b, K}(F)$. Then we have a bijection $$X_{G,\sG,M_b, K}(F)\overset{\cong}\to X_{G,\sG, M_1,K'}(F)$$ where to define the last set, we identify $\Rep(G)\to\text{Mod}_{ff}(\Q)\;;\; V\mapsto M_1(V)_B$ with $V\mapsto V$ by fixing $\theta=\theta_1$ for $M_1$, and we define $K'= g_1Kg_1^{-1}$ where $g_1=({\bf A}^f_\Q\otimes \theta_1)\circ \tilde \la_1\in \sG({\bf A}^f_\Q)$.   The bijection is defined by $\text{class}(M, \xi, \la)\mapsto \text{class}(M, \xi_1^{-1}\circ \xi, \la\circ\la_1^{-1})$ where $\la\circ\la_1^{-1}$ denotes the mod $K'$-class of $\tilde \la \circ (\tilde \la_1)^{-1}$.

Via the maps $X_{G,\sG, M_b,K}(F)\to X_{G,\sG, M_{b,H},K}(\C)$ and $X_{G, \sG, M_1,K'}(F)\to X_{G,\sG,M_{1,H},K'}(\C)$ in \ref{FtoC}, this bijection is compatible with the isomorphism $$X_{G,\sG,M_{b,H},K}(\C)= \sG(\Q)\bs (D(G,\sG,M_{b,H})\times (\sG({\bf A}^f_\Q)/K))\overset{\cong}\to $$ $$X_{G,\sG, M_{1,H}, K'}(\C)=\sG(\Q)\bs (D(G, \sG, M_{1,H})\times (\sG({\bf A}^f_\Q)/K'))\;;\;\text{class}(H,g)\mapsto \text{class}(H, gg_1^{-1}).$$

(2) In \ref{XFA}, let $F_0\subset \C$ be a number field and assume that $X_{G, \Upsilon,K}(F_0)$ is not empty. Fix an element $\text{class}(M_1, \la_1)\in X_{G, \Upsilon, K}(F_0)$ and identify $\Rep(G)\to\text{Mod}_{ff}(\Q)\;;\; V\mapsto M_1(V)_B$ with $V\mapsto V$ by fixing $\theta=\theta_1$ for $M_1$. Then we have a bijection $$X_{G,\Upsilon,K}(F)\overset{\cong}\to X_{G,G, M_1, K'}(F)\;;\;\text{class}(M,\la) \mapsto \text{class}(M,\xi, \la\circ\la_1^{-1})$$
where $K'= g_1Kg_1^{-1}$ with $g_1=({\bf A^f_\Q}\otimes   \theta)\circ \tilde \la_1$, $\xi$ is the evident one,  and $\la\circ\la_1^{-1}$ is the mod $K'$-class  of $\tilde \la \circ (\tilde \la_1)^{-1}$. 

Via the maps $X_{G,\Upsilon,K}(F)\to X_{G,\Upsilon,K}(\C)$ and $X_{G,G, M_1,K'}(F)\to X_{G,G, M_{1,H},K'}(\C)$  in \ref{FtoC}, this bijection is compatible with the isomorphism $$X_{G,\Upsilon, K}(\C)= G(\Q)\bs (D(G, \Upsilon) \times (G({\bf A}^f_\Q)/K))\overset{\cong}\to $$ $$X_{G,G,M_{1,H},K'}(\C)= G(\Q)\bs (D(G, G,M_{1,H}) \times (G({\bf A}^f)/K))\;;\; \text{class}(H, g)\mapsto (H, gg_1^{-1}).$$

\end{sbpara}

\begin{sbprop}\label{Aut2}  The automorphism group of a  triple $(M,\xi,\la)$ as in  \ref{XFB}  is trivial.  
\end{sbprop}

This follows from the Hodge version \ref{Aut1}.

By \ref{Aut2}, we will identify a triple $(M, \xi, \la)$ as in \ref{XFB} with its class in $X(F)$. 

We will often denote $(M,\xi,\la)\in X(F)$ simply as $M\in X(F)$.

\begin{sbpara} We expect that the canonical map $X(F)\to X(\C)$ is always injective. 

We expect that if $F'$ is a finite Galois extension of $F$ in $\C$, the map $X(F)\to X(F')$ induces a bijection from $X(F)$ to the $\Gal(F'/F)$-fixed part of $X(F')$.

\end{sbpara}

\begin{sbpara}\label{H1Gal} (1) Note that
$M\in X_{G,\Upsilon,K}(F)$ defines a $K$-conjugacy class of a continuous homomorphism $\Gal(\bar F/F)\to K\;;\;\sigma\mapsto k(\sig)$, that is, an element of $H^1_{\text{cont}}(\Gal(\bar F/F), K)$ where $\Gal(\bar F/F)$ acts on $K$ trivially.

(2) Note that 
$M\in X_{G,\sG, M_b, K}(F)$ defines the class of a continuous $1$-cocycle $\sig\mapsto k(\sig)$ in $ H^1_{\text{cont}}(\Gal(\bar F/F), K)$ where $\sig\in \Gal(\bar F/F)$ acts on $K$ as $k\mapsto \sig \circ k \circ \sig^{-1}$ (here $\sig$ is regarded as an automorphism of $(V\mapsto M_b(V)_{et}= {\bf A}^f_\Q\otimes V)$).

\end{sbpara}

\begin{sbpara}

{\bf Questions.} 

(1) Is the map $X_{G,\Upsilon,K}(F)\to \Hom_{\text{cont}}(\Gal(\bar F/F), K)$ in \ref{H1Gal} (1) always injective?

(2) Is the map $X_{G,\sG,M_b,K}(F)\to H^1_{\text{cont}}(\Gal(\bar F/F), K)$ in \ref{H1Gal} (2) always injective?

\end{sbpara}

\begin{sbpara}\label{canpol2} For $(M, \xi, \la)\in X(F)$, we define the canonical polarization on $\gr^W_wM(V_0)$ for each $w\in \Z$ in the similar way to the Hodge version \ref{canpol1}.

Take $\theta$ in \ref{XFB}. 
Let $g\in \sG({\bf A}^f_\Q)$ be the composition $(\bf A^f_\Q)\otimes \theta)\circ \tilde \la$, and let $c= \prod_p |\eta(g)_p|_p\in \Q_{>0}$. Then the canonical polarization is either $c^wp_w$ or $(-c)^wp_w$ where $p_w$ is as in the condition (iii) in \ref{XFA}.

\end{sbpara}

\begin{sbpara} We define a subset  $X_{\text{gen}}(F)=X_{G,\sG, M_b, K,\text{gen}}(F)$ of $X(F)=X_{G,\sG,M_b,K}(F)$
consisting of {\it generic} elements. Here $(M. \xi, \la)\in X(F)$ is 
 generic means that 
 
 (i) There is no algebraic subgroup $G'$ of $G$ such that $\dim(G')<\dim(G)$ and such that the exact $\otimes$-functor $\Rep(G)\to MM(F)\;;\;V\mapsto M(V)$ comes from an exact $\otimes$-functor $\Rep(G')\to MM(F)$. 
 
 \medskip

Consider the following conditions (ii)--(iv). 
We have (iii) $\Rightarrow $ (ii) $\Rightarrow$ (i), and (iv) $\Rightarrow$ (i). 
In (ii) and (iii), we consider the homomorphism $\Gal(\bar F/F)\to G({\bf A}^f_\Q)$ induced by a representative $\tilde \la$ of $\la$. 
\medskip

(ii) The image of $\Gal(\bar F/F)$ in $G(\Q_p)$ is open for some prime number $p$,

(iii) The image of $\Gal(\bar F/F)$ in $G(\Q_p)$ is open for any prime number $p$.

(iv) There is no algebraic subgroup $G'$ of $G$ such that $\dim(G')<\dim(G)$ and such that the exact $\otimes$-functor $\Rep(G)\to\Q\text{MHS}\;;\;V\mapsto M(V)_H$ comes from an exact $\otimes$-functor $\Rep(G')\to \Q\text{MHS}$.

\medskip

By the philosophy of Mumford-Tate groups, we expect that these (i)--(iv) are equivalent. 
\end{sbpara}

 \begin{sblem}\label{genlem} Let $M\in X(F)=X_{G,\sG, M_b, K}(F)$ and assume that the functor $M:\Rep(G) \to MM(F)$ is fully faithful. Then $M\in X_{\text{gen}}(F)$. 
 
 \end{sblem}

 \begin{pf} Let $G'$ be a linear subgroup of $G$ and assume 
  that $M_b$ is isomorphic to the composition $\Rep(G)\to \Rep(G')\overset{a}\to MM(F)$ for an exact $\otimes$-functor $a: \Rep(G')\to MM(F_0)$. We prove $G'=G$.
  
  Let $\sP$ be the smallest full subcategory of $MM(F)$ which contains the image of $a$ and which is stable under $\otimes$, $\oplus$, the dual, and subquotients. Let $P$ be the Tannakian group of $\sP$ associated to the fiber functor $\sP\to \text{Mod}_{ff}(\Q): M\mapsto M_B$. Then $a$ induces  a homomorphism $P\to G'$, and the composition $P\to G'\to G$ is faithfully flat because the corresponding functor $M_b: \Rep(G)\to\sP$ is fully faithful (\cite{DM}, Proposition 2.20 (a)). Hence $G'=G$. \end{pf}

\begin{sbpara}\label{gen2} 
 We show that any point of $X_{G,\sG, M_b, K}(F)$ comes from a point of $X_{G',\sG', M'_b, K'}(F)$ which is generic for some algebraic subgroup $G'$ of $G$ and for some $\sG'$, $M'_b$, $K'$. 
 
 Let $(M_1,\xi_1,\la_1)\in X_{G,\sG,M_b,K}(F)$. 
 Let $\sC$ be the smallest full subcategory of $MM(F)$ which contains any subquotients of $M_1(V)$ for any $V\in \Rep(G)$.  Let $G'$ be the Tannakian group of $\sC$ with respect to the fiber functor $\sC\to \text{Mod}_{ff}(\Q)\;;\; M\mapsto M_B$. Then $\sC\cong \Rep(G')$. Let $M'_b$ be the composition $\Rep(G')\cong \sC \overset{\subset}\to MM(F)$. Then the fiber functor $\Rep(G')\to \text{Mod}_{ff}(\Q)\;;\; V\mapsto M'_b(V)_B$ coincides with the evident fiber functor.  
 Fix $\theta_1=\theta$ of $M_1$. 
 The composition $\Rep(G)\overset{M_1}\to \sG\cong \Rep(G')\to \text{Mod}_{ff}(\Q)$ is $V\mapsto M_1(V)_B$, and by identifying this fiber functor 
 with the evident fiber functor  $\Rep(G) \to \text{Mod}_{ff}(\Q)$ by using $\theta_1$, we obtain a 
 homomorphism $G'\to G$ which induces the composition $\Rep(G) \overset{M_1}\to \sC\cong \Rep(G')$. This homomorphism $G'\to G$ is 
 a closed immersion by the fact that any object of $\Rep(G')\cong \sC$ is isomorphic to a subquotient of an 
 object which comes from $\Rep(G)$, and by \cite{DM}, Proposition 2.20 (b). 
  Let $\sG'=\sG\cap G'$,
  and let $K'= 
 g_1 K g_1^{-1}\cap \sG'({\bf A}^f_\Q)$ where $g_1= ({\bf A}^f_\Q \otimes \theta_1)\circ \tilde \la\in \sG({\bf A}^f_\Q)$. 
 Consider the set $X'(F)= X_{G', \sG', M'_b, K'}(F)$. 
We have a map $$X'(F)\to X(F):=X_{G, \sG, M_b, K}(F)\;;\; (M. \xi, \la)\mapsto (M, \xi_1\circ \xi, \la \circ \la_1)$$ which sends $M'_b\in X'(F)$ to $(M_1, \xi_1, \la_1)\in X(F)$. By \ref{genlem}, $M'_b$ is generic in $X'(F)$.
 \end{sbpara}

\subsection{Examples}

We give some examples. In \ref{HEx3}, we explain that $X(\C)$ and $X(F)$ in Part I can be regarded as special cases of $X(\C)$ and $X(F)$ of this Part II, respectively.  In \ref{HEx4}, by using \cite{KNU2},  
we explain that a higher Albanese manifold (\cite{Ha}) is an example of $X(\C)$, and then in \ref{MEx4}, we give an arithmetic version of a higher Albanese manifold as an example of $X(F)$.

  \begin{sbpara}\label{HEx2} Shimura variety.
  
  Assume $G$ is reductive. If $\Upsilon$ satisfies the conditions of Deligne \cite{De71} to define a Shimura variety, $X_{G,\Upsilon,K}(\C)$ is the space of $\C$-points of the Shimura variety associated to $(G, \Upsilon)$  of level $K$ (\cite{De71}).  

If $\Upsilon$ defines a Shimura variety of Hodge type (\cite{Mi} Section 7), 
we have a canonical map $X'(F)\to X_{G,\Upsilon,K}(F)$, 
where $X'$ is the  Shimura variety associated to $(G,\Upsilon)$  of level $K$ (\cite{De71})
and $X'(F)$ is the set of $F$-points of $X'$. This map is defined by sending a class of an abelian variety $A$ to  the class of  the motive $H^1(A)$ ($X'$ is a moduli space of abelian varieties). We expect that this map $X'(F)\to X_{G,\Upsilon,K}(F)$ is bijective.

\end{sbpara}

 \begin{sbpara}\label{HEx1} The period domain associated to a mixed Hodge structure (this is the mixed Hodge version of the Mumford-Tate domain associated to pure Hodge structures in \cite{GGK}).

 Let $H_0$ be a $\Q$-mixed Hodge structure and assume that a polarization $p_w: \gr^W_wH_0\otimes \gr^W_wH_0 \to \Q(-w)$ of $\gr^W_wH_0$ is given for each $w\in \Z$.

Let $\sC$ be the smallest full subcategory of $\Q$MHS which contains $H_0$ and $\Q(-1)$ and is stable under $\otimes$, $\oplus$, the dual, and subquotients. Let $G$ be the Tannakian group of $\sC$ associated to the fiber functor $H\mapsto H_\Q\;;\;\Q$MHS $\to$ $\text{Mod}_{ff}(\Q)$. It is the linear algebraic group over $\Q$ defined as the automorphism group of this fiber functor. Then this fiber functor induces an equivalence of categories $\sC\overset{\simeq}\to \text{Rep}(G)$. Let $H_b$ be the composite functor $\Rep(G) \overset{\simeq}\to \sC \overset{\subset}\to \Q$MHS.

  Let $\sC_{\text{red}}$ be the full subcategory of $\sC$ consisting of objects which belong to $\Q$HS. Then the Tannakian group of $\sC_{\text{red}}$ associated to the fiber functor $H\mapsto H_Q$ is identified with $G_{\text{red}}$. The weight filtrations of objects of $\sC_{\text{red}}$ which uniquely split define a homomorphism $w:{\bf G}_m\to G_{\text{red}}$ and this homomorphism $w$ satisfies the condition in \ref{Gandw}.

 Let $V_0= H_{0,\Q}$ with the action of $G$, and for $w\in \Z$, let $\langle\;,\;\rangle_{0,w}: \gr^W_wV_0\times \gr^W_wV_0\to \Q\cdot (2\pi i)^{-w}$ be the pairing induced by the polarization $p_w$. Let $\eta: G\to {\bf G}_m$ be the homomorphism defined by the action of $G$ on $\Q(-1)_\Q=\Q\cdot (2\pi i)^{-1}$. Then the conditions in \ref{1.2.1} and \ref{1.2.1b} are satisfied.

 Let  $G_{\Phi}$ be the algebraic subgroup of $\text{Aut}_\Q(H_{0,\Q},W)\times {\bf G}_m$ consisting of all $(g,t)$ such that $\langle gx, gy\rangle_w = t^w \langle x, y\rangle_w$ for all $x,y$ and $w$. Let $\Upsilon$ (resp. $\Upsilon_{\Phi}$) be the $G(\R)$ (resp. $G_{\Phi}(\R)$)-conjugacy class of the homomorphism $S_{\C/\R}\to G_{\red,\R}$ (resp. $S_{\C/\R}\to G_{\Phi,\red, \R}$) associated to $\gr^W H_0$. Then 
 the representation $G$ on $H_{0,\Q}$ induces an injective homomorphism $G\to G_{\Phi}$ and induces an injective morphism $D(G, \Upsilon) \to D(G_{\Phi}, \Upsilon)$.

For an algebraic normal subgroup $\sG$ of $G$ and for an open compact subgroup $K$ of $\sG({\bf A}^f_\Q)$ satisfying the neat condition (\ref{neat}), we have the complex analytic manifold $X_{G,\sG, H_b, K}(\C)$.

\end{sbpara}

 \begin{sbpara}\label{MEx1} Motive version of \ref{HEx1}.

 Let $F_0\subset \C$ be a number field, let $M_0$ be a $\Q$-mixed motive over $F_0$,  and assume that a polarization $p_w: \gr^W_wM_0\otimes \gr^W_wM_0 \to \Q(-w)$ of $\gr^W_wM_0$ is given for each $w\in \Z$.

Let $\sC$ be the smallest full subcategory of $MM(F_0)$ which contains $M_0$ and $\Q(-1)$ and is stable under $\otimes$, $\oplus$, the dual, and subquotients. Let $G$ be the Tannakian group of $\sC$ associated to the fiber functor $M\mapsto M_B\;;\;MM(F_0) \to\text{Mod}_{ff}(\Q)$. 
Then this fiber functor induces an equivalence of categories $\sC\overset{\simeq}\to \text{Rep}(G)$. Let $M_b$ be the composite functor $\Rep(G) \overset{\simeq}\to \sC \overset{\subset}\to \Q$MHS.

  Let $\sC_{\text{red}}$ be the full subcategory of $\sC$ consisting of objects which are direct sums of pure objects. Then the Tannakian group of $\sC_{\text{red}}$ associated to the fiber functor $M\mapsto M_B$ is identified with $G_{\text{red}}$. The weight filtrations of objects of $\sC_{\text{red}}$ which uniquely split define a homomorphism $w:{\bf G}_m\to G_{\text{red}}$ and this homomorphism $w$ satisfies the condition in \ref{Gandw}.

 Let $V_0= M_{0,B}$ with the action of $G$, and for $w\in \Z$, let $\langle\;,\;\rangle_{0,w}: \gr^W_wV_0\times \gr^W_wV_0\to \Q\cdot (2\pi i)^{-w}$ be the pairing induced by the polarization $p_w$. Let $\eta: G\to {\bf G}_m$ be the homomorphism defined by the action of $G$ on $\Q(-1)_B=\Q\cdot (2\pi i)^{-1}$. Then the conditions in \ref{1.2.1}, \ref{1.2.1b}, and \ref{XFsetting2} are satisfied. 
 
 For an algebraic normal subgroup $\sG$ of $G$, for an open compact subgroup $K$ of $\sG({\bf A}^f_\Q)$ satisfying the neat condition (\ref{neat}), and for a finite extension $F$ of $F_0$ in $\C$, we have the set $X_{G,\sG, M_b, K}(F)$.

\end{sbpara}

 \begin{sbpara}\label{HEx4} Mixed Hodge structures with fixed pure graded quotients.
 
 Let the notation be as in \ref{HEx1} except that we denote $G$ in \ref{HEx1} by $G_0$ here. Assume $H_0$ is endowed with a structure of $\Z$-Hodge structure, that is, a $\Z$-structure $H_{0,\Z}$ of $H_{0,\Q}$ is given. 
 We show that as an example of $X(\C)$, we have the classifying space of $\Z$MHS $H$ such that $\gr^WH=\gr^WH_0$.

 Let $\sQ=G_{0,\red}$, let $G$ be the inverse image of $G_{0,\red}\subset G_{\Phi,\red}$ in $G_{\Phi}$ (hence $G_0\subset G$), let $\sG= G_{0,u}= G_{\Phi,u}$, and let $K= \{g\in \text{Aut}_{\hat \Z}(\hat \Z\otimes_{\Z} H_{0,\Z}, W)\;|\; \gr^W(g)=1\}$. We denote the composition $\Rep(G)\to \Rep(G_0) \to \Q MHS$, where the last arrow is $H_b$ in \ref{HEx1}, by the same letter $H_b$.

  As a subset of $D(G_{\Phi}. \Upsilon_{\Phi})$, $D(G, \sG, H_b)$ is identified with the set of decreasing filtrations $\fil$ on $H_{0,\C}$ such that $\gr^W\fil=\gr^WH_{0,\C}$. This induces a bijection between 
 $X(G, \sG, H_b, K)$ and the set of all isomorphism classes of a $\Z$MHS such that 
 $\gr^W H= \gr^W H_0$.

  \end{sbpara}

 \begin{sbpara}\label{MEx4} Mixed motives with fixed pure graded quotients.

 Let the notation be as in \ref{MEx1} except that we denote $G$ in \ref{MEx1} by $G_0$ here. Assume $M_0$ is endowed with a structure of $\Z$-mixed motive, that is, a $\Z$-structure $M_{0,B,\Z}$ of $M_{0,B}$ such that $\hat \Z \otimes_{\Z} M_{0,B, \Z}$ is stable in $M_{\text{et}}= \hat \Z \otimes_\Z M_{0,B}$ under the action of $\Gal(\bar F_0/F_0)$.  
 We show that as an example of $X(F)$, we have the classifying space of mixed motives $M$ over $F$ with $\Z$-coefficients  such that $\gr^WH=\gr^WH_0$.

 Let $\sQ=G_{0,\red}$, let $G$ be the inverse image of $G_{0,\red}\subset G_{\Phi,\red}$ in $G_{\Phi}$ (hence $G_0\subset G$), let $\sG= G_{0,u}= G_{\Phi,u}$, and let $K= \{g\in \text{Aut}_{\hat \Z}(\hat \Z\otimes_{\Z} H_{0,\Z}, W)\;|\; \gr^W(g)=1\}$. We denote the composition $\Rep(G)\to \Rep(G_0) \to MM(F_0)$, where the last arrow is $M_b$ in \ref{MEx1}, by the same letter $M_b$.

 Then $X_{G, \sG, M_b, K}(F)$ is identified with the set $E$ of isomorphism classes of mixed motives over $F$ with $\Z$-coefficients such that $\gr^WM=\gr^W M_0$. In fact, the map $X_{G, \sG, M_b, K}(F)\to E$ is given by $M\mapsto M(M_{0,B})$. The converse map $E\to X_{G,\sG,M_b, K}(F)
  \;;\;\newline M\mapsto (\tilde M, \xi, \la)$ is obtained as follows. Let $\sC'$ be the smallest full subcategory of $MM(F)$ which contains $M$ and $\Q(-1)$ and which is stable under $\oplus$, $\oplus$, the dual, and subquotients. Let $G'$ be the Tannakian group of $\sC'$ with respect to the fiber functor $M'\mapsto M'_B$. Then $\sC'_{\red}=\sC_{\red}$ and  $(G')_{\red}=G_{0,\red}$. Hence we have a canonical injective homomorphism $G'\to G$. We define $\tilde M$ as the composite functor $\Rep(G)\to \Rep(G')\cong \sC'\subset MM(F)$. The composition $\Rep(\sQ)\to \Rep(G)\to MM(F)$ induced by $\tilde M$ coincides with $\Rep(\sQ) \overset{\cong}\to \sC_{\red} \overset{\subset}\to  MM(F)$ and hence coincides with the functor induced by $M_b$. This gives $\xi$. For $V\in \Rep(G)$, $\tilde M(V)_B$ is identified with $V$ and hence $\tilde M(V)_{et}$ is identified with ${\bf A}^f_\Q\otimes_\Q V$. This identification defines $\la$.

  \end{sbpara}

\begin{sbpara}\label{HEx3} Relation to Part I. 

We show that the period domain $X(\C)$ and the set of motives $X(F)$ in Part I is regarded as an example of $X_{G,\Upsilon,K}(\C)$ and $X_{G,\Upsilon,K}(F)$ of this Part II, respectively.  

As in Section 2.2 of Part I, let $\Phi=((h(w,r))_{w,r\in \Z}, H_{0,\Q}. W, (\langle\;,\;\rangle_{0,w})_{w\in \Z})$, where $h(w,r)\in \Z_{\geq 0}$ satisfying $h(w,r)=0$ for almost all $(w,r)$ and $h(w,r)=h(w,w-r)$ for all $w,r$, $H_{0,\Q}$ is a $\Q$-vector space of dimension $\sum_{w,r} h(w,r)$, $W$ is an increasing filtration on $H_{0,\Q}$ such that $\dim_\Q(\gr^W_w)= \sum_r h(w,r)$ for all $w$, and $\langle\;,\;\rangle_{0,w}$ is a non-degenerate $\Q$-bilinear form $\gr^W_wH_{0,\Q}\times \gr^W_wH_{0,\Q}\to \Q\cdot (2\pi i)^{-w}$ for each $w\in \Z$ which is symmetric if $w$ is even and is anti-symmetric if $w$ is odd. 

As in Part I, Section 2.2, define the linear algebraic group $G$ over $\Q$ as
$G=\{(g,t)\in \text{Aut}(H_{0,\Q}, W)\times {\bf G}_m\;|\;\langle gx,gy\rangle_{0,w}=t^w \langle x,y\rangle_{0,w} \forall x,y\in \gr^W_wH_{0,\Q}\:\}$. Let $K$ be an open subgroup of $G({\bf A}^f_\Q)$ satisfying the neat condition. 
As in Part I, Section 2.2, let  $D^{\pm}$ be the set of decreasing filtrations $\fil$ on $H_{0,\C}$ such that either $(\gr^W_wH_{0,\Q}, \gr^W_w\fil, \langle\;,\;\rangle_{0,w})$ is a polarized Hodge structure of weight $w$ for any $w\in \Z$ or 
$(\gr^W_wH_{0,\Q}, \gr^W_w\fil, (-1)^w\langle\;,\;\rangle_{0,w})$ is a polarized Hodge structure of weight $w$ for any $w$. Let $X_{\Phi,K}(\C)$ (resp. $X_{\Phi,K}(F)$) be $X(\C)$ (resp. $X(F)$) in Part I, Section 2.2.

Since $G(\R)G_u(\C)$ acts transitively on $D^{\pm}$ (Part I, Section 2.2), $\gr^W_w\fil$ ($w\in \Z$) for $\fil\in D^{\pm}$ give a $G_{\red}(\R)$-conjugacy class $\Upsilon$ of $S_{\C/\R}\to G_{\red\R}$. 

We show $$D^{\pm}=D(G, \Upsilon),\quad X_{\Phi,K}(\C)=X_{G,\Upsilon,K}(\C),\quad X_{\Phi,K}(F)=X_{G,\Upsilon,K}(F).$$

We have a canonical map $D(G, \Upsilon) \to D^{\pm}\;;\;H\mapsto H(H_{0,\Q})$. Since $G(\R)G_u(\C)$ acts on $D(G,\Upsilon)$ and $D^{\pm}$ transitively and since this map is compatible with these actions, this map is surjective. The injectivity follows from the fact that $G\to \text{Aut}(H_{0,\Q})\times {\bf G}_m$ is injective and from \ref{faith}.

Hence the space $X_{\Phi,K}(\C)= G(\Q)\bs (D^{\pm}\times (G({\bf A}^f_\Q)/K)$  is identified with  $X_{G, \Upsilon, K}(\C)=G(\Q)\bs (D(G, \Upsilon) \times (G({\bf A}^f_\Q)/K))$.

We have a canonical map
$$X_{G,\Upsilon,K}(F)\to X_{\Phi,K}(F)\;;\;\text{class}(M,\la)\mapsto \text{calss}(M', \la_1, \theta_1, \la_2, \theta_2)$$
where 
$M'=M(H_{0,\Q})$. $\la_1= \tilde \la_{H_{0,\Q}}$. $\theta_1=\theta_{H_{0,\Q}}$. $\la_2=\theta\circ \tilde \la$. $\theta_2=1$.
The converse map $$X_{\Phi,K}(F) \to X_{G,\Upsilon, K}(F)\;;\text{class}(M, \la_1, \la_2, \theta_1,\theta_2)\mapsto (\tilde M, \la)$$
is defined as follows. Let $\sC$ be the smallest full subcategory of $MM(F)$ which contains $M$ and $\Q(-1)$ and which is stable under $\otimes$, $\oplus$, the dual, and subquotients, and let $G'$ be the Tannakian group of $\sC$ with respect to the fiber functor $(-)_B$. 
Then $\sC\cong \Rep(G')$. We have a canonical  injective homomorphism $G'\to G$. Let $\tilde M$ be the composite functor  $\Rep(G)\to \Rep(G')\cong \sC\subset MM(F)$. Then for $V\in \Rep(G)$, $\tilde M(V)_B$ is identified with $V$. We define $\la$ by $\tilde \la= (\theta_1\la_1, \theta_2\la_2)\in G({\bf A}^f_\Q)$.

\end{sbpara}

\begin{sbpara}\label{halbH} Higher Albanese manifolds.

Here by using \cite{KNU2},
we show that the higher Albanese manifold of Hain \cite{Ha} is regarded as an example of $X(\C)$. 
We also correct a mistake in \cite{KNU2}. 

Let $Z$ be a connected smooth quasi-projective algebraic variety over $\C$. 
Fix $b\in Z$. Let $J$ be the augmentation 
ideal $\text{Ker}(\Q[\pi_1(Z,b)]\to \Q)$ of the group ring $\Q[\pi_1(Z, b)]$. Fix $n\geq 0$, and let $\Gamma=\Gamma_n$ be the image of $\pi_1(X,b)\to \Q[\pi_1(Z,b)]/J^{n+1}$. Then $\Gamma$ is a finitely generated torsion-free nilpotent group.

Let $\sG$ be the unipotent algebraic group over $\Q$ 
whose Lie algebra is defined as follows. Let $I$ be the augmentation 
ideal $\text{Ker}(\Q[\Gamma]\to \Q)$ of $\Q[\Gamma]$.
Then 
$\Lie(\sG)$ is the $\Q$-subspace of 
$\Q[\Gamma]^{\wedge}:=\varprojlim_r \Q[\Gamma]/I^r$ generated by all 
$\log(\gamma)$ ($\gamma \in \Gamma$).

We have 
$$\Lie(\sG)= \{h \in \Q[\Gamma]^{\wedge}\;|\; \Delta(h)= h \otimes 1 + 
1 \otimes h\},$$ 
$$\sG (R)= \{g\in (R[\Gamma]^{\wedge})^\times\;|\; \Delta(g)= g \otimes g\}$$ 
for any commutative ring $R$ over $\Q$, where $\Delta: 
R[\Gamma]^{\wedge}\to R[\Gamma\times \Gamma]^{\wedge}$ is the ring 
homomorphism induced by the ring homomorphism $R[\Gamma]\to 
R[\Gamma\times \Gamma]\;;\;\gamma \mapsto \gamma\otimes \gamma$ 
($\gamma \in \Gamma$). The Lie product of $\Lie(\sG)$ 
is defined by $[x,y]= xy-yx$. 
We have $\Gamma \subset \sG(\Q)$.

 By the work \cite{HZ} of Hain-Zucker, 
 $\Lie(\sG)$ is regarded as a polarizable mixed $\Q$-Hodge structure. The Lie product $\Lie(\sG)\otimes \Lie(\sG)\to \Lie(\sG)$ is a homomorphism of mixed Hodge structures.

The $n$-th higher Albanese manifold $\Alb_{Z,n}(\C)$ of $Z$ is as follows. 
Let $\fil^0\sG(\C)$ be the algebraic subgroup of $\sG(\C)$ over $\C$ corresponding to the Lie subalgebra $\fil^0\Lie(\sG)_\C$ ($\fil^0$ here denotes the Hodge filtration) of $\Lie(\sG)_\C$. Then
$$\Alb_{Z, n}(\C):= \Gamma\bs \sG(\C)/\fil^0\sG(\C).$$

This is a complex analytic manifold but usually it is not an algebraic variety.

In \cite{KNU2}, we took as $\Gamma$ any quotient group  of $\pi_1(Z, b)$ which is  nilpotent and torsion-free  (we did not assume $\Gamma=\Gamma_n$). This was a mistake because for such general $\Gamma$, the Lie algebra $\Lie(\sG)$  need not have a mixed Hodge structure. The authors of \cite{KNU2} correct this mistake by adding in 5.1.1 of \cite{KNU2} the assumption that $\Lie(\sG)$ has a mixed Hodge structure which is a quotient of the case $\Gamma=\Gamma_n$ for some $n$. With this correction, all arguments and results of \cite{KNU2} work. (This correction will be written also in the joint paper \cite{KNU} Part V.)

 By the work \cite{KNU2}, the higher Albanese manifold 
   $\Alb_{Z,n}(\C)$ is interpreted as an example of $X(\C)$ as follows.

       Let $\sC$ be a full subcategory of $\Q$MHS which contains the mixed Hodge structures $\Lie(\sG)$ and $\Q(-1)$, which is stable under $\otimes$, $\oplus$, the dual,  and subquotients, and which is of finite type as a Tannakian category.

      We define categories  $\sC_Z$ and $\sC'_Z$.

   Let $\sC_Z$ be the category of variations of mixed $\Q$-Hodge structure $\sH$ on $Z$ satisfying the following conditions (i)--(iv).
   
   \medskip
   
   (i) All graded quotients $\gr^W_w\sH$ for the weight filtration are constant mixed Hodge strutures.
   
   (ii) The monodromy actions of $\pi_1(Z, b)$ on the fiber $\sH_{\Q,b}$ at $b$ factors through the projection $\pi_1(Z,b)\to\Gamma$.  
   
   (iii) The fiber $\sH(b)$ of $\sH$ at $b$ belongs to $\sC$. 

   (iv) $\sH$ is  good at the boundary of $Z$ in the sense of \cite{HZ}.

Let  $\sC'_Z$ be the category of $h\in \sC$
    which is endowed with a morphism  $\Lie(\sG)\otimes h \to h$ in $\Q$MHS which is a Lie action of $\Lie(\sG)$ on $h$.

   By \cite{HZ}, we have an  equivalence of categories 
   
   \medskip
   
   $\sC_Z\overset{\simeq}\to \sC'_Z$ 
      
   \medskip
   \noindent
which sends $\sH\in \sC_Z$ to the fiber $\sH(b)$ of $\sH$ at $b$ with the action of $\Lie(\sG)$ induced by the monodromy action of $\Gamma$ on $\sH_{\Q,b}$.

       Let $\sQ_{\sC}$ be the Tannakian group of $\sC$ associated to the fiber
   functor $\sC\to \text{Mod}_{ff}(\Q)\;;\; H\mapsto H_\Q$.  By the mixed $\Q$-Hodge structure on $\Lie(\sG)$, we have an action of $\sQ_{\sC}$ on the Lie algebra $\Lie(\sG)$. This induces an action of $\sQ_{\sC}$ on the algebraic group $\sG$. Let $G_{\sC}$ be the semi-direct product of $\sG$ and $\sQ_{\sC}$ in which $\sG$ is a normal subgroup of $G_{\sC}$ and in which the inner-automorphism action of $\sQ_{\sC}$ on $\sG$ is the action which we just defined. Then we have  equivalences of categories
   $$\sC\simeq \Rep(\sQ_{\sC}), \quad \sC'_Z\simeq \Rep(G_{\sC}),$$
where  the second equivalence is  because a linear representation of $\sG$ and a representation of the Lie algebra $\Lie(\sG)$ are equivalent. 
   The functor $\Rep(\sQ_{\sC}) \to \Rep(G_{\sC})$ induced by $G_{\sC}\to \sQ_{\sC}$ corresponds to the functor $\sC\to \sC'_Z$ to give the trivial action of $\Lie(\sG)$. The functor $\Rep(G_{\sC})\to \Rep(\sQ_{\sC})$ induced by the inclusion map $\sQ_{\sC}\to G_{\sC}$ corresponds to the functor $\sC'_Z\to \sC$ to forget the action of $\Lie(\sG)$. 
  Let $H_b:\Rep(G_{\sC})\to \Q$MHS be the composition $\Rep(G_{\sC})\to \Rep(\sQ_{\sC})\cong \sC\overset{\subset}\to \Q$MHS. 
  Then the action of $\sG(\C)$ on $D=D(G_{\sC},\sG, H_b)$ is transitive and $\fil^0\sG(\C)$ coincides with the isotropy group of $\text{class}(H_b)\in D$ in $\sG(\C)$. Hence we have $D\cong \sG(\C)/\fil^0\sG(\C)$. 

Let 
$K$ is the profinite completion of $\Gamma$, which is naturally regarded as an open compact subgroup of $\sG({\bf A}^f_\Q)$. 
Then 
$\sG(\Q)/\Gamma\overset{\cong}\to \sG({\bf A}^f_\Q)/K$. Hence 
$$X_{G_{\sC},\sG, H_b, K}(\C)=\sG(\Q)\bs (D\times (\sG({\bf A}^f)/K))=\sG(\Q)\bs(D\times (\sG(\Q)/\Gamma))= \Gamma\bs D =\Alb_{Z,n}(\C).$$
By this, $\Alb_{Z,n}(\C)$ is identified with $X_{G_{\sC},\sG, H_b,K}(\C)$.

   The higher Albanese map $Z(\C) \to \Alb_{Z, n}(\C)$ of Hain is interpreted as follows.
   Let $X_{\sC}(\C)$ be the set of isomorphism classes of triples $(H, \xi, \la)$, where $H$ is an exact $\otimes$-functor $\sC_Z\to \Q$MHS, $\xi$ is an isomorphism of $\otimes$-functors from $\Q$MHS to $\Q$MHS
   $$(h \mapsto H(h_Z))\overset{\cong}\to (h\mapsto h),$$
   where $h_Z$ denotes the constant varitation of mixed $\Q$-Hodge structure on $Z$ associated to $h$, and $\la$ is a mod $\Gamma$-class of an isomorphism of $\otimes$-functors from $\sC_Z$ to $\text{Mod}_{ff}(\Q)$
   $$(\sH\mapsto \sH(b)_\Q)\overset{\cong}\to (\sH\mapsto H(\sH)_\Q) \quad \text{preserving the weight filtrations}$$
   such that the mod $K$-class of the isomorphism of $\otimes$-functors from $\Q$MHS to $\text{Mod}_{ff}(\Q)$
   $$(h\mapsto h_\Q)\overset{\cong}\to (h\mapsto H(h)_\Q)$$
   induced by $\xi^{-1}$ coincides with that induced by $\la$.

The commutative diagram of categories
$$\begin{matrix} \sC & = & \sC & \simeq & \Rep(\sQ)\\
\downarrow &&\downarrow&&\downarrow \\
\sC_Z &\overset{\simeq}\to & \sC'_Z & \simeq & \Rep(G_{\sC}),
\end{matrix}$$
where the vertical arrows are the pullback functors, induces a bijection
$$X_{G_{\sC}, \sG, H_b,K}(\C) \overset{\cong}\to X_{\sC}(\C).$$

The higher Albanese map $Z(\C) \to \Alb_{Z, n}(\C)$ of Hain is interpreted as the map
  $Z(\C) \to X_{\sC}(\C)$, which sends $s\in Z$ to roughly speaking, the functor to take the fiber at $s$.  Precisely speaking, it sends $s$ to  the class of $(H, \xi, \la)$, where $H$ is 
   the 
  exact  $\otimes$-functor $\sC_Z\to \Q\text{MHS}\;;\;\sH\mapsto \sH(s)$, $\xi$ is the evident identification $h(s)=h$ ($h\in \Q$MHS), and $\la$ is the mod $\Gamma$ class of $\sH(s)_\Q \cong \sH(b)_\Q$ of the local system $\sH_\Q$.

\end{sbpara}

\begin{sbpara}\label{halbM} Motive version of \ref{halbH}.

We consider the arithmetic version  of the higher Albanese manifold, which is obtained as  $X(F)$ and relate it to the Selmer variety of Kim \cite{Kim} which is also an arithmetic version of a higher Albanese manifold.

Let $F_0\subset \C$ be a number field and let $Z$ be a geometrically connected smooth quasi-projective algebraic variety over $F_0$. Assume we are given $b\in Z(F_0)$. Let $n\geq 0$, and define $\sG$ and $K$
 as in \ref{halbH} by taking  $Z(\C)$  as $Z$ in \ref{halbH}.

By Deligne-Goncharov \cite{DG}, 3.12, we  have a mixed motive 
 $\Lie(\sG)_{\text{mot}}$ with $\Q$-coefficients over $F_0$ whose Hodge realization is  the mixed Hodge structure $\Lie(\sG)$ in \ref{halbH}. We have a morphism $\Lie(\sG)_{\text{mot}}\otimes \Lie(\sG)_{\text{mot}}\to \Lie(\sG)_{\text{mot}}$ of $MM(F_0)$ which induces the Lie product $\Lie(\sG)\otimes \Lie(\sG)\to \Lie(\sG)$. 
 
 We define a set $\Alb_{Z, n}(F)$ for a finite extension $F$ of $F_0$ in $\C$. 
 
 Let $\sC$ be a full subcategory of $MM(F_0)$, which contains $\Lie(\sG)_{\text{mot}}$ and $\Q(-1)$ and which is stable under $\otimes$, $\oplus$, the dual, and subquotients, and which is of finite type as a Tannakian category.

Let $\sQ_{\sC}$ be the Tannakian group of $\sC$ with respect to the fiber functor $\sC\to \text{Mod}_{ff}(\Q)\;;\; M\mapsto M_B$. Then $\sC\simeq \Rep(\sQ_{\sC})$.
Since $\Lie(\sG)_{\text{mot}}\in \sC$, $\sQ$ acts on $\Lie(\sG)_{\text{mot},B}=\Lie(\sG)$ preserving the Lie product.  Hence $\sQ$ acts on the algebraic group $\sG$. Let $G_{\sC}$ be the semi-direct product of $\sG$ and $\sQ$ in which $\sG$ is a normal subgroup and in which the inner-automorphism action of $\sQ$ on $\sG$ is the one just defined. 
Let $M_b:\Rep(G)\to MM(F_0)$ be the composition $\Rep(G)\to \Rep(\sQ) \cong \sC \overset{\subset}\to MM(F_0)$. 
 Define
 $$\Alb_{Z, n,\sC}(F):= X_{G_{\sC}, \sG, M_b, K}(F).$$

 We have a canonical map $\Alb_{Z,n,\sC}(F)\to \Alb_{Z,n}(\C)$ defined as follows. 
Let $\sC'$ be the Tannakian subcategory of $\Q$MHS generated by Hodge realizations of objects of $\sC$. Let $G_{\sC'}\supset \sG$ be the $G_{\sC'}$ in \ref{halbH} associated to $\sC'$. Then we have a canonical homomorphism $G_{\sC'}\to G_{\sC}$ which induces the identity map of $\sG$.
 Denote by $H_b$ the functor $\Rep(G_{\sC})\to  \Q$MHS induced by $M_b$, and denote by
$H'_b$  the canonical functor $\Rep(G_{\sC'})\to \Q$MHS. 
Then we can show that the canonical map
$X_{G_{\sC'}, \sG, H'_b, K}(\C)\to X_{G_{\sC},\sG, H_b, K}(\C)$ is an isomorphism. Hence we have a map
$\Alb_{Z,n,\sC}(F)=X_{G_{\sC},\sG, M_b, K}(F) \to X_{G_{\sC},\sG, H_b, K}(\C)  \cong X_{G_{\sC'}, \sG, H_b, K}(\C)= \Alb_{Z, n}(\C)$. 

We have a canonical map from $\Alb_{Z, n,\sC}(F)$ to the Selmer vareity  of M. Kim by  \ref{H1Gal} (2).

The author can not show that $\Alb_{Z, n,\sC}(F)$ is independent of $\sC$. We define  $\Alb_{Z,n}(F)$ as the inverse limit of $\Alb_{Z,n,\sC}(F)$ for all $\sC$.

Concerning the higher Albanese map from $Z(F)$, we may define the motive version of $\sC_Z$ of \ref{halbH} and a motive version $X_{\sC}(F)$ of $X_{\sC}(\C)$ of \ref{halbH} and define a map   $Z(F) \to X_{\sC}(F)$ as in \ref{halbH}. But the author can not show that the canonical map $\Alb_{Z,n,\sC}(F) \to X_{\sC}(F)$ is bijective. We do not discuss these points in this paper. 

\end{sbpara}

\section{Curvature forms and Hodge theory}

\subsection{Reviews on curvature forms of line bundles}\label{s2.3.0}

\begin{sbpara}\label{231} Let $Y$ be a complex analytic manifold and let $L$ be a line bundle on $Y$. A metric $|\;|_L$ on $L$ is said to be $C^{\infty}$ if $|e|_L$ is a $C^{\infty}$ function for a local basis  $e$ of $L$. A $C^{\infty}$-metric on $L$ exists if $Y$ is paracompact. 

\end{sbpara}

\begin{sbpara}\label{232} Let $Y$ be as in \ref{231} and let $L$ be a line bundle on $Y$
endowed with a $C^{\infty}$ metric $|\;|_L$. Then the curvature form $\kappa(L)$ of $L$ is a $C^{\infty}$ $(1,1)$-form on $Y$ defined by
$$\kappa(L)=\kappa(L, |\;|_L):= \partial {\bar \partial}  \log(|e|_L^2)= \partial {\bar \partial}  \log ((e,e)_L)$$
where $e$ is a local basis of $L$ and  $(\;,\;)_L$ is the Hermitian form on $L$ corresponding to $|\;|_L$. 
 (Recall that if  $Y$ is $n$-dimensional and $(z_j)_{1\leq j\leq n}$ is a local coordinate of $Y$, for a $C^{\infty}$-function $g$ on  $Y$, $\partial {\bar \partial} g=\sum_{j,k}  (\frac{\partial}{\partial  z_j}\frac{\partial}{\partial \bar z_k} g)dz_j\wedge d\bar z_k$.)

This $\kappa(L)$ does not depend on the local choice of a basis $e$ and hence defined globally.

\end{sbpara}

\begin{sbpara}\label{233} Let $Y$ and $L$ be as in \ref{232}. Then the curvature form $\kappa(L)$ 
gives an Hermitian form on the tangent bundle $T_Y$ of $Y$. Locally, using a local coordinate $(z_j)_j$ of $Y$, 
this Hermitian form $(\;,\;)$ is given by
 $$(\frac{\partial}{\partial z_j}, \frac{\partial}{\partial z_k})= \frac{\partial}{\partial z_j}\frac{\partial}{\partial \bar z_k}\log((e,e)_L).$$

\end{sbpara}

\begin{sbpara}\label{234} Let $Y$ and $L$ be as in \ref{232} and assume that $Y$ is a connected compact Riemann surface. Then the theory of Chern forms and Chern classes shows
$$\frac{1}{2\pi i} \int_Y \kappa(L)= \text{deg}(L).$$

\end{sbpara}

\subsection{Review on the result of Griffiths on curvature forms}

  \begin{sbpara}\label{321}
Let $Y$ be a complex analytic manifold.  Let $\sH$ be a variation of polarized $\Q$-Hodge structure of weight $w$ on $Y$. Let $r\in \Z$, and consider the line bundle $\det(\fil^r\sH_{\sO})$ with the Hodge metric. The result of Griffiths on the curvature form $\kappa(\det(\fil^r\sH_{\sO}))$ is as follows. (Here $\fil$ is the Hodge filtration.)

Let $r\in \Z$. We have a homomorphism
$$h_r: T_Y\to \Hom_{\sO_Y}(\gr^r\sH_{\sO}, \gr^{r-1}\sH_{\sO})$$
which sends $\alpha\in T_Y$ to the composition $\gr^r\sH_{\sO}\overset{\nabla}\to \gr^{r-1}\sH_{\sO}\otimes_{\sO_Y} \Omega^1_Y\overset{\alpha}\to \gr^{r-1}\sH_{\sO}$. Here $\gr^r=\fil^r/\fil^{r+1}$, the first arrow is induced by the connection $\sH_{\sO}=\sO_Y\otimes_\Q \sH_\Q \to \sH_{\sO}\otimes_{\sO_Y} \Omega^1_Y=  \Omega^1_Y\otimes_\Q \sH_\Q$ which kills $\sH_{\Q}$, and the second arrow  is by the duality between $T_Y$ and $\Omega^1_Y$.

\medskip

Theorem (Griffiths (\cite{Gri})):
The Hermitian form $\kappa(\text{det}(\fil^r \sH_{\sO}))$ on $T_Y$ is described  as $(\alpha,\beta)\mapsto (h(\alpha), h(\beta))$ where the last $(\;,\;)$ denotes the Hodge metric of $(\gr^r\sH_{\sO})^* \otimes \gr^{r-1}\sH_{\sO}$.

\medskip

\end{sbpara}

\begin{sbpara}\label{h^*}  As is explained in Part I, 2.4.3, the curvature form $\kappa(\det(\fil^r\sH_{\sO}))$ is independent of the polarization. This can be seen also as follows. 

Let 
$$A_{\R,Y} \quad (\text{resp.} \; A_{\C,Y})$$
be the sheaf of $\R$ (resp. $\C$) valued $C^{\infty}$ functions on $Y$.

Let $\lambda_r$ be the composition of the bijections $$A_{\C,Y}\otimes_{\sO_Y} \gr^r \sH_{\sO}\cong (A_{\C,Y} \otimes_{\sO_Y} \sH_{\sO})^{r,w-r} \to (A_{\C,Y}\otimes_{\sO_Y} \sH_{\sO})^{w-r,r} \cong A_{\C,Y}\otimes_{\sO_Y} \gr^{w-r}\sH_{\sO}$$
where $(-)^{p,q}$ denotes the Hodge $(p,q)$-component and the middle bijection is induced by the complex conjugation $A_{\C,Y}\otimes_{\Q} \sH_{\Q} \to A_{\C,Y} \otimes_{\Q} \sH_{\Q}\;;\;a\otimes b \mapsto {\bar a}\otimes b$.

Let  $$h^*_r: T_Y\to A_{\C,Y}\otimes_{\sO_Y} \sH om_{\sO_Y}(\gr^{r-1} \sH_{\sO}, \gr^r\sH_{\sO})$$
be the map which sends $\alpha\in T_Y$ to the following composite map: 
$$A_{\C,Y}\otimes_{\sO_Y} \gr^{r-1}\sH_{\sO}\overset{\lambda_{r-1}}\to A_{\C,Y}\otimes_{\sO_Y} \gr^{w-r+1}\sH_{\sO}\overset{h_{w-r+1}(\alpha)}\longrightarrow A_{\C,Y}\otimes_{\sO_Y} \gr^{w-r}\sH_{\sO}\overset{\lambda_{w-r}}\to A_{\C,Y}\otimes_{\sO} \gr^r\sH_{\sO}.$$

Then $$(h_r(\alpha), h_r(\beta))=\langle h_r(\alpha), h^*_r(\beta)\rangle$$
where $\langle \;,\;\rangle$ denotes the duality $(\gr^r \sH_{\sO})^*\otimes \gr^{r-1}\sH_{\sO} \times 
 (\gr^{r-1}\sH_{\sO})^* \otimes \gr^r \sH_{\sO}\to \sO_Y$. The right hand side of this equation is defined without using the polarization. 
 \end{sbpara}

\begin{sbpara}\label{322} Consider the period domain $X(\C)= X_{G,\sG, H_b, K}(\C)$ with $G$ reductive.

Let $\sH_{X(C)}$ be universal object on $X(\C)$, let $V\in \Rep(G)$, and let $r\in \Z$. Consider  the restriction $\kappa(\text{det}(\fil^r \sH_{X(\C)}(V)_{\sO}))_{\hor}$ to $T_{X(\C),\hor}$ of the curvature form 
$\kappa(\text{det}(\fil^r \sH_{X(\C)}(V)_{\sO}))$ of the line bundle $\text{det}(\fil^r\sH_{X(\C)}(V)_{\sO})$ with the Hodge metric of a polarization (\ref{XCpol})  (the curvature form is independent of the choice of the polarization). Recall that $\sH_{X(\C)}(V)$ need not be a variation of Hodge structure (it need not satisfy the Griffiths transversality). Recall that
$T_{X(\C),\hor}=\gr^{-1}\sH_{X(\C)}(Lie(\sG))$. The action of $G$ on $V$ induces a Lie action $\Lie(\sG)\otimes V \to V$ of $\Lie(\sG)$ and hence a homomorphism $\gr^{-1}\sH_{X(\C)}(\Lie(\sG))_{\sO}\otimes \gr^r\sH_{X(\C)}(V)_{\sO}\to \gr^{r-1}\sH_{X(\C)}(V)_{\sO}$. 
Let $h_r: T_{X(\C),\hor}= \gr^{-1}\sH_{X(\C)}(\Lie(\sG))_{\sO}\to (\gr^r\sH_{X(\C)}(V)_{\sO})^* \otimes \gr^{r-1}\sH_{X(\C)}(V)_{\sO}$ be the induced homomorphism. Just as in \ref{h^*}, we have a homomorphism $h^*_r: T_{X(\C),\hor}\to A_{\C,X(\C)}\otimes_{\sO_{X(\C)}} (\gr^{r-1}\sH_{X(\C)}(V)_{\sO})^* \otimes \gr^r\sH_{X(\C)}(V)_{\sO}$.

\medskip

Theorem \ref{322}.1.  The Hermitian form $\kappa(\text{det}(\fil^r  \sH_{X(\C)}(V)_{\sO})_{\hor}$  on
 $T_{X(\C),\hor}$ coincides with $(\alpha,\beta)\mapsto (h_r(\alpha), h_r(\beta))=\langle h_r(\alpha), h_r^*(\beta)\rangle$ where the second $(\;,\;)$ is the Hodge metric of $(\gr^r \sH_{X(\C)}(V)_{\sO})^*\otimes \gr^{r-1}\sH_{X(\C)}(V)_{\sO}$ and $\langle \;,\;\rangle$ is the paring $(\gr^r\sH_{X(\C)}(V))_{\sO})^* \otimes \gr^{r-1}\sH_{X(\C)}(V)_{\sO} \times (\gr^{r-1}\sH_{X(\C)}(V))_{\sO})^* \otimes \gr^r\sH_{X(\C)}(V)_{\sO}\to \sO_{X(\C)}$. 
 
\medskip

This is proved in the same way as the theorem of Griffiths (\ref{321}), and is also reduced to that theorem by the following fact. For any $x\in X(\C)$ and $v\in T_{x, X(\C),\hor}$, there are a 
one-dimensional complex analytic space $Y$, $y\in Y$,  and a horizontal morphism $f: Y\to X(\C)$ such that $f(y)=x$ and $v$ is in the image of $T_{y,Y}\to T_{x,X(\C),\hor}$.

\end{sbpara}

\begin{sbpara} The pullback of $\kappa(\text{det}(\fil^r  \sH_{X(\C)}(V)_{\sO})_{\hor}$ to $T_{D(G,\sG, H_b),\hor}$ is invariant under the action of $\sG(\R)$. Recall that by our assumption $G$ is reductive, $D(G, \sG, H_b)$ is a finite disjoint union of $\G(\R)$-orbits (\ref{forb}).  

\end{sbpara}

\subsection{$X(\C)$ is like a hyperbolic space in the case $G$ is reductive}

\begin{sbpara}\label{hypb1} A complex analytic manifold $Y$ is said to be Brody hyperbolic if any morphism $\C\to Y$ from the complex plane $\C$ is a constant map. If $Y$ is hyperbolic in the sense of Kobayashi, then $Y$ is Brody hyperbolic. The converse is true in the case $Y$ is compact (\cite{Br}). 

Kobayashi conjectures (\cite{Kob} page 370) that if $Y$ is a compact hyperbolic complex manifold of dimension $n$,  the line bundle  $\wedge^n \Omega^1_Y$ (the canonical bundle) of $Y$ is ample. The author thinks that the experts believe the following non-compact version of this conjecture is true: 

Conjecture \ref{hypb1}.1. If $Y$ is Brody hyperbolic and  if $Y=\bar Y\smallsetminus D$ for a compact complex manifold  $\bar Y$ of dimension $n$ and for a normal crossing divisor $D$ on $Y$, then the line bundle  $\wedge^n \Omega^1_{\bar Y}(\log D)$ (the log canonical bundle) on $\bar Y$ is ample. 

The result \ref{Sh0} below shows that in the case $G$  is reductive, the period domains $Y=D(G, \sG, H_b)$ and $Y=X_{G, \sG, H_b, K}(\C)$ are ''like Brody hyperbolic'', which means that any horizontal morphisms from $\C$ to $Y$ are constant (notice that the condition ''horizontal'' is put on morphisms).  This is reduced to the case $Y$ is the classical period domain of Griffiths classifying polarized Hodge structures, and the author believes that this case of \ref{Sh0} is well known to the experts. Another result \ref{spadepos}, which is deduced from  \ref{322}.1, shows that if $G$ is reductive, $X_{G, \sG, H_b, K}(\C)$ satisfies  something like the ampleness condition in the above Conjecture \ref{hypb1}.1 (see \ref{amplog}).

\end{sbpara}

\begin{sbpara}\label{kspade}

Let $\kappa_{X(\C),\spadesuit}=\kappa(\text{det}(\Omega^1_{X(\C)}))_{\hor}$  be the Hermitian form on $T_{X(\C), \hor}$ with respect to the Hodge metric defined by the identification 
$\Omega^1_{X(\C)}$ as the dual of 
$\sH_{X(\C)}(Lie(\sG))_{\sO})/\fil^0)$ and by polarizations (\ref{XCpol}) on $\gr^W_w\sH_{X(\C)}(\Lie(\sG))$ for $w\in \Z$.
Here the curvature form is independent of the choices of polarizations (\ref{h^*}).

This $\kappa_{X(\C),\spadesuit}$ is regarded as an Hermitian form on the quotient $T_{X_{\red}(\C),\hor}$ of $T_{X(\C),\hor}$. 
\end{sbpara}

\begin{sbthm}\label{spadepos}  Assume $G$ is reductive. Then  $\kappa_{X(\C), \spadesuit}$ is positive definite as an Hermitian form on $T_{X(\C), \hor}$.

\end{sbthm}

\begin{pf} 
We have  $$\kappa(\text{det}(\Omega^1_{X(\C)}))=-\kappa(\text{det}(\sH_{X(\C)}(\Lie(\sG))_{\sO}/\fil^0))$$ $$= \kappa(\text{det}(\fil^0\sH_{X(\C)}(\Lie(\sG))_{\sO}))- 
\kappa(\text{det}(\sH_{X(\C)}(\Lie(\sG))_{\sO}))$$ and 
$\kappa(\text{det}(\sH_{X(\C)}(\Lie(\sG))_{\sO}))= \kappa(\sO_{X(\C)} \otimes_{\Q} \text{det}_\Q(\sH_{X(\C)}(\Lie(\sG))_\Q))=0$. Hence 
$$\kappa(\text{det}(\Omega^1_{X(\C)}))_{\hor}=  \kappa(\text{det}(\fil^0 \sH_{X(\C)}(\Lie(\sG))_{\sO}))_{\hor}.$$
By Theorem \ref{322}.1,  the Hermitian form $\kappa(\text{det}(\fil^0 \sH_{X(\C)}(\Lie(\sG))_{\sO}))_{\hor}$ on $T_{X(\C),\hor}=\gr^{-1}\sH_{X(\C)}(\Lie(\sG))_{\sO}$
is the pullback of a positive definite Hermitian form on \newline
$\sH om(\gr^0\sH_{X(\C)}(\Lie(\sG))_{\sO}, \gr^{-1}\sH_{X(\C)}(\Lie(\sG))_{\sO})$ by the map $\gr^{-1}\sH_{X(\C)}(\Lie(\sG))_{\sO}\to \sH om(\gr^0\sH_{X(\C)}(\Lie(\sG))_{\sO}, \gr^{-1}\sH_{X(\C)}(\Lie(\sG))_{\sO})$. Hence it is 
sufficient to prove that the last map is injective.  
Since $G$ is reductive, there is an isomorphism of Lie algebras $\Lie(G)\cong \Lie(\sG)\times \Lie(\sQ)$ which is compatible with the inclusion map $\Lie(\sG)\to \Lie(G)$. Hence it is sufficient to prove that the map 
$$\gr^{-1}\sH_{X(\C)}(\Lie(G))_{\sO} \to \sH om(\gr^0 \sH_{X(\C)}(\Lie(G))_{\sO}, \gr^{-1}\sH_{X(\C)}(\Lie(G))_{\sO})$$
is injective. Let $S_{\C/\R}\to G_\R$ be a homomorphism associated to 
$\sH$, let $h: \C=\Lie(S_{\C/\R})\to \Lie(G)_\R$ be the induced homomorphism of 
Lie algebras, and consider $h(i)\in \Lie(G)_\R$. We have 
$\Lie(G)_\R=\oplus_{r\in \Z} \Lie(G)_\R^{(r)}$ where 
$\Lie(G)_\R^{(r)}=\{x\in \Lie(G)_\R\;|\; [h(i),x]=rx\}$. 
It is sufficient to prove that $$\Lie(G)_\R^{(-1)}\to \Hom_\R(\Lie(G)_\R^{(0)}, \Lie(G)_\R^{(-1)})\;;\; x \mapsto (y\mapsto [x,y])$$ is injective. But this is seen by 
$h(i)\in \Lie(G)_\R^{(0)}$ and by the fact any $x\in \Lie(G)_\R^{(-1)}$ satisfies $[x, h(i)]= x$. 
\end{pf}

\begin{sbprop}\label{Sh0} Assume $G$ is reductive. Then 
any horizontal holomorphic map $\C\to X(\C)$ is a constant map. 
\end{sbprop}

\begin{pf}  By Schmid nilpotent orbit theorem \cite{SW}, this holomorphic map $f$ gives a nilpotent orbit at $\infty\in {\bf P}^1(\C)\supset \C$. But since $\pi_1(\C)=\{1\}$, this nilpotent orbit has no local monodromy and hence no degeneration. Hence we get $f: {\bf P}^1(\C)\to X(\C)$. But as is well know, a variation of pure Hodge structure $\sH$ on ${\bf P}^1(\C)$ is constant. (This is proved by using the fact that for any $r\in \Z$, $\text{deg}(\fil^r\sH_{\sO})\geq 0$ (a consequence of the theorem of Griffiths in \ref{321}) and hence $\text{dim}(H^0({\bf P}^1(\C), \fil^r\sH_{\sO}))\geq \text{rank}(\fil^r \sH_{\sO})$ by Riemann-Roch.) Hence $f$ is constant.
\end{pf}

\subsection{Height pairings in Hodge theory and curvature forms}\label{ss2.4}

\begin{sbpara}\label{htp1} Let $Y$ be a complex analytic manifold. 
Assume we are given a variation of $\Z$-Hodge structure $\sH_0$ of weight $-1$ on $Y$ and variations of mixed $\Z$-Hodge structure $\sH_1$ and $\sH_2$ on $Y$ with exact sequences $$\text{(i)} \quad 0 \to \sH_0\to \sH_1\to  \Z\to  0\quad \text{and} \quad \text{(ii)} \;  0\to \sH_0^*(1)  \to \sH_2\to \Z\to 0$$ where $\sH_0^*$ denotes the $\Z$-dual of $\sH_0$. Note that the exact sequence (ii) corresponds to an exact sequence (iii) $0\to \Z(1) \to \sH_2^*(1) \to \sH_0\to 0$.

Then we have a line bundle on $Y$ with a $C^{\infty}$-metric as follows (see for example \cite{BP}).

 Locally on $Y$,  we have a variation $\sH$ of mixed Hodge structure on $Y$ which satisfies $\gr^W_{-1}\sH= \sH_0$, $\sH/W_{-2}\sH=\sH_1$ (with the exact sequence (i)), $W_{-1}\sH=\sH_2^*(1)$ (with the exact sequence (iii)). 
 
 The isomorphism classes of such $\sH$ form a ${\bf G}_m$-torsor on $Y$ (defined globally on $Y$).  The action of ${\bf G}_m$ is as follows. A local section $s$ of ${\bf G}_m$ sends the class of $\sH$ to the Baer sum of the following two extensions of $\sH_1$ by $\Z(1)$. One is $\sH$, and the other is the extension of $\sH_1$ by $\Z(1)$ which is obtained from the extension of $\Z$ by $\Z(1)$ corresponding to $s$ and from the projection $\sH_1\to \Z$. 
 
 Hence we obtain a line bundle $L(\sH_1, \sH_2)$ on $Y$ associated to this ${\bf G}_m$-torsor. This line bundle has the following $C^{\infty}$ metric $\text{class}(\sH)\mapsto |\sH|$ characterized by the property that locally there is 
  a lifting $v$ of $1\in \gr^W_0\sH_\R$ in $A_{\R,Y} \otimes_\R \sH_\R$ such that $v+ \log(|\sH|)\in A_{\C,Y} \otimes_{\sO_Y} \fil^0\sH_{\sO}$ in 
  $A_{\C,Y}\otimes_{\sO_Y}\sH_{\sO}=A_{\C,Y}\otimes_{\R} \sH_\R$. Here $\log(|\sH|)\in A_{\R,Y}$ is  regarded as a local section of $A_{\C,Y}\otimes_\R W_{-2}\sH_\R$ 
  via the isomorphism 
  $A_{\C,Y}\cong A_{\C,Y}\otimes_\R W_{-2}\sH_{\sO}$ which is induced by $\Z \cdot 2\pi i\cong W_{-2}\sH_\Z$. (Hence $\log(|\sH|)$ is regarded as a local section of  $i\cdot A_{\R,Y}\otimes_\R W_{-2}\sH_\R$. )

 In this Section \ref{ss2.4}, we compute the curvature form of $L(\sH_1, \sH_2)$.   

This line bundle $L(\sH_1, \sH_2)$ with metric is related to the theory of height pairing as is explained in the Section \ref{ss2.5} later.

 \end{sbpara}

 \begin{sbpara}\label{htp2}

(1) We have an $\sO_Y$-homomorphism $$h:T_Y\to A_{\C,Y}\otimes_{\sO_Y} \gr^{-1}\sH_{0,\sO}$$ defined as follows.
For $\alpha\in T_Y$, consider the composition $\gr^0 \sH_{1,\sO}\overset{\nabla}\to \gr^{-1}\sH_{1,\sO}\otimes_{\sO_Y} \Omega^1_Y\overset{\alpha}\to \gr^{-1}\sH_{1,\sO_Y}=\gr^{-1}\sH_{0,\sO}$. Define $h(\alpha)$ as the image of $1\in \R$ under $\R\to \A_{\C,Y}\otimes_{\sO_Y}\gr^0\sH_{1,\sO}\overset{\alpha}\to A_{\C,Y}\otimes_{\sO_Y} \gr^{-1}\sH_{0,\sO}$ where the first  arrow is induced by the inverse of the isomorphism $A_{\C,Y} \otimes_{\sO_Y} \fil^0\sH_{1, \sO} \cap A_{\R,Y} \otimes_\R \sH_{1,\R}\overset{\cong}\to A_{\R,Y} \otimes_{\R} \gr^W_0\sH_{1,\R}$.

(2) We have $$h^*: T_Y\to A_{\C,Y}\otimes_{\sO_Y} \gr^0 \sH_0^*(1)_{\sO}$$ 
as the composition of $T_Y\overset{h}\to A_{\C,Y}\otimes_{\sO_Y} \gr^{-1}\sH_0^*(1)_{\sO} \overset{\lambda_{-1}}\to A_{\C,Y}\otimes_{\sO_Y} \gr^0\sH_0^*(1)_{\sO}$ where $h$ is defined as in (1) replacing $0\to \sH_0\to \sH_1\to \Q\to 0$ by $0\to \sH_0^*(1) \to \sH_2\to \Q\to 0$, and $\lambda_{-1}$ is as in \ref{h^*} for $\sH_0^*(1)$.

\end{sbpara}

\begin{sbthm}\label{htp3}
The Hermitian form $\kappa(L(\sH_1, \sH_2))$ on $T_Y$ coincides with $(\alpha, \beta)\mapsto  \langle h(\alpha), h^*(\beta)\rangle$ where $\langle\,\;\rangle$ denotes the canonical pairing $A_{\C,Y}\otimes_{\sO_Y} \gr^{-1}\sH_{0,\sO}\times A_{\C,Y}\otimes_{\sO_Y}\gr^0 \sH_0^*(1)_{\sO}\to A_{\C,Y}$.

\end{sbthm}

\begin{sbpara}\label{htp4}

The proof of \ref{htp3} is the reduction to the case of a period domain.

Fix a free $\Z$-module  $H_{0,\Z}$ of finite rank, and let $H_{1,\Z}=\Z e_1 \oplus H_{0,\Z}$ and $H_{2,\Z}= \Z e_2 \oplus H_{0,\Z}^*(1)$ be free $\Z$-module of rank $\text{rank}_\Z(H_{0,\Z})+1$. Here $H_{0,\Z}^*(1)= \Hom_\Z(H_{0,\Z}, \Z\cdot 2\pi i)$.
Define an increasing filtrations on $H_{1,\Q}=\Q\otimes_\Z H_{1,\Z}$ and on $H_{2,\Q}= \Q\otimes_\Z H_{2,\Z}$ by 
$$W_{-2}=0\subset W_{-1}= H_{0, \Q} \subset W_0= H_{1,\Q}, \quad W_{-2}=0\subset W_{-1}= H_{0,\Q}^*(1) \subset W_0= H_{2,\Q}.$$
Fix integers $h(r)\in \Z_{\geq 0}$ for $r\in \Z$ such that $\sum_r h(r)=\dim_\Q H_{0,\Q}$ and $h(-1-r)=h(r)$ for all $r$. Let $D$ be the set of all pairs $(\fil_{(1)}, \fil_{(2)})$, where $\fil_{(i)}$ is a  decreasing filtration on $H_{i,\C}$ satisfying the following conditions (i)--(iii). 

\medskip

(i) The restriction of $\fil_{(2)}$ to  $H_{0,\Z}^*(1)_{\C}$ coincides with the filtration induced by the restriction of $\fil_{(1)}$ to $H_{0,\C}$. 

(ii) For $i=1,2$, $(H_{i, \Z}, W,  \fil_{(i)})$ is a mixed Hodge structures. 

(iii) $\dim \gr^rW_{-1}\fil_{(i)}= h(r)$ for all $r$ ($i=1,2$).

\medskip

Then $D$ is naturally regarded as a complex analytic object. On $D$, we have the universal objects $\sH_{0,D}$, $\sH_{1,D}$, $\sH_{2,D}$ of $\Q$MHS$(D)$ (\ref{Drep}) (but these need not belong to $\Q$VMHS$(D)$ (\ref{varH})) with exact sequences $0\to \sH_{0,D}\to \sH_{1,D}\to \Z\to 0$ and $0\to \sH_{0,D}^*(1) \to \sH_{2,D}\to \Z\to 0$. 
By the method of \ref{htp1}, we have a line bundle $L(\sH_{1,D}, \sH_{2,D})$ on $D$ with a $C^{\infty}$ metric. 

The tangent bundle of $D$ is identified with $\sE /\fil^0$, where $\sE$ is a part of $ \sE nd_{\sO_D}(\sH_{1, D, \sO}, W) \times \sE nd_{\sO_D}(\sH_{2,D,\sO},W) $ consisting of all pairs $(a,b)$ such that the restriction of $a$ to $\sH_{0,D,\sO}$ and the restriction of $b$ to $\sH_{0,D}^*(1)_{\sO}$ are induced from the other, 
and $\fil$ is the Hodge filtration.

The horizontal tangent bundle $T_{D,\hor}$ of $D$ is defined to be $\fil^{-1}/\fil^0$ of $\sE$. 

Let $\kappa(L(\sH_{1,D}, \sH_2))$ be the curvature form of the metric of $L(\sH_{1,D}, \sH_{2,D})$ and let $\kappa(L(\sH_{1,D}, \sH_{2,D}))_{\hor}$ be its restriction to $T_{D,\hor}$. 

By the method of \ref{htp2}, we have
$$h: T_{Y, \hor} \to A_{\C,D} \otimes_{\sO_Y} \gr^{-1}\sH_{0,D, \sO}, \quad h^*:T_{Y,\hor}\to A_{\C,D} \otimes_{\sO_D} \gr^0\sH_{0,D,\sO}. $$

\end{sbpara}

\begin{sbprop}\label{htp5} The Hermitian form $\kappa(L(\sH_{1,D}, \sH_{2,D}))_{\hor}$ on  $T_{D,\hor}$ coincides with   $(\alpha,\beta)\mapsto \langle h(\alpha), h^*(\beta)\rangle$.

\end{sbprop}

\begin{sbpara}
Theorem \ref{htp3} is reduced to \ref{htp5} as follows. We may assume that $Y$ is connected. Take $b\in Y$ and  define the above period domain $D$ by  
taking the stalks $\sH_{j,\Z,y}$ at $b$ as $H_{j,\Z}$ ($j=0,1,2$)  and by taking the Hodge numbers of $\sH_0$ as $h(r)$. Then we have the period map $\tilde Y\to D$ from the 
universal covering $\tilde Y$ of $Y$ to $D$ associated to $(\sH_0, \sH_1, \sH_2)$, and the pullback of $\kappa(L(\sH_1, \sH_2))$ on $\tilde Y$ coincides with the pullback 
of $\kappa(L(\sH_{1,D}, \sH_{2,D}))_{\hor}$. Hence \ref{htp3} follows from \ref{htp5}. 
\end{sbpara}

\begin{sbpara} For the proof of \ref{htp5}, we use several period  domains  related to $D$. 

We first define the period domain $\tilde D$ and describe the ${\bf G}_m$-torsor on $D$ associated to $L(\sH_{1,D}, \sH_{2,D})$ explicitly. Let $H_{0,\Z}$ be as in \ref{htp4} and consider a free $\Z$-module $H_{\Z}:= \Z e \oplus H_{0, \Z} \oplus \Z e'$ of rank $\text{rank}_{\Z}(H_{0,\Z})+2$.
Define the increasing filtration $W$ on $H_\Q= \Q\otimes_\Z H_{0,Z}$ by 
$$W_{-3}=0\subset W_{-2}=\Q e' \subset W_{-1}= H_{0,\Q} \oplus \Q e' \subset W_0=H_{0,\Q}.$$
Let $\tilde D$ be the set of all decreasing filtrations $\varphi$ on $H_\C=\C\otimes_\Z H_\Z$ such that
$(H_\Z, W, \varphi)$ is a mixed Hodge structure and such that $\dim \gr^r \gr^W_{-1}\varphi= h(r)$ for all $r$. 
Then $\tilde D$ is naturally regarded as a complex analytic manifold. 

We have a surjective morphism 
$$\tilde D \to D\;;\;\varphi\mapsto (\fil_{(1)}, \fil_{(2)}).$$ Here  $\fil_{(1)}$ is the filtration induced on $H_\C/W_{-2}H_\C$ by $\varphi$ where we identify $e$ with $e_1$ and $\fil_{(2)}$ is the filtration induced by $\varphi$ by identifying $H_{2,\Z}^*(1)_\C$ with $W_{-1}H_\C$ where we identify $e'$ with $e_2^* \otimes 2\pi i$. 
Consider the action of the additive group $\C$ on $H_\C$ as follows: $z\in \C$ sends $e$ to $e+z\cdot (2\pi i)^{-1}e'$ and fix all elements of $W_{-1}H_\C$. This action induces an action of $\C$ on $\tilde D$ and induces an isomorphism $\C\bs \tilde D \overset{\cong}\to D$ and $\tilde D$ is a $\C$-torsor over $D$. Via $\C/\Z(1) \cong \C^\times\;;\; z \mapsto \exp(z)$ ($\Z(1)=\Z \cdot 2\pi i$), $\Z(1)\bs \tilde D$ becomes a $\C^\times$-torsor over $D$.

This $\Z(1)\bs \tilde D \to D$ is  the ${\bf G}_m$-torsor associated to the line bundle $L(\sH_{1,D}, \sH_{2,D})$ on $D$. Its metic  is $\tilde D\ni \varphi \mapsto |\varphi|$ where $|\varphi|\in \R_{>0}$ is characterized by the following property.  There is 
an element $v\in H_{0,\R}$ whose image in $\gr^W_0H_\R$ is the class of $e$ such that $v+\log(|\varphi|)\cdot (2\pi i)^{-1}e' \in \varphi^0$.

\end{sbpara}

\begin{sbpara} We define period domains $D'$, $D^{(i)}$ ($i=1,2$), and $D(s)$. Let the notation be as in \ref{htp4}.

Let $D'$ be the set of all decreasing filtrations $\varphi$ on $H_{0,\C}$ such that $(H_{0,\Z}, \varphi)$ is a Hodge structure of weight $-1$ and $\dim_\C(\gr^r\varphi)= h(r)$ for all $r$. We have a projection $p: D\to D'\;;\;(\fil_{(1)}, \fil_{(2)}) \mapsto \varphi$ 
where $\varphi$ is the restriction of $\fil_{(1)}$ to $H_{0,\C}$. For $s\in D'$, let $D(s)\subset D$ be the inverse image of $s$ in $D$. 

For $i=1,2$, let $D^{(i)}$ be the set of all decreasing filtrations $\varphi$ on $H_{i,\C}$ such that $(H_{i, \Z}, W, \varphi)$ is a mixed Hodge structure and $\dim_\C(\gr^r\gr^W_{-1}\varphi)= h(r)$ for all $r$. Let $p_i: D\to D^{(i)}$ be the map $(\fil_{(1)}, \fil_{(2)}) \to \fil_{(i)}$.

 For $x=(\fil_{(1)}(x), \fil_{(2)}(x))\in D$, we have a map $f_{x,i}: D^{(i)}\to D\;;\;\varphi\mapsto (\fil_{(1)},\fil_{(2)})$ defined as follows. Let $S^{(i)}$ be the one-dimensional $\R$-space $\fil_{(i)}(x)^0 \cap H_{i,\R}$. 
 Then for $\varphi\in D^{(i)}$, $\fil_{(i)}=\varphi$ and $\fil_{(3-i)}$ is characterized by the property that 
   $\fil_{(3-i)}^0=\varphi^0\cap W_{-1}H_{3-i,\C}\oplus \C\otimes_\R S^{(3-i)}$. 
The composition $D^{(i)}\overset{f_{x,i}}\to D \overset{p_i}\to D^{(i)}$ is the identity map, and  we have $f_{x,i}(p_i(x))=x$. 
 
 Let $g_{x,i}:D'\to D^{(i)}$ be the following map $\varphi' \mapsto \varphi$. $\varphi$ is characterized by the properties that $W_{-1}\varphi=\varphi'$ iand $\varphi^0 =(\varphi')^0\oplus \C\otimes_\R S^{(i)}$. 
 Then $g_{x,i}$ sends  $p(x)\in D'$ to $p_i(x)$. 
We have $f_{x,1}\circ g_{x,1}=f_{x,2}\circ g_{x,2}: D'\to D$. We denote this map $D'\to D$ by $h_x$.  
The composition $D'\overset{h}\to D \overset{p}\to D'$ is the identity map and $h_x(p(x))=x$.

We consider the tangent bundles and horizontal tangent bundles of  these period domains.

We have $T_{D'}= (\sO_{D'} \otimes_\C \text{End}_{\C}(H_{0,\C}))/\fil^0$ where $\fil$ is the Hodge filtration. Define $T_{D',\hor}= \fil^{-1}/\fil^0\subset T_{D'}$. We have $T_{D^{(i)}}= (\sO_{D^{(i)}}\otimes_\C \text{End}_{\C}(H_{i,\C}))/\fil^0$ where $\fil$ is the Hodge filtration. Define $T_{D^{(i)},\hor}= \fil^{-1}/\fil^0\subset T_{D^{(i)}}$. We have $T_{D(s)}= (W_{-1}$ part of the pullback of $T_D$ to $D(s)$). Let $T_{D(s),\hor}$ be the $W_{-1}$ part of the pullback of $W_{-1}T_{D,\hor}$ to $D(s)$.

For $x\in D$ with  $s:=p(x)\in D'$, we have isomorphisms
$$T_{s,D'} \oplus T_{x,D(s)}\overset{\cong}\to T_{x,D}, \quad T_{s,D',\hor}\oplus T_{x,D(s),\hor}\overset{\cong}\to T_{x,D,\hor}$$
induced by the map $h_x:D'\to D$ and the inclusion map $D(s)\to D$.

Hence for the proof of \ref{htp5}, it is sufficient to prove the pullbacks of 
$\kappa(L(\sH_{1,D}, \sH_{2,D}))$ to $T_{D^{(i)},\hor}$ ($i=1,2$, under $f_{x,i}$) and to $T_{D(s),\hor}$ are described as in \ref{htp5}.

\end{sbpara}

\begin{sbpara}
We consider first $D^{(i)}$. As is easily seen, under the map $f_{x,i}:D^{(i)}\to D$, the pullback of $L(\sH_{1,D}, \sH_{2,D})$ with the metric is isomorphic to the trivial line bundle $\sO_{D^{(i)}}$ with the standard metric. Hence its curvature form is zero. On the other hand, the pullback of the map $h$  on $T_{D^{(2)},\hor}$ under $f_{x,2}$ is the zero map and the pullback of the map $h^*$ on $T_{D^{(1)},\hor}$ under $f_{x,1}$ is the zero map.

\end{sbpara}

\begin{sbpara} Fix $s\in D'$. Let $\tilde D(s)$ be the inverse image of $D(s)$ in $\tilde D$. 
Let $U$ be a $\C$-subspace of $H_{0,\C}$ such that $H_{0,\C}=\fil(s)^0\oplus U$ and let $U'$ be a $\C$-subspace of $H^*_{0,\Z}(1)_\C$ such that $H^*_{0,\Z}(1)_\C=(\fil(s))^*(1)^0\oplus U'$. Then we have  isomorphisms 
$$U\times U' \overset{\cong}\to D(s)\;;\; (z,w) \mapsto (\fil(z), \fil(w)),$$
$$U\times U' \times \C\overset{\cong}\to \tilde D(s)\;;\; (z, w, u) \mapsto \varphi(z,w,u).$$ 
Here    $\fil(z)$ on $H_{1,\C}$ extends $\fil(s)$ on $H_{0,\C}$ by $e+z\in \fil(z)^0$, $\fil(w)$ on $H_{2,\C}$ extends $\fil(s)^*(1)$ on $H_{0,\Z}^*(1)_\C$ by $e+w\in \fil(w)^0$, and $\varphi(z,w,u)$ extends the filtration $\fil(w)^*(1)$ on $W_{-1}H_\C$ by $e+z+u\cdot (2\pi i)^{-1}e' \in \varphi(z,w,u)^0$.

These isomorphisms are compatible with the projection $\tilde D\to D$.

\end{sbpara}

\begin{sbpara}\label{hp1}  By using  $H_{0,\C}=  \fil(s)^0 \oplus H_{0, \R}$, write
$z=z_1+z_2$, $(z\in H_{0,\C}$, $z_1\in \fil(s)^0$. $z_2\in H_{0,\R}$).  
Similarly, by using $H_{0,\Z}^*(1)_\C= \fil(s)^*(1)^0 \oplus H_{0,\Z}^*(1)_\R$, write
$w=w_1+w_2$, ($w\in H_{0,\Z}^*(1)_\C$, $w_1\in \fil(s)^*(1)^0$, $w_2\in H_{0,\Z}^*(1)_\R$). 

We also use the notation
$$h= \text{Re}(h) + i \text{Im}(h) \quad \text{for}\; h\in H_{0,\C}, \text{where Re}(h), \text{Im}(h)\in H_{0,\R}.$$

Let $\langle\;,\;\rangle: H_{0,\R}\times H_{0,\Z}^*(1)_\R \to \R(1)$ be the canonical pairing. We denote the induced $\C$-linear pairing $H_{0,\C}\times H_{0,\Z}^*(1)_\C \to \C$ also by $\langle\;,\;\rangle$. 
\end{sbpara}

\begin{sblem}\label{hp2}  $\log(|\varphi(z,w,u)|)= i \langle\text{Im}(z), w_2\rangle+ \text{Re}(u).$

(Note that $\langle\text{Im}(z), w_2\rangle\in \R(1)$ and hence $ i \langle\text{Im}(z), w_2\rangle\in \R$.)

\end{sblem}

\begin{pf}  $\varphi(z,w,u)^0$ is generated by  $e+z+ u\cdot (2\pi i)^{-1}e'$ and $\lambda- \langle \lambda, w\rangle\cdot (2\pi i)^{-1} e'$ ($\lambda \in \fil(s)^0$). In particular, $z_1-\langle z_1, w\rangle \cdot (2\pi i)^{-1}e'\in \varphi(z,w,u)^0$. Hence $e+z_2+ (\langle z_1, w\rangle+u)\cdot (2\pi i)^{-1}e'\in \varphi(z,w,u)^0$.
This shows $\log(|\varphi|)= \text{Re}(z_1, w\rangle + \text{Re}(u)$. 
We have $\langle z_1,w_1\rangle=0$ because $\langle \fil(s)^0, (\fil(s)^*(1))^0\rangle=0$. Hence
$\text{Re}(\langle z_1, w\rangle)= \text{Re}(\langle z_1, w_2\rangle)=i \langle \text{Im}(z_1), w_2\rangle= i\langle \text{Im}(z), w_2\rangle$. \end{pf}

\begin{sbpara} On $\tilde D(s)$, we have $\partial {\bar \partial} \text{Im}(u)=0$. Hence $2\partial {\bar \partial}\ \log(|\varphi|)$ descents to $D(s)$ and this is the pullback $\kappa(L(\sH_{1,D(s)}, \sH_{2,D(s)}))$ of $\kappa(L(\sH_{1,D}, \sH_{2,D}))$ to $D(s)$ where $\sH_{j,D(s)}$ ($j=1,2$) denotes the pullback of $\sH_{j,D}$ to $D(s)$. 
It is equal to $2i\partial {\bar \partial}   \langle \text{Im}(z), w_2\rangle$. 
\end{sbpara}

\begin{sblem}\label{le4} On $D(s)$, we have  $\partial {\bar \partial}  z_1=\partial {\bar \partial}   z_2=0$. $\partial {\bar \partial} w_1=\partial {\bar \partial}  w_2=0$

\end{sblem}

\begin{pf} Write the $\R$-linear map $\text{Im}: \fil^0H_{0,\C}\to H_{0,\R}$ by $\ell$. Then $\ell$ is a bijection. Since $z_1=\ell^{-1}\text{Im}(z_1)=\ell^{-1}\text{Im}(z)$,  we have $\partial {\bar \partial}z_1= \ell^{-1}(\partial {\bar \partial} \text{Im}(z))=0$. The statement for $z_2$ follows from this by $\bar \partial z=0$. Results for $w_i$ are proved similarly. 
\end{pf}

\begin{sbpara}\label{hp3} By \ref{hp2} and \ref{le4}, we have
$2{\partial \bar \partial}\log(|\varphi|)=  \langle \partial z, \bar \partial w_2\rangle-\langle \bar \partial \bar z, \partial w_2\rangle$ in $A^2_{\C,D}$. Here $A^p_{\C,D}$ denotes the sheaf of complex valued $C^{\infty}$ $p$-forms on $D$, and $\langle\;,\;\rangle$ denotes the pairing
$A^1_{\C,D}\otimes_\R H_{0,\R}\otimes A^1_{\C, D}\otimes_\R H_{0,\Z}^*(1)_\R\to A_{\C,D}\;;\;(\omega\otimes h, \omega' \otimes h')\mapsto (\omega\wedge \omega') \otimes \langle h, h'\rangle$.

\end{sbpara}

\begin{sbpara} Let $\Omega^1_{\hor}$ be the quotient of $\Omega_{D(s)}^1$ corresponding to the subbundle $T_{D(s),\hor}$ of $T_{D(s)}$ by duality. Let $A^1_{\C,D,\hor}$ be the quotient
$A_{\C,D} \otimes_{\sO_D} \Omega^1_{hor} \oplus (A_{\C,D} \otimes_{\sO_D} \Omega^1_{\hor})^{-}$
of $A^1_{\C,D}=A_{\C,D} \otimes_{\sO_D} \Omega^1_D \oplus (A_{\C,D} \otimes_{\sO_D} \Omega^1_D)^{-}$ where $(\;)^-$ denotes the complex conjugate. 
Consider
$$\nabla_{\hor}: \sH_{1,D(s), \sO}=\sO_{D(s)}\otimes_\R H_{1,\R} \to \sH_{1,D(s),\sO}\otimes_{\sO_{D(s)}} \Omega^1_{\hor}= \Omega^1_{\hor} \otimes_\R H_{1,\R}$$
which satisfies the Griffiths transversality. This induces
$$\nabla_{\hor}=(\partial, \bar \partial): A_{\C, D(s)} \otimes_\R H_{1,D(s),\R} \to A^1_{\hor} \otimes_\R \sH_{0,D,\R}$$ $$= ((A^1_{\C, D(s)} \otimes_{\sO_D} \Omega^1_{\hor}) \oplus (A_{\C, D(s)}\otimes_{\sO_D} \Omega^1_{\hor})^{-})\otimes_\R H_{0,\R}.$$

  In the rest of the proof of \ref{htp5} below, $\nabla_{\hor}$, $\partial$, $\bar\partial$ ($\nabla=\partial+\bar \partial$) are considered by using $\Omega^1_{\hor}$ and $A^1_{\hor}$, not using $\Omega^1_{D(s)}$ and $A^1_{\C,D}$. 

\end{sbpara}
 
 \begin{sblem}\label{zwhor} On $D(s)$, we have:
 
 (1) $\partial z\in \fil^{-1}\sH_{0,D,\sO}\otimes_{\sO_D} \Omega^1_{\hor}$. 
 
 (2) $\nabla_{\hor} z_1\in \fil^0\sH_{0,D,\sO} \otimes_{\sO_D} A^1_{\hor}$.
 
 (3) $\partial z_2\in \fil^{-1}\sH_{0,D,\sO} \otimes_{\sO_D} A^1_{\hor}$. 
 
 (4) $\bar \partial z_2 \in \fil^0\sH_{0,D,\sO} \otimes_{\sO_D} A^1_{\hor}$.
 
 (5) The images of $\partial z$ and $\partial z_2$ in $\gr^{-1}\sH_{0,D,\sO}\otimes_{\sO_D} A^1_{\hor}$ coincide.

 We have similar results for $w,w_1,w_2$. 
 \end{sblem}

\begin{pf}
Since $e+z\in \fil^0\sH_{1,D(s), \sO}$ and $\nabla_{\hor}(e+z)= \nabla z$,  we have (1) by Griffiths transversality for $\Omega^1_{\hor}$. 
 (2) follows from the fact that the Hodge filtration of $\sH_{0,D(s)}$ is constant. By (1) and (2), we have (3). 
 (4) follows from (2) and $\bar \partial z=0$. (5) follows from (2). 
 \end{pf}

\begin{sbpara} By \ref{hp3} and by (4) and (5) of \ref{zwhor}, $2i\partial {\bar \partial} \log(|\varphi|)= 
\langle (\partial z_2)^{-1,0}, (\bar \partial w_2)^{0,-1}\rangle - \langle (\bar \partial z_2)^{(0,-1)}, (\partial w_2)^{-1,0}\rangle$.
Here $(\;)^{p,q}$ denotes the Hodge $(p,q)$-component.
Since $z_2$ and $w_2$ are  real, $(\bar \partial w_2)^{0,-1}$ is the complex conjugate of $(\partial w_2)^{-1,0}$ and $(\bar \partial z_2)^{0,-1}$ is the complex conjugate of $(\partial z_2)^{-1,0}$.

This proves \ref{htp5}. 

\end{sbpara}

\begin{sbpara}\label{QandZ} We can develop similarly the theory of ${\bf G}_m\otimes \Q$-torsors with a $C^{\infty}$-metric and its curvature forms. A ${\bf G}_m\otimes \Q$-torsor corresponds to a pair $(L, n)$ with $L$ a line bundle and with $n$ an integer $\geq 1$, which we regard as $L^{\otimes 1/n}$. 

We have the evident variant of  \ref{htp3} for a variation of $\Q$-Hodge structure  $\sH_0$ of weight $-1$ and for variations of $\Q$-MHS $\sH_1$, $\sH_2$ ($\Z$-coefficients are replaced by $\Q$-coefficients) and for the associated ${\bf G}_m\otimes \Q$-torsor which has a $C^{\infty}$-metric.  This variant can be proved by the following simple reduction to \ref{htp3}:   Locally, such $\sH_0$ has a $\Z$-structure. For some $n\geq 1$, for $i=1,2$, replace  $\sH_1$  (resp. $\sH_2$) by its $n$ fold Baer sum as  an extension of $\Q$ by $\sH_0$ (resp. $\sH_0^*(1)$) which has $\Z$-coefficients. Then $n^{-2}$ times the curvature form of the associated ${\bf G}_m$-torsor coincides with the curvature form of the original ${\bf G}_m\otimes \Q$-torsor. 

\end{sbpara}
\begin{sbpara}\label{XCdw1} Consider our period domain  $X(\C)$. We define an Hermitian form
$\kappa_{X(\C), \diamondsuit, w,1}$ on $T_{X(\C), \hor}$.

Consider the variation of $\Q$-MHS 
$\sH_0= (\gr^W_w\sH_{X(\C)}(V_0))^*\otimes\gr^W_{w-1}\sH_{X(\C)}(V_0)$ on $X(\C)$ of weight $-1$. The exact sequence $$0\to \gr^W_{w-1}\sH_{X(\C)}(V_0) \to W_w\sH_{X(\C)}(V_0)/W_{w-2}\sH_{X\C)}(V_0) \to \gr^W_w\sH_{X(\C)}(V_0)\to 0$$
gives an exact sequence
$0\to \sH_0 \to \sH_1\to \Q\to 0$ in $\Q$MHS$(X(\C))$. By the isomorphism  $\sH_0\cong \sH_0^*(1)$ obtained by the canonical polarizations (\ref{XCpol}) of  $\gr^W_w\sH_{X(\C)}(V_0)$ and $\gr^W_{w-1}\sH_{X(C)}(V_0)$, we have an exact sequence $0\to \sH_0^*(1)\to \sH_1 \to \Q\to 0$.

By \ref{QandZ},  we obtain a ${\bf G}_m\otimes \Q$-torsor $L(\sH_1, \sH_1)$ with $C^{\infty}$-metric.  We define 
$$\kappa_{X(\C), \diamondsuit, w,1}:= \kappa(L(\sH_1, \sH_1))_{\hor},$$
the restriction  of the curvature form of $L(\sH_1, \sH_1)$ to $T_{X(\C), \hor}$. 

To describe $\kappa_{X(\C), \diamondsuit, w,1}$, let 
$$h: T_{X(\C),\hor}\to A_{\C,X(\C)}\otimes_{\sO_{X(\C)}}\gr^{-1}\sH_{0,\sO}$$
be the map which sends $\alpha\in T_{X(\C), \hor}$ to the image of $1\in \R$ under
$$\R\to A_{\C,X(\C)}\otimes_{\sO_{X(\C)}} \gr^0\sH_{1,\sO}\overset{\alpha}\to A_{\C,X(\C)}\otimes_{\sO_{X(\C)}} \gr^{-1}\sH_{1,\sO} = \gr^{-1}\sH_{0,\sO}.$$
Here the first arrow comes from the inverse of the isomorphism 
$$A_{\C,X(\C)}\otimes_{\sO_{X(\C)}} \fil^0\sH_{1,\sO}\cap A_{\R,X(\C)} \otimes_\R \sH_{1,\R}\overset{\cong}\to 
A_{\R,X(\C)}\otimes_\R \gr^W_0\sH_{1,\R}.$$  
\end{sbpara}

\begin{sbprop}\label{htp7}
(1) The Hermitian form $\kappa_{X(\C), \diamondsuit,w,1}$ coincides with $(\alpha, \beta)\mapsto (h(\alpha), h(\beta))$ where the last 
$(\;,\;)$ is the Hodge metric on $\gr^{-1}\sH_{0,\sO}$.

(2) This Hermitian form $\kappa_{X(\C),\diamondsuit,w,1}$ on $T_{X(\C)}$ is semi-positive definite. It comes from an Hermitian form on the quotient $\gr^{-1}((W_0/W_{-2} \;\text{of}\; \sH_{X(\C)}(\Lie(\sG))_{\sO})$ of $T_{X(\C),\hor}=\gr^{-1}\sH_{X(\C)}(\Lie(\sG))_{\sO}$. The Hermitian form $\sum_{w\in \Z} \kappa_{X(\C), \diamondsuit, w,1}$ on \newline $\gr^{-1}\gr^W_{-1}\sH_{X(\C)}(\Lie(\sG))_{\sO}$   is positive definite. 

\end{sbprop}
\begin{pf} (1) follows from the ${\bf G}_m\otimes \Q$-version (\ref{QandZ}) of \ref{htp3}  by using horizontal morphisms $Y\to X(\C)$ from one-dimensional complex analytic manifolds $Y$ (see the proof of \ref{322}.1 in \ref{322})
and by using $\langle h(\alpha),  h^*(\beta)\rangle=(h(\alpha), h(\beta))$.

(2) follows from (1).    

\end{pf}

 \section{Height functions}\label{htf}
 
\subsection{The setting}\label{s2.1}

 \begin{sbpara}\label{2.1.2} We consider the following three situations (I), (II), (III).
 
 In Situation (I), we consider a number field $F_0\subset \C$.
 
 In Situation (II), we consider  a connected one-dimensional complex analytic manifold $B$ endowed with a finite flat morphism $\Pi: B\to \C$.
 
 In Situation  (III), we consider a connected projective smooth curve $C$ over $\C$.

\end{sbpara}

\begin{sbpara}\label{2.1.3}  

Let $G$, $w: {\bf G}_m\to G_{\text{red}}$ be as in \ref{Gandw}.  Let $\sG$ be a normal algebraic subgroup of $G$. Let $\sQ=G/\sG$.

In Situations (I), we assume that we are given $M_b: \Rep(G) \to MM(F_0)$ as in \ref{XFsetting2} and we assume \ref{XFsetting2}. Let $X(F)=X_{G,\sG,M_b,K}(F)$.

 In Situations (II) and (III), we assume we are given $H_b:\Rep(G)\to \Q$MHS  as in \ref{HAB} and we assume \ref{1.2.1} and \ref{1.2.1b}.  Let $X(\C)= X_{G,\sG,H_b,K}(\C)$ and let $\sM_{\hor}(Y, X(\C))$ ($Y=B$ in Situation (II) and $Y=C$ in Situation (III))  be as  in \ref{hormer}.

Let $K$ be an open compact subgroup of $\sG({\bf A}^f_K)$ which satisfies the neat condition (\ref{neat}).

\end{sbpara}

 \begin{sbpara}　
 Let 
$$\La=((V_i, w(i), s(i), c(i))_{1\leq i\leq m}, (t(w,d))_{w,d\in\Z, d\geq 1}),$$ where $m\geq 0$, $V_i\in \Rep(G)$, $w(i), s(i)\in \Z$, $c(i)\in \R$, and $t(w,d)\in \R$. 

\end{sbpara}

\begin{sbpara} The organization of Section 3 is as follows. 
 
 In Section \ref{s2.2}, we define  the height function 
 $$H_{\La} : X(F) \to \R_{>0}$$
in Situation  (I) for a finite extension $F$ of $F_0$ in $\C$  fixing somethings and assuming somethings, 
 and we also define the height function 
$$h_{\La} : \sM_{\text{hor}}(C,  X(\C)) \to \R$$
 in Situation  (III). 
 
In Section \ref{s2.3.2}, we review height functions in the usual Nevanlinna theory as a preparation for Sections \ref{s2.3}.  In Section \ref{s2.3}, we define the
 height function
 $$T_{f,\La}(r)\in \R \;\;(r\in \R_{\geq 0})\quad \text{for}\;\; f\in \sM_{\text{hor}}(B, X(\C))$$
 in Situation  (II).

  In Section \ref{s2.5}, we consider special cases 
  $h_{\spadesuit}$, $H_{\spadesuit}$, $T_{f, \spadesuit}(r)$ of 
   $h_{\La}$, $H_{\La}$, $T_{f, \La}(r)$, respectively, and then define
  a height function $h_{\heartsuit}: \sM_{\text{hor}}(C, X(\C))\to \R$ in Situation  (III),  and its analogues $H_{\heartsuit,S}$ for Situation  (I) and $N_{f,\heartsuit}(r)$ for Situation  (II), and we give some comments on height functions.  

  In Section \ref{s2.6}, we explain that the height functions of the different situations (I), (II), (III) are connected via asymptotic behaviors.  
 
   In Section \ref{s2.7}, we describe the relations of these height functions to the geometry of a toroidal partial compactification $\bar X(\C)$ of $X(\C)$.

\end{sbpara}

\begin{sbpara}\label{BC}
Situation (II) and Situation (III) are directly related when  $B=C\smallsetminus R$ with $B$ as in Situation (II), $C$ as in Situation (III),  and $R$ a finite subset of $C$. In this case, by the ''Great Picard Theorem'' applied to small neighborhoods of points of $R$ in $C$, $\Pi: B\to \C$  extends to a morphism  $C\to {\bf P}^1(\C)$ for which $R$ is the inverse image of $\infty\in {\bf P}^1(\C)$. 

\end{sbpara}

\subsection{Height functions $h_{\La}(\sH)$ and $H_{\La}(M)$}\label{s2.2}
 
 \begin{sbpara}
Assume we are in Situation (III).

We define the height function $h_{\La}:\sM_{\text{hor}}(C, X(\C))\to \R$ as 
 $$h_{\La}(\sH):= \sum_{i=1}^n c(i) \text{deg}(\gr^{s(i)}\gr^W_{w(i)}\sH(V_i)) + \sum_{w\in \Z, d\geq 1}  t(w,d)h_{\diamondsuit,w,d}(\sH(V_0)),$$
 where $h_{\diamondsuit, w,d}(\sH(V_0))$ is defined as in Part I, Section 1.6 by using the canonical polarization \ref{canpol1}. 
  
    \end{sbpara}

 \begin{sbpara} Let $F$ be a finite extension of $F_0$ in $\C$.

 To define the height function $H_{\La}: X(F)\to \R_{>0}$, we fix a $\Z$-lattice $V_{i,\Z}$ in $V_i$ ($1\leq i\leq m$) such that $\hat \Z\otimes_\Z V_{i,\Z}$ is stable under the action of $K$. 

\end{sbpara}

\begin{sbpara} 
Define
$$H_{\La}(M)=\prod_i H_{s(i)}(\gr^W_{w(i)}M(V_i))^{c(i)}\cdot \prod_{w,d} H_{\diamondsuit,w,d}(M(V_0))^{t(w,d)}$$
using Part I. 

Here in the definition of $H_{\diamondsuit, w,d}(M(V_0))$, we assume the conjectures in Part I, Section 1.7, which we assumed to define $H_{\diamondsuit, w,d}$ there. We use the canonical polarization on $\gr^W_wM(V_0)$ (\ref{canpol2})  for the definition of $H_{\diamondsuit, w,d}$. 

\end{sbpara}
\subsection{Reviews on height functions in Nevanlinna theory}\label{s2.3.2}

See \cite{LC} and \cite{Vo} for example.

\begin{sbpara}\label{238} Let $Y$ be a compact complex analytic manifold.

 Let $E$ be a divisor on $Y$. A Weil function of $E$ is an $\R$-valued continuous function on $Y\smallsetminus \text{Supp}(E)$ which is written locally on $Y$ as $$W=-\log(|g|)+\text{a continuous function},$$
where $g$ is a meromorphic function (found locally) such that $\text{div}(g)=-E$. 

A Weil function of $E$ exists.

\end{sbpara}

\begin{sbpara}\label{237} Let $B$ be a connected one-dimensional complex analytic manifold endowed with a finite flat morphism $\Pi:B\to \C$. Let $E$ be a divisor on $Y$ and assume $f(B)$ is not contained in $\text{Supp}(E)$. Assume that a Weil function $W$ of $E$ is given. Then the height function $T_{f,E}(r)$ $(r\in \R_{\geq 0})$
 for $f$ and $E$ in Nevanlinna theory  with respect to $W$ is defined by 
$$T_{f,E}(r):= m_{f,E}(r)+N_{f,E}(r)$$
where $m_{f,E}(r)$ and $N_{f,E}(r)$ are defined as follows. First, 
$$m_{f,E}(r):= \frac{1}{2\pi} \int_0^{2\pi} (\Pi_*f^*W)(re^{2\pi i \theta}) d\theta.$$
Here $\Pi_*$ is the trace map associated to $\Pi$. 
Writing $f^*E= \sum_{x\in B} n(x)x$ ($n(x)\in \Z$, note that this can be an infinite sum), let 
$$N_{f, E}(r) := \sum_{x\in B, 0<\Pi(x)|<r} n(x)\log(r/|\Pi(x)|)+\sum_{x\in B, \Pi(x)=0} n(x)\log(r).$$

If  $W'$ is another Weil function of $E$, then $W'-W$ is an $R$-valued  continuous function on the compact space $Y$, and hence 
there is a constant $c\in \R_{\geq 0}$ such that $|W-W'|\leq c$. Hence in this case, if we denote by $m'_{f,E}$ and $T'_{f,E}$ the $m_{f,E}$ and $T_{f,E}$ defined using $W'$ in place of $W$, respectively, we have $$|m_{f,E}(r)-m'_{f,E}(r)|\leq c, \quad |T_{f,E}(r)-T'_{f,E}(r)|\leq c.$$

\end{sbpara}

\begin{sbpara} In our Hodge-Nevanlinna theory, it is better to formulate the Nevanllina height functions by using the curvature forms of line bundles with $C^{\infty}$-metrics as below.

\end{sbpara}

\begin{sbpara}\label{236}

Let $L$ be a line bundle on $Y$. Then a $C^{\infty}$-metric on $L$ exists. 

Let $B$ be as in \ref{237} and let $f: B \to Y$ be a holomorphic map. Take a $C^{\infty}$-metric $|\;|_L$ on $L$ and 
define the height function of $(f, L)$ by
$$T_{f,L}(r)=T_{f,L |\;|_L}(r) :=(2\pi i)^{-1}  \int_0^r (\int_{B(t)} f^*\kappa(L)) \frac{dt}{t} \in \R\quad \text{for}\;r\in \R_{\geq 0}$$
where  $B(t)=\{x\in B\;|\; |\Pi(x)|<t\}$ and $\kappa(L)$ is the curvature form of $(L, |\;|_L)$.

\end{sbpara}

\begin{sbpara}\label{239} Let the notation be as in \ref{237}. 
Assume  $L=\sO_Y(E)$ for a divisor $E$.  Then $W:= - \log(|1|_L)$ is a Weil function of $E$.
The height functions in \ref{237} and \ref{236} (the former is defined by using $W$ and the latter is defined by using $|\;|_L$ ) are related as 
$$T_{f,L}(r)= T_{f,E}(r) -c$$
where $c$ is the constant $\lim_{z\to 0} (\Pi_*f^*W- n \log(|z|))$ with $n$ the coefficient of the divisor $-\Pi_*f^*E$ on $\C$ at $z=0$. 

This formula follows from Stokes' formula. 
\end{sbpara}

\subsection{Height functions $T_{f,\La}(r)$}\label{s2.3}

\begin{sbpara} Assume we are in Situation (II).  Let $f\in \sM_{\hor}(B, X(\C))$. 
We denote  $f$ also by $\sH$ when we regard it as the exact $\otimes$-functor $\Rep(G)\to \Q$VMHS$_{\log}(B)$ (\ref{hormer}).

\end{sbpara}

\begin{sbpara}\label{TLa}
Let $f\in \sM_{\text{hor}}(B,  X(\C))$. For $r\in \R_{>0}$. we define
$$T_{f, \La}(r):= T_{f, \La, \red}(r)+ \sum_{w,d} t(w,d)T_{f,\diamondsuit, w, d}(r),$$
where each term of the right hand side is defined in  \ref{Tka4}, \ref{defdia1}, \ref{defdia2} below. 
 
\end{sbpara}

 \begin{sbpara}\label{Tka4} 
 We define the first term $T_{f,\La,\red}(r)$ in the definition of $T_{f,\La}(r)$  in \ref{TLa}.

 Let $$\kappa_{X(\C), \La, \red}:= \sum_i c(i)\kappa(\text{det}(\gr^{s(i)}\gr^W_{w(i)}\sH_{X(\C)}(V_i)_{\sO}))_{\hor}$$ where $\kappa(\text{det}(\gr^{s(i)}\gr^W_{w(i)}\sH_{X(\C)}(V_i)_{\sO}))_{\hor}$ is the restriction to $T_{X(\C),\hor}$ of the curvature form of the line bundle $\text{det}(\gr^{s(i)}\gr^W_{w(i)} \sH_{X(\C)}(V_i)_{\sO})$ with the Hodge metric. The Hodge metric is defined by a polarization (\ref{XCpol}) but the curvature form  is independent of the choice of the polarization. This $\kappa_{X(\C),\La, \red}$ comes from an Hermitian form on $T_{X_{\red}(\C), \hor}$ via $T_{X(\C), \hor}\to T_{X_{\red}(\C), \hor}$. 
 
 Define 
 $$T_{f,\La,\red}(r):= \frac{1}{2\pi i}\int_0^r (\int_{B(t)} \kappa_{f,\La, \red}) \frac{dt}{t} \quad \text{with} \;\kappa_{f,\La,\red}:= f^*\kappa_{X(\C), \La,\red}.$$

 Here we have to be careful that not as in \ref{236}, the above differential $2$-form $\kappa_{f,\La,\red}$ is $C^{\infty}$ on $U=B\smallsetminus R$, where $R$ is the set of points at which $\sH(V_i)$ for some $i$ has singularity and it is a discrete subset of $B$, and $\kappa_{f,\La,\red}$  may have singularities at $R$. But  at a singular point,  $\kappa_{f,\La,\red}$ has the shape 
 $$O(\log(|z|)^{-2})d\log(z) \wedge d\log(\bar z)$$
 for a local coordinate function $z$, and hence 
 the integration on the right hand side in the definition of $T_{f,\La,\red}(r)$ converges. See \ref{Deg3}.

 \end{sbpara}

\begin{sbpara}\label{defdia1}

We define $T_{f, \diamondsuit, w,1}(r)$. 

Let $\kappa_{X(\C), \diamondsuit, w, 1}$ be the Hermitian form on $T_{X(\C), \hor}$ defined in \ref{XCdw1}. 
 Define
$$T_{f,\diamondsuit, w,1}(r):= \frac{1}{2\pi i}\int_0^r (\int_{B(t)} \kappa_{f,\diamondsuit,w,1})\frac{dt}{t}\quad \text{with}\;\;
\kappa_{f,\diamondsuit, w,1}:=f^* \kappa_{X(\C), \diamondsuit, w, 1}.$$

This $\kappa_{f,\diamondsuit, w,1}$ is $C^{\infty}$ on $B\smallsetminus R$ for some discrete subset $R$ of $B$. It may have singularities at $R$. But at a point of $R$ of singularity of the form
$$O(\log(|z|)^{-2}) d\log(z) \wedge d\log(\bar z).$$ Hence  the integration on the right hand side of the definition of $T_{f, \diamondsuit,w,1}(r)$ converges. See \ref{Deg4}. 
\end{sbpara}

\begin{sbpara}\label{defdia2}

Assume $d\geq 2$. 
We define 
$$T_{f,\diamondsuit,w,d}(r):= m_{f,\La}(r)+ N_{f, \diamondsuit,w,d}(r),$$
where: 
$$m_{f,\La}(r)=(2\pi)^{-1}\int_0^{2\pi} g(re^{\theta i}) d\theta
\quad \text{with} \; g=\Pi_*(\langle \delta_{w,d}, \delta_{w,d}\rangle^{1/d}),$$
where $\delta_{w,d}$ is that of 
$\sH(V_0)$ in Part I, Section 1.6, $\Pi_*$ is the trace map for $\Pi: B\to \C$,  
$$N_{f,\diamondsuit,w,d}(r):= \sum_{x\in B, 0<|\Pi(x)|<r}\langle N_{x,w,d}, N_{x,w,d}\rangle^{1/d}_{N_{x,0}}\cdot \log(r/|\Pi(x)|)+ \sum_{x\in B, \Pi(x)=0} \langle N_{x,w,d}, N_{x,w,d}\rangle_{N_{x,0}}^{1/d}\cdot \log(r),$$ 
where $N_{x,w,d}$ and $N_{x,0}$ are those of $\sH(V_0)$ in Part I, Section 1.6.

The integral $m_{f, \La}(r)$ exists by \ref{Deg7}.

\end{sbpara}

 \subsection{$\spadesuit$-height functions, $\heartsuit$-height functions, complements on height functions}\label{s2.5}
 
We consider height functions  $H_{\spadesuit}$ in Situation (I), $T_{f,\spadesuit}(r)$ in Situation (II), and $h_{\spadesuit}$ in situation (III). These are special cases of $H_{\La}$, $T_{f,\La}(r)$, $h_{\La}$, respectively.  We define height functions $H_{\heartsuit,S}$ in Situation (I), $N_{f,\heartsuit}(r)$ in Situation (II), $h_{\heartsuit}$ in Situation (III). We give  complements to height functions.

\begin{sbpara}\label{Spade} The following height functions $h_{\spadesuit}(\sH)$, $H_{\spadesuit}(M)$, $T_{f,\spadesuit}(r)$, which are special cases of the height functions $h_{\La}(\sH)$, $H_{\La}(M)$, $T_{f,\La}(r)$, will be important in Section 4:
$$h_{\spadesuit}(\sH):=-\text{deg}(\sH(\Lie(\sG))_{\sO}/\fil^0)\quad \text{in Situation  (III)},$$
$$H_{\spadesuit}(M):=\prod_{r<0} H_r(M(\Lie(\sG)))^{-1}\quad \text{in Situation  (I)},$$
$$T_{f,\spadesuit}(r):=\int_0^r (\int_{B(t)} f^*\kappa_{X(\C),\spadesuit})\frac{dt}{t} \quad \text{in Situation  (II)}. $$
Here $G$ acts on $\Lie(\sG)$ by the adjoint action. See \ref{kspade} for the definition of 
$\kappa_{X(\C), \spadesuit}$. 

This height function is the special case $\La=\spadesuit$ of the height functions in Section \ref{s2.2} and \ref{s2.3}, 
where $\La=\spadesuit$ means that $\La$ is such that $(w(i), s(i))$ ($1\leq i\leq m$) are all different pairs such that $\gr^{s(i)}\gr^W_{w(i)}H_b(\Lie(\sG))\neq 0$ and $s(i)<0$, $V_i= \Lie(\sG)$ for $1\leq i\leq m$, $c(i)=-1$ for all $i$, and $t(w,d)=0$ for all $w,d$.

 \end{sbpara}

\begin{sbpara} For a one-dimensional complex analytic manifold $Y$, for $\sH\in \sM_{\hor}(Y, X(\C))$, and for 
$x\in Y$, we define $e(x)=h_{\heartsuit, x}(\sH)\in \Z_{\geq 0}$ as follows.

Let $N'_x: \sH(V)_{\Q,x}\to \sH(V)_{\Q,x}$ be the local monodromy operator at $x$. It is the logarithm of the action of the canonical generator of the local monodromy group at $x$. Here $\sH(V)_{\Q,x}$ denotes the stalk of $\sH_\Q$ at a point $\neq x$ of $Y$ which is near to $x$. It belongs to $\sH(\Lie(\sG))_{\Q,x}$. 

By level structure, we have $N'_x\in {\bf A}^f_\Q\otimes_\Q \Lie(\sG)$ mod the adjoint action of $K$ on ${\bf A}^f_\Q\otimes_\Q \Lie(\sG)$. We define $e(x)=0$ if $N'_x=0$. if $N'_x\neq 0$, 
$e(x)$ is defined by  $\Z e(x)^{-1}=\{b\in \Q\;|\; \exp(bN'_x)\in K\}$. 
This is well defined because the local monodromy commutes with $N'_x$. 
  
\end{sbpara}
 
 \begin{sbpara} In Situation  (III), for $\sH\in \sM_{\hor}(C, X(\C))$, define
 $$h_{\heartsuit}(\sH):=\sum_{x\in C} h_{\heartsuit,x}(\sH)\in \Z_{\geq 0}.$$

 \end{sbpara}
 
 \begin{sbpara} We have
 $h_{\heartsuit}(\sH)\geq \sharp(\Sig(C,\sH))$, where $\Sig(C,\sH)$ denotes the set of points of $C$ at which $\sH$ has singularity. 
 \end{sbpara}

\begin{sbpara} In Situation  (II), for $f\in \sM_{\hor}(B, X(\C))$, define 
$$N_{f,\heartsuit}(r) := \sum_{x\in B, 0<|\Pi(x)|<r}  h_{\heartsuit, x}(\sH)\log(r/|\Pi(x)|))+\sum_{x\in B, \Pi(x)=0} h_{\heartsuit,x}(\sH) \log(r).$$

\end{sbpara}

\begin{sbpara}\label{N^1} Let 
$$N^{(1)}_{f, \heartsuit}(r):=\sum_{x\in \Sig(B, \sH), 0<|\Pi(x)|<r} \log(r/|\Pi(x)|) + \sum_{x\in \Sig(B, \sH), \Pi(x)=0} \log(r).$$
 Here $\Sig(B, \sH)$ is the set of points of $B$ at which $\sH$ has singularity. 

We have $N_{f,\heartsuit}(r) \geq N^{(1)}_{f,\heartsuit}(r)$.
\end{sbpara}

\begin{sbpara} We now consider Situation  (I). For a finite extension $F$ of $F_0$ in $\C$, for $M\in X(F)$ and for a non-archimedean place $v$ of $F$, we say $M$ is of good reduction at $v$ if $M(V)$ is of good reduction at $v$ for any $V\in \Rep(G)$. We say $M$ is of bad reduction at $v$ if it is not of good reduction at $v$. 

For each $M\in X(F)$, there are only finitely many non-archimedean places $v$ at which $M$ is of bad reduction. This is because for a faithful representation $V_1\in \Rep(G)$, by \ref{faith},  $M$ is good reduction at $v$ if and only if the mixed motive $M(V_1)$ is of good reduction at $v$. 

\end{sbpara}

\begin{sbpara}\label{Heart}

Let $F$ be a finite extension of $F_0$ in $\C$ and let $S$ be a finite set of places of $F$ which contains all archimedean places of $F$ and all non-archimedean places of $F$ at which $M_b$ is of bad reduction.
We define $H_{\heartsuit,S}(M)$ for $M\in X(F)$. The method is the same as that of 3.5.8 of \cite{Ka3}.

To define $H_{\heartsuit, S}:X(F)\to \R_{>0}$, we fix 
$$(L_i)_{1\leq i\leq m},\;\; (L'_i)_{1\leq i\leq m}, \;\;\text{and}\;\;  J$$
as in (i) below, and we assume that $K$ satisfies (ii) below.

\medskip

(i) $L_i$ and $L'_i$ are $\Z$-lattices in $V_i$ such that $L_i\supset L'_i$ and  $J$ is a subgroup of $\oplus_{i=1}^m (\Z/n_i\Z)^\times$ where $n_i$ is the smallest integer $>0$ which kills $L_i/L'_i$. 
  
  \medskip
  
(ii) $K$ coincides with the group of all elements  $g$ of $\sG({\bf A}^f_\Q)$ satisfying the following conditions (ii-1) and (ii-2).

(ii-1). $g\hat L_i=\hat L_i$, $g\hat L'_i=\hat L'_i$ for $1\leq i\leq m$.
Here $\hat L_i=\hat \Z\otimes_\Z L_i$ and $\hat L'_i$ is defined similarly.
 
 (ii-2) For some $h\in J$, the action of $g$ on $L_i/L_i'=\hat L_i/\hat L'_i$ coincides with the scaler action of $h$ via $J\to (\Z/n_i\Z)^\times$ for any $i$. 
   \medskip
   
   For $M=(M, \xi, \la)\in X(F)$, define
   $$H_{\heartsuit,S}(M):= \prod_{v\notin S} H_{\heartsuit,v}(M), \quad H_{\heartsuit,v}(M)=\sharp(\F_v)^{e(v)}$$
   where $v$ ranges over all places of $F$ outside $S$ and $e(v)\in \Q_{\geq 0}$ are defined as follows.

   $T_i= \tilde \la(L_i)$, $T'_i=\tilde \la(L_i')$, These are $\hat \Z$-lattices in $M(V_i)_{et}$ and are independent of the choice of the representative $\tilde \la$ of $\la$. 
    For a prime number $\ell$, let $T_{i,\Z_{\ell}}$ be the $\ell$-adic component of $T_i$ and define $T'_{i,\Z_{\ell}}$ similarly. 
 Let 
 $p=\text{char}(\F_v)$ ($\F_v$ is the residue field of $v$)  and let $D=D_{\pst, W(\bar \F_v)}$ (Part I, 1.3.7). Let $N'_v$ on $M(V_i)_{et, \Q_{\ell}}$ for $\ell\neq p$ and $N'_v$ on $D_{\text{pst}}(F_v, M(V_i)_{et,\Q_p})$ be the local monodromy operators at $v$.

   Let $I$ be the subgroup of $\Q$ consisting of all $b\in \Q$ satisfying the following condition:

For any $\ell\neq p$, the action of $\exp(bN_v')$ on $M(V_i)_{et, \Q_{\ell}}$ satisfies
$\exp(bN'_v)T_{i,\Z_{\ell}}=T_{i,\Z_{\ell}}$ and $\exp(bN_v')T'_{i,\Z_{\ell}}= T'_{i, \Z_{\ell}}$ for $1\leq i\leq m$, 
the action of $\exp(bN_v')$  on $\hat F_{v,0,\text{ur}}\otimes_{F_{v,0, \text{ur}}} D_{\text{pst}}(F_v,M(V_i)_{et, \Q_p})$ satisfies $\exp(bN_v')D(T_{i, \Z_p})= D(T_{i, \Z_p})$  and $\exp(bN'_v)D(T'_{i,\Z_p})= D(T'_{i, \Z_p})$ for $1\leq i\leq m$, and there is $h\in J$ such that for $1\leq i\leq m$,  the action of $\exp(bN_v')$ on $T_{i, \Z_{\ell}}/T'_{i, \Z_{\ell}}$ 
coincides with the scaler action of $h$ via $J\to (\Z/n_i\Z)^\times$ for any $\ell\neq p$, 
 and the action of $\exp(bN_v')$ on 
 $D(T_{i,\Z_p})/D(T'_{i,\Z_p})$ coincides with the scaler action of $h$ via $J\to (\Z/n_i\Z)^\times$.

Then $I$ is non-zero. $I=\Q$ if $M$ is of good reduction at $v$.

We define $e(v)=0$ if $I$ is not isomorphic to $\Z$, and we define $e(v)\in \Q_{>0}$ by $I=e(v)^{-1}\Z$ if $I\cong \Z$. 

We expect that the following (*) is true. 

\medskip

(*)  $e(v)\in \Z$. $e(v)=0$ if and only if $M$ is of good reduction at $v$. 

\medskip

If this (*) is true, we have 
 $$\prod_{v\in \Sig_S(M)} \sharp(\F_v) \leq H_{\heartsuit,S}(M)$$
where $\Sig_S(M)$ is the set of all places of $F$ which do not belong to $S$ and at which $M$ has bad reductions. 

\end{sbpara}

\begin{sbpara} 
We give comments on height functions.

(1) 
The bijections in \ref{XFtoXF} (1) do not change the height functions in (I), (II), (III). 

(2) In this Section 3, we considered the height functions on $X_{G,\sG, M_b, K}(F)$ but not yet on  $X_{G,\Upsilon,K}(F)$.  We define height functions on $X_{G,\Upsilon,K}(F)$  by using the bijection $X_{G, \Upsilon, K}(F)\cong X_{G, G, M_1, K'}(F)$ in \ref{XFtoXF} (2) by choosing $M_1\in X_{G,\Upsilon,K}(F)$, and using the height functions on the latter. Then they are independent of the choice of $M_1$. 
\end{sbpara}

\begin{sbprop}\label{semi+}

Assume $G$ is reductive and assume $\kappa_{X(\C),\La, \red}$ (\ref{Tka4}) is positive semi-definite. Then 

(1) In Situation  (II), $T_{f\in \La}(r) \geq 0$ for any $f\in \sM_{\hor}(B, X(\C))$ and any $r\in \R_{\geq 0}$.

(2) In Situation  (III), $h_{\La}(\sH)\geq 0$ for any $\sH\in \sM_{\hor}(C, X(\C))$.

\end{sbprop}
\begin{pf} (1) is evident. (2) follows from $(2\pi i)^{-1}\int_{C} \kappa_{\La}(\sH)= h_{\La}(\sH)$ (\cite{Pe}). \end{pf}

\begin{sbpara}\label{ample}  
We say $\La$ is ample if the Hermitian form $\kappa_{X(\C), \La,\red}$ on $T_{X_{\text{red}}(\C), \hor}$ (\ref{Tka4}) is positive definite and 
$t(w,d)>0$ for any $w,d$. 

For example, by the theorem of Griffiths in \ref{321}, $\La$ is ample if $V_i=V_0$ for all $i$, $c(i)= s(i)$ for all $i$, $\{(w(i), s(i))\; |\;i\in \Z\}$ covers all pairs $(w,s)$ such that $\gr^s\gr^W_wV_0\neq 0$, and $t(w,d)>0$ for all $w,d$. 
\end{sbpara}

\begin{sbprop}\label{HN+} Assume we are in Situation (II). Assume $\La$ is ample and $G$ is reductive.

Let $f \in \sM_{\text{hor}}(B, X(\C))$ and assume that $T_{f,\La}(r)=o(\log(r))$.  Then  $f$ is constant.

\end{sbprop}

\begin{pf} If $f$ is not constant,  the Hermitian form $\kappa_{f, \La,\red}$ is positive definite and hence  there are constants $b,c>0$ such that $(2\pi i)^{-1}\int_{B(t)} \kappa_{f, \La,\red}\geq c$ for any $t\geq b$. For any $r\geq b$, we  have $T_{f,\La}(r) \geq \int_b^r c dt/t = c(\log(r)-\log(b))$.
\end{pf}
The author expects that the assumption $G$ is reductive in \ref{HN+} is unnecessary.

\begin{sbprop}\label{H+}  Assume we are in Situation (III). Assume $\La$ is ample and let $\sH\in \sM_{\hor}(C, X(\C))$. Then $h_{\La}(\sH)\geq 0$.  If $h_{\La}(\sH)= 0$, then $\sH$ is constant. 

\end{sbprop}

\begin{pf} This is proved in the same way as Part I, Proposition 1.6.16. 
\end{pf}

The proof of the following Proposition is easy.

\begin{sbprop} (1) Assume we are in Situation (I) and let $F'$ be a finite extension of $F$ of degree $m$. Let $M\in X(F)$ and let $M'$ be the image of $M$ in $X(F')$. Then we have 
$$H_{\La}(M')= H_{\La}(M)^{[F':F]}, \quad H_{\heartsuit,S'}(M')=H_{\heartsuit,S}(M)^{[F':F]}.$$
Here when we consider $H_{\heartsuit,S}$, we assume what we assumed in Section \ref{s2.5} to define it, and $S'$ denotes the set of all places of $F'$ lying over $S$. 

(2) Assume we are in Situation (II) and let $B'$ be a connected smooth curve over $\C$ endowed with a finite flat morphism $B'\to B$ of degree $[B':B]$. Let $f\in \sM_{\hor}(B, X(\C))$ and let $f'\in \sM_{\hor}(B', X(\C))$ be the composition of $f$ and $B'\to B$. Then 
$$T_{f'. \La}(r)= [B':B]\cdot T_{f,\La}(r), \quad N_{f',\heartsuit}(r)= [B':B]\cdot N_{f,\heartsuit}(r).$$

(3) Assume we are in Situation (III) and let $C'$ be a proper smooth curve over $\C$ endowed with a finite flat morphism $C'\to C$ of degree $[C':C]$. Let $\sH\in \sM_{\hor}(C, X(\C))$ and let $\sH'\in \sM_{\hor}(C', X(\C))$ be the pull back of $\sH$. Then 
$$h_{\La}(\sH')= [C':C]\cdot h_{\La}(sH), \quad h_{\heartsuit}(\sH')= [C':C]\cdot h_{\heartsuit}(\sH).$$

\end{sbprop}

\subsection{Degeneration}\label{ss2.5}

\begin{sbpara}\label{Deg1} Concerning \ref{Deg1}--\ref{Deg3}, see \cite{Pe}.   

Let $Y$ be a one-dimensional complex manifold, let $R$ be a discrete subset of $Y$, let $L$ be a line bundle on $Y$, and let $|\;|$ be a $C^{\infty}$ metric on the restriction of $L$ to $Y\smallsetminus R$. We say the metric $|\;|$ is good at $R$ if the following condition (*) is satisfied at each $x\in R$. 
Let $z$ be a local coordinate function on $Y$ at $x$ such that $z(x)=0$. 

(*) There are an open neighborhood $U$ of $x$ on which $|z|<1$, a basis $e$ of $L$ on $U$, and  constants $C_1,C_2, C_3, C_4, c\in \R_{> 0}$ such that we have the following (i)--(iii) on $U\cap (Y\smallsetminus  R)$.

(i) $C_1\log(1/|z|)^{-c} \leq |e| \leq C_2\log(1/|z|)^c$.  

(ii) $|z\frac{\partial}{\partial z}\log|e| |\leq C_3 \log(1/|z|)^{-1}$. 

(iii) $|z^2 \frac{\partial}{\partial z}\frac{\partial}{\partial\bar z}  \log(|e|)| \leq C_4 \log(1/|z|)^{-2}$. 

\end{sbpara}

\begin{sbpara}\label{Deg2} In \ref{Deg1}, consider the case $Y$ 
is a connected projective smooth curve $C$ over $\C$. Assume and the metric $|\;|$ is good at $R$. Let $\kappa(L)$ be the associated curvature form on $C\smallsetminus R$. Then by \cite{Mum}, 
$\kappa(L)$ is integrable on $C$, and $(2 \pi i)^{-1}\int_C \kappa(L) = \text{deg}(L)$.

\end{sbpara}

\begin{sbpara}\label{Deg3} Let $Y$ and $R$ be as in \ref{Deg1}, and let $\sH\in \Q$VMHS$_{\log}(Y)\cap \Q$VMHS$(Y\smallsetminus R)$. Assume $\sH$ is pure and polarizable. Then 
the Hodge metric on $\text{det}(\gr^r \sH_{\sO})$ given by a polarization of $\sH$ is good at $R$  by \cite{Pe}.

\end{sbpara}

\begin{sbpara}\label{Deg4} Let  $Y$ and $R$ be as in \ref{Deg1}. Assume we are given  $\sH_0$, $\sH_1$, $\sH_2$ on $Y\smallsetminus R$ as in Section 2.4 (we replace $Y$ there by $Y\smallsetminus R$ here). Assume $\sH_0$ is polarizable and assume that $\sH_i\in \Q VMHS_{\log}(Y)$ for $i=0,1,2$. Then by \cite{BP},  
the  ${\bf G}_m\otimes \Q$-torsor $L(\sH_1, \sH_2)$ on $Y\smallsetminus R$ (Section 2.4) extends uniquely to a ${\bf G}_m\otimes \Q$-torsor on $Y$ and its metric on $Y\smallsetminus R$ is good at $R$.

\end{sbpara}

\begin{sbpara}\label{Deg5}  In \ref{Deg4}, assume that $Y$ is a connected projective smooth curve over $\C$.  Then by \cite{BP}, 
$\text{deg}(L(\sH_1,\sH_2))\in \Q$ is equal to the height paring $\text{ht}(\sH_1, \sH_2)$. By \ref{Deg2} and \ref{Deg4}, this $\text{ht}(\sH_1, \sH_2)$ coincides with  $(2\pi i)^{-1} \int_C \kappa(L(\sH_1, \sH_2))$.

Hence the positive definite property of $\kappa(L(\sH_1, \sH_2))$ (\ref{htp3}) gives a new proof of the positive definite property of the height pairing $\text{ht}(\sH_1, \sH_2)$. 

\end{sbpara}

\begin{sbpara}\label{Deg6} We give a property  $\delta_{w,d}$ in degeneration. 

Let $\Delta= \{q\in \C\;|\; |q|<1\}$ and let $\Delta^*= \Delta \smallsetminus \{0\}$. 
Then the upper half plane $\frak H$ is the universal covering of $\Delta^*$ with the covering map $x+iy \mapsto q= e^{2\pi i(x+iy)}$. Endow $\Delta$ with the log structure associated to $\{0\}\subset \Delta$ and let $\Delta^{\log}$ be the space 
over $\Delta$ in \cite{KN} (see also \cite{KU}, \cite{KNU}). We have $\Delta^{\log}= |\Delta|\times {\bf S}^1$ 
where $|\Delta|:=\{r\in \R\;|\; 0\leq r<1\}$ and ${\bf S}^1:= \{u\in \C^\times\;|\;|u|=1\}$, and
the map $\Delta^{\log}\to \Delta$ is given by $(r,u)\mapsto ru$. Let $\tilde \Delta^{\log}= |\Delta|\times \R$ be the universal covering of $\Delta^{\log}$ with the covering map $\tilde \Delta^{\log}\to \Delta^{\log}$  
induced by $\R\to {\bf S}^1\;;\;x\mapsto e^{2\pi i x}$. Then $\frak H$ is identified with an open set of $\tilde \Delta^{\log}$ via the embedding $x+iy\mapsto (|q|, x)=(e^{-2\pi y},x)$.

Assume we are given  $\sH\in \Q$VMHS$_{\log}(\Delta) \cap \Q$VMHS$(\Delta^*)$ and assume that $\gr^W_w\sH$ are polarized for all $w\in \Z$. We consider the behavior of $(\delta_{w,d}(\sH(q)), \delta_{w,d}(\sH(q)))$ when $q\in \Delta^*$ converges to $0$, where $\sH(q)$ denotes the fiber of $\sH$ at $q$ and $(\;,\;)$ denotes the Hodge metric of $(\gr^W_w\sH)^*\otimes \gr^W_{w-d}\sH$.

Fix an element $\tilde 0$ if $\tilde \Delta^{\log}$ lying over $0\in \Delta$. The local system $\sH_{\Q}$ on $\Delta^*$ extends uniquely to a local system on $\Delta^{\log}$
 which we denote also by $\sH_{\Q}$. Let $H_{0,\Q}$ be the stalk 
 of $\sH_{\Q}$ at the image of $\tilde 0$ in $\Delta^{\log}$. 
 Then $H_{0,\Q}$ has the weight filtration and the polarization $\gr^W_w \sH_\Q \times \gr^W_w \sH_\Q \to \Q \cdot (2\pi i)^{-w}$ induces a $\Q$-bilinear form 
 $\langle\;,\;\rangle_w$ on
 $\gr^W_w H_{0,\Q}$.  Let $h(w,r)$ be the rank of $\gr^r \gr^W_w\sH_{\sO}$.

  Let $D$ be the set of all 
 descending filtrations $\fil$ on $H_{0,\C}$ such that 
 $(H_{0,\Q}, W, (\langle\;,\;\rangle_w)_w, \fil)$ is a mixed Hodge structure 
  with polarized pure graded quotients and such that 
  $\dim_\C(\gr^r\gr^W_w\fil)=h(w,r)$ for all $w,r$.

We have the period map $p: \frak H \to D$ associated to $\sH$. The theory of associated SL$(2)$-orbit in \cite{KNU} Part II, Part III gives the following $((\rho_w)_{w\in \Z}, s, {\bf r})$.  $\rho_w$ 
for each $w\in \Z$ is a  homomorphism 
$\text{SL}(2)_\R\to \text{Aut}_\R(\gr^W_wH_{0,\R}, \langle \;,\;\rangle_w)$ of algebraic groups over $\R$. 
The canonical $\R$-splitting of $W$ (\cite{KNU} part I, \S4, Part II, \S1.2) associated to $p(\tilde q)$ converges 
when $\tilde q\in \frak H$ converges to $\tilde 0$. 
We denote the limit splitting of $W$ over $\R$ by $s$. 
Let $\tau: {\bf G}_{m,\R}\to \text{Aut}_\R(H_{0,\R},W)$ be the 
homomorphism which corresponds via $s$ to the direct sum of 
$t^w \rho_w\begin{pmatrix} 1/t  & 0\\ 0&t\end{pmatrix}$ on $\gr^W_wH_{0,\R}$ (\cite{KNU} Part II, 2.3.5). $\tau(y^{1/2})p(\tilde q)\in D$ converges in $D$ 
when $\tilde q\in \frak H$ converges to $\tilde 0$. We denote the limit by ${\bf r}\in D$. (This ${\bf r}$ is written as $\exp(iN)\hat F$ in \cite{KNU} Part IV, Theorem 6.2.4.)

\end{sbpara}
\begin{sbprop}\label{Deg7} 
 When $q\in \Delta$ converges to $0$,  $y^{-d}(\delta_{w,d}(\sH(q)), \delta_{w,d}\sH(q)))$ converges to $(\delta_{w,d}({\bf r}), \delta_{w,d}({\bf r}))$.

\end{sbprop}

Here $(\;,\;)$ is the Hodge metric.

\begin{pf} 
We have 
$$y^{-d}(\delta_{w,d}(\sH(q)), \delta_{w,d}(\sH(q)))= (\tau(y^{1/2})\delta_{w,d}(p(q)), \tau(y^{1/2})\delta_{w,d}(p(q)))$$ $$=(\delta_{w,d}(\tau(y^{1/2})p(q)), \delta_{w,d}(\tau(y^{1/2})p(q))\to (\delta_{w,d}({\bf r}), \delta_{w,d}({\bf r})).$$ 
\end{pf}

\subsection{Asymptotic behaviors}\label{s2.6}

In the following \ref{Tandh} and \ref{Tandh2}, we describe relations between the height functions in  Situation  (II) and the  height functions in Situation  (III). At the end of this Section \ref{s2.6}, we shortly describe the relations of between the height functions in Situation (I) and the height functions in (II) or (III), but the details of them will be given elsewhere.  

\begin{sbprop}\label{Tandh}  Assume $B=C\smallsetminus R$ for a finite subset $R$ of $C$ (that is, we are in the situation of \ref{BC}).  Let $\sH\in \sM_{\hor}(C, X(\C))$ and let $f\in \sM_{\hor}(B, X(\C))$ be the restriction of $\sH$ to $B$. 

(1)  $T_{f,\La}(r)/\log(r)$ converges  when $r\to \infty$. 

(2) If either one of the following conditions (i)--(iii) is satisfied, $T_{f,\La}(r)/\log(r)$ converges to $h_{\La}(\sH)$. 

(i)  $G$ is reductive.

(ii) $t(w,d)=0$ for any $d\geq 2$.

(iii) $\sH$ has no degeneration at $R$. 
\end{sbprop}

\begin{pf} (1) follows from \ref{Deg7}. 

We prove (2). Case (i) follows from Case (ii). Assume (ii). Then $T_{f,\La}(r)= \int_0^r (\int_{B(t)} \omega)dt/t$ with $\omega= \kappa_{f,\La,\red}+\sum_w t(w,1)\kappa_{f,\diamondsuit,w,1}$. 
Let $A(t)= (2\pi i)^{-1}\int_{B(t)} \kappa$. Then when $t\to \infty$, $A(t)$ converges to 
$A:=(2\pi i)^{-1}\int_C \kappa=h_{\La}(\sH).$ 
Take $\epsilon>0$. Then there is $r_0>0$ such that $|A(r)- A|\leq \epsilon$ if $r\geq r_0$. Let $c:= \int_0^{r_0} A(t)dt/t$. We have
$T_{f, \La}(r)=c+\int_{r_0}^r A(t) dt/t= c + \int_{r_0}^r  (A + u(t))dt/t$ with $|u(t)|\leq \epsilon$, and hence 
$|T_{f,\La}(r)- c - A \log(r/r_0)|\leq \epsilon \int_{r_0}^r dt/t= \epsilon \log(r/r_0)$. Hence if $r$ is sufficiently large, we have $|T_{f,\La}(r)/\log(r) - A|\leq 2 \epsilon$. 

Assume  (iii). In this case,  $m_{f,\diamondsuit, w,d}$ is bounded when $r\to \infty$ and hence $$T_{f, \diamondsuit, w,d}(r)/\log(r) \to 
\sum_{x\in B} \langle N_{x,w,d}, N_{x,w,d}\rangle^{1/d}_{N_{x,0}}=h_{\diamondsuit, w,d}(\sH).$$
\end{pf}

\begin{sbprop}\label{Tandh2} Let the assumption be as in \ref{Tandh}. Let $\sH\in \sM_{\hor}(C, X(\C))$.

(1) $\lim_{r\to \infty} N_{f, \heartsuit}(r)/\log(r)= \sum_{x\in B} h_{\heartsuit,x}(\sH).$

(2) $\lim_{r\to \infty} N^{(1)}_{f, \heartsuit}(r)/\log(r) = \sharp\{x\in B\;|\;\sH \;\text{has singularity at $x$}\}.$

\end{sbprop}

\begin{pf} This is clear. \end{pf} 
 
 \begin{sbrem} (1) In \ref{Tandh2}, if  $\sH$ has no degeneration at $R$, the rhs becomes $h_{\heartsuit}(\sH)$ in (1) $\sharp\{x\in C\;|\; \sH \;\text{has singularity at $x$}\}$ in (2).

 (2)
 \ref{Tandh} and \ref{Tandh2} tell that the study of height functions in Situation  (III) is essentially reduced to that in Situation  (II).

\end{sbrem}

\begin{sbpara}\label{Handh} In the above, we considered relations between Situation (II) and Situation (III). 

We have relations between Situation (I) and Situation (III), and also relations between Situation (I) and Situation (II). We give only rough stories here. We hope to discuss more precise things elsewhere. 

Let $C_0$ be a projective smooth curve over a number field $F_0$.  We can define the set $X(F_0(C_0))$, where $F_0(C_0)$ denotes the function field of $C_0$, in the similar way as the definition of $X(F)$.  For $M\in X(F_0(C_0))$, for a finite extension $F$ of $F_0$ in $\C$ and for $x\in C_0(F)$ at which $M$ does not have  singularity, we can define the specialization $M(x)\in X(F)$ of $M$ at $x$. 

Concerning the relation between (I) and (III), we can show that
$$ \log(H_{\La}(M(x))/\log(H(x))\;\;\text{in (I)}\quad  \text{and} \quad h_{\La}(M_H)\;\;\text{in (III)}$$ are closely related. Here $M_H$ is the Hodge realization of $M$ on $C:=C_0(\C)$ and $H(x)$ is the height of $x$ as an element of $C_0(F)$ defined by a ${\bf G}_m\otimes \Q$-torsor on $C_0$ of degree $1$. In some cases, the right hand side here is the limit of the left hand side. 

Concerning the relation between (I) and (II), we can show that for an affine dense open set $B_0$ of $C_0$, if we denote by $f$ the element of $\sM_{\hor}(B, X(\C))$ ($B:=B_0(\C)$) induced by $M$, $$
\log(H_{\La}(M(x)))\;\;\text{in (I)}\quad \text{and}\quad T_{f, \La}(H(x))\;\;\text{in (II)}$$ are closely related.

\end{sbpara}
 
 \subsection{Toroidal partial compactifications  and height functions}\label{s2.7}

Here we describe the toroidal partial compactification $\bar X(\C)$ of $X(\C)$ and its relation to height functions.
 
This $\bar X(\C)$ is a generalization of $\bar X(\C)$ in Part I. It is given by the work of Kerr and Pearlstein \cite{KP} in the case $G$ is reductive and $\sG=G$. For general $G$, it is described in \cite{KNU2} shortly in the case $\sG=G$. The details of the general case will be given in a forth-coming paper \cite{KNU}, Part V. 

In the rest of this paper,  this Section \ref{s2.7} serves to make our philosophy clearer, and we do not use the contents of this Section \ref{s2.7} except the following points: \ref{barX3} is used in Remark \ref{Qvo}.9. \ref{amplog} is used to interpret the result \ref{spadepos}. 
. 

\begin{sbpara}\label{Sig}

Let $\Sig$ be the set of cones in $\Lie(\sG)_\R$ of the form $\R_{\geq 0}N$ for some element $N$ of $\Lie(\sG)$ satisfying the following conditions (a) and (b) for any $V\in \Rep(G)$. 

\medskip

(a) The image of $N$ under $\Lie(G)\to \text{End}(V)$ is a nilpotent operator on $V$. 

(b) There is a relative monodromy filtration of $N: V\to V$ with respect the weight filtration $W_{\bullet}V$. 

\medskip

If $V_1\in \Rep(G)$ is a faithful representation, (a) and (b) are satisfied for any $V\in \Rep(G)$ if they are satisfied by $V=V_1$ (\ref{faith}).

\end{sbpara}

\begin{sbpara} Let $D_{\Sig}$ be the set of pairs $(\sig, Z)$, where $\sig\in \Sig$ and $Z$ is a non-empty subset of $\check{D}(G, \Upsilon)$ satisfying the following (i)--(iv).

\medskip

(i)  Write $\sig=\R_{\geq 0}N$ with $N\in \Lie(\sG)$. Then $Z$ is an $\exp(\C N)$-orbit in $\check{D}(G, \Upsilon)$.

(ii) The image of $Z$ in $\check{D}(\sQ, \Upsilon_{\sQ})$ is $\text{class}(H_{b,\sQ})\in D(\sQ, \Upsilon_{\sQ})$, 

(iii) If $H\in Z$, $N$ belongs to $\fil^{-1}H(\Lie(\sG))_\C$,  where $G$ acts on $\Lie(\sG)$ by the adjoint action and $\fil^{-1}$ is the Hodge filtration. 

(iv) Let $H\in Z$. Then $\exp(zN)H\in D(G, \Upsilon)$ if $z\in \C$ and $\text{Im}(z)$ is sufficiently large. 

\medskip
 
Let $\bar X(\C):= \sG(\Q) \bs (D_{\Sig} \times (\sG({\bf A}^f_\Q)/K))$ 

\end{sbpara}

\begin{sbpara}

$\bar X(\C)$ has a structure of  a logarithmic manifold 
 which extends the complex analytic structure of $X(\C)$.
(Logarithmic manifold is a generalization of complex analytic manifold (\cite{KU}, 3.5.7). It is like a complex manifold with slits.)

\end{sbpara}

 \begin{sbpara} For a one-dimensional complex analytic manifold $Y$, for a discrete subset $R$ of $Y$,  and for a horizontal holomorphic map $f:Y\smallsetminus R\to X(\C)$, the following (i)--(iii) are equivalent. (i) 
 $f$ is meromorphic on $Y$ (\ref{hormer}). (ii)  $f$ extends to a morphism $Y\to \bar X(\C)$ of locally ringed spaces over $\C$. (iii) $f$ extends to a morphism $Y\to \bar X(\C)$ of logarithmic manifolds. If these equivalent conditions are satisfied, the extensions of $f$ in (ii) and (iii) are unique.

\end{sbpara}

\begin{sbpara} We describe the relations of $h_{\spadesuit}$ and $h_{\heartsuit}$ to the extended period domain $\bar X(\C)$.

Around here, $Y$ is a one-dimensional complex analytic manifold and $f\in \sM_{\hor}(Y, X(\C))$ and let $\sH$ be the corresponding exact $\otimes$-functor $\Rep(G)\to \Q$VMHS$_{\log}(Y)$. We denote the morphism $Y\to \bar X(\C)$ which extends $f$ by the same letter $f$ 
\end{sbpara}

\begin{sbpara}\label{barX1} For $V\in \Rep(G)$, we have a universal log mixed Hodge structure $\sH_{\bar X(\C)}(V)$ on $\bar X(\C)$.

 We have 
$$\sH(V)= f^* \sH_{\bar X(\C)}(V).$$

\end{sbpara}

\begin{sbpara}\label{barX20} We have the sheaf $\Omega^1_{\bar X(\C)}(\log)$ of differential forms on $\bar X(\C)$ with log poles outside $X(\C)$ and the sheaf $\Omega^1_{\bar X(\C)}\subset \Omega^1_{\bar X(\C)}(\log)$ of differential forms on $\bar X(\C)$. These are vector bundles on $\bar X(\C)$. 

\end{sbpara}

\begin{sbpara}\label{barX4}
Let $I_{\bar X(\C)}$ be the invertible ideal of $\sO_{\bar X(\C)}$ which is locally generated by generators of log structures. For $x\in \bar X(\C)$, the stalk $I_x$ coincides with $\sO_{\bar X(\C),x}$ if and only if $x\in X(\C)$.

We have an exact sequence
$$0 \to \Omega^1_{\bar X(\C)}\to \Omega^1_{\bar X(\C)}(\log)\to \sO_{\bar X(\C)}/I_{\bar X(\C)}\to 0,$$
where the central arrow sends $d\log(q)$ for a local generator $q$  of the log structure to $1$.

\end{sbpara}

\begin{sbpara}\label{barX2} The vector bundle $\sH_{\bar X(\C)}(\Lie(\sG))_{\sO}/\fil^0$ on $\bar X(\C)$  is canonically isomorphic to the logarithmic tangent bundle of $\bar X(\C)$, that is, it is canonically isomorphic to the dual vector bundle of $\Omega^1_{\bar X(\C)}(\log)$. 

\end{sbpara}

\begin{sbpara}\label{barX3} In Situation (III), for $Y=C$, 
by \ref{barX1} and \ref{barX2}, we have
$$h_{\spadesuit}(\sH)= \text{deg}(f^*\Omega^1_{\bar X(\C)}(\log)).$$

\end{sbpara}

\begin{sbpara}\label{barX5} In Situation (II) (resp. (III)), for $Y=B$ (resp. $Y=C$)  and for $x\in Y$, 
$h_{\heartsuit,x}(\sH)= e(x)$, where $e(x)$ is the integer $\geq 0$ such that $f^*I_{\bar X(\C)}=m_{Y,x}^{e(x)}$ in $\sO_{Y,x}$.

\end{sbpara}

\begin{sbpara}\label{barX6}
Assume $Y=C$. 

By \ref{barX5}, we have
$$h_{\heartsuit}(\sH)=- \text{deg}(f^*I_{\bar X(\C)}).$$

By this and by \ref{barX3}, we have
$$h_{\spadesuit}(\sH)-h_{\heartsuit}(\sH)= \text{deg}(f^*\Omega^1_{X(\C)}).$$
\end{sbpara}

 \begin{sbpara}

In the comparison of (I) and (1) in Section 0, $\bar X(\C)\supset X(\C)$ is like $\bar V \supset V$ such that $\bar V$ is smooth and $E:= \bar V\smallsetminus V$ is a divisor on $\bar V$ with normal crossings. Basing on \ref{barX3} and following the analogy between (I) and (III), we think that the height function 
$H_{\spadesuit}$ for (I)  is similar to the height function $H_{K+E}$ for (1), where $K$ is the canonical divisor of $\bar V$. 
\end{sbpara}

\begin{sbpara}\label{amplog}  It can be shown that the Hodge metric on $\text{det}(\Omega^1_{X(\C)})$ 
 (defined by a polarization on $\sH_{X(\C)}(\Lie(\sG))$ (\ref{XCpol}) and the duality  
 between $\Omega^1_{X(\C)}$ and $T_{X(\C)}= \sH_{X(\C)}(\Lie(\sG))/\fil^0$) 
 extends to a metric on $\text{det}(\Omega^1_{\bar X(\C)}(\log))$ with at worst log 
 singularity. Hence the positivity of the curvature form $\kappa(\text{det}(\Omega^1_{X(\C)}))_{\hor}$ (\ref{spadepos}) tells that $\text{det}(\Omega^1_{\bar X(\C)}(\log))$ is something like an ample line bundle.

\end{sbpara}

\section{Speculations}
We extend speculations in Part I to the setting of this Part II.

Let the setting be as in \ref{2.1.2}--\ref{2.1.3}.

\subsection{Speculations on positivity}\label{ss:fin}

\begin{sbpara}\label{Q1} {\bf Question 1.} Assume $G$ is reductive. 

Are the following (i)--(v) equivalent?

\medskip

(i) The Hermitian form  $\kappa_{\X(\C), \La,\red}$ on $T_{X(\C),\hor}$   is positive semi-definite.

(ii) In Situation (III), we  have $h_{\La}(\sH)\geq 0$ in the situation (III). for any  $C$  and for any $\sH \in \sM_{\text{hor}}(C, X(\C))$.

(iii) In Situation (II), we have  $T_{f,\La}(r) \geq 0$ for any $B$, any $f\in \sM_{\text{hor}}(B, X(\C))$ and any $r\in \R_{\geq 0}$. 

(iv) In Situation (I), if $F$ is a finite extension of $F_0$ in $\C$, there is $c>\R_{>0}$ such that $H_{\La}(M)\geq c$ for any $M\in X(F)$.

(v) In Situation (I), there is $c\in \R_{>0}$ such that for any number field $F\subset \C$ such that $F_0\subset F$ and for any $M\in X(F)$, we have $H(M)\geq c^{[F:F_0]}$.

\medskip

Remark \ref{Q1}.1.  We have 
(i) $\Rightarrow$ (iii) $\Rightarrow$ (ii). In fact, (i) $\Rightarrow $ (ii) and (i) $\Rightarrow $ (iii) follow from \ref{semi+}, (iii) $\Rightarrow$ (ii) follows from \ref{Tandh}.

\end{sbpara}

\begin{sbpara}\label{Q2}

{\bf Question 2.} Assume $G$ is reductive. 

 Are the following (i)--(v) equivalent?
 
 \medskip

(i) $\La$ is ample (\ref{ample}).

(ii)  In Situation (III),  if $\sH\in \sM_{\text{hor}}(C, X(\C))$ and if $h_{\La}(\sH)\leq 0$, then $\sH$ is constant.

(iii) In situation (II), if $f\in \sM_{\text{hor}}(B, X(\C))$ and if there is $c\in \R$ such that $T_{f, \La}(r)\leq c$ for any $r\in \R_{\geq 0}$, then $f$ is constant.

(iv) In Situation (I), for any finite extension $F$ of $F_0$ in $\C$ and for any $c\in \R_{>0}$, the set $\{m\in X(F)\;|\; H_{\La}(M) \leq c\}$ is finite.

(v) In Situation (I), for any $d\geq 1$ and $c\in \R_{>0}$, there are finitely many pairs $(F_i, M_i)$ ($1\leq i\leq n$), where $F_i$ are finite extensions of $F_0$ in $\C$  and $M_i\in X(F_i)$, satisfying the following condition. If $F\subset \C$ is a finite extension of $F_0$ and $M\in X(F)$ and if $[F:F_0]\leq d$ and $H_{\La}(M)\leq c$, there is $i$ ($1\leq i\leq n$) such that $F_i\subset F$ and such that $M$ is the image of $M_i$ under $X(F_i)\to X(F)$. 

\medskip

Remark \ref{Q2}.1. We have (i) $\Rightarrow$ (ii) by \ref{H+}, (iii) $\Rightarrow$ (ii) by \ref{Tandh}, and (i) $\Rightarrow$ (iii) by \ref{HN+}.

\end{sbpara}

Now we do not assume $G$ is reductive. 

\begin{sbpara}\label{Q3} {\bf Question 3.}
Assume $\La$ is ample (\ref{ample}). Are the following (1) and (2) true?

\medskip

(1) For any finite extension $F$ of $F_0$ in $\C$, $\{M\in X(F)\;|\; H_{\La}(M)\leq C\}$ is finite.

(2) (A stronger form of (1).) Fix $d\geq 1$ and $c\in\R_{>0}$. Then 
there is a  finite number of pairs pairs $(F_i, M_i)$ ($1\leq i\leq n$), where $F_i$ are finite extensions of $F_0$ in $\C$  and $M_i\in X(F_i)$ satisfying the following condition. If  $F$ is a finite extension of $F_0$ in $\C$ such that $[F:F_0]\leq d$ and if $M_i\in X(F)$ and $H_{\La}(M)\leq c$, then there is $i$ ($1\leq i\leq n$) such $F_i\subset F$ and such that $M$ is the image of $M_i$ under $X(F_i)\to X(F)$. 

\medskip
Remark \ref{Q3}.1.  Question 3 asks whether the motive version of the finiteness theorem of Northcott on the usual height is true.  

\medskip

Remark \ref{Q3}.2. 
It seems that through the analogies between Situation (I) and Situation (II) (resp. Situation (I) and Situation (III)), the result \ref{HN+} (resp. \ref{H+}) supports the answer Yes to Question 3.

\end{sbpara}

\subsection{Speculations on Vojta conjectures}\label{ss:Voj}  

\begin{sbpara} In Situation (I), for a finite extension $F$ of $F_0$ and for a finite set $S_0$ of places of $F_0$ containing all archimedean places of $F_0$, let $D_{S_0}(F/F_0)$ be the norm of the component of different ideal of $F/F_0$ outside $S_0$.

In Situation (II), we consider the Nevanlinna analogue of it
 $$N_{\text{Ram}(\Pi)}(r):= \sum_{x\in B, 0<|\Pi(x)|<r} (e_x(\Pi)-1)\log(r/|\Pi(x)|) + \sum_{x\in B,\Pi(x)=0} (e_x(\Pi)-1)\log(r),$$ where $e_x(\Pi)$ is the ramification index of $B$ over $\C$ at $x$.
\end{sbpara}

\begin{sbpara}\label{Qvo}  {\bf Question 4 (1).} Assume we are in Situation (I). Assume $\La$ is ample (\ref{ample}). Is the following statement true?

 Fix a finite set of places $S_0$ of $F_0$ which contains all archimedean places of $F_0$.
Then there is a constant $c\in \R_{>0}$ such that 
$$(\prod_{v\in\Sig_S(M)} \sharp(\F_v))\cdot D_{S_0}(F/F_0) \geq c^{[F:F_0]}\cdot  H_{\spadesuit}(M)H_{\La}(M)^{-1}$$ 
for any finite extension $F$ of $F_0$ in $\C$ and for any $M\in X_{\text{gen}}(F)$.
Here $S$ denotes the set of all places of $F$ lying over $S_0$, $\Sig_S(M)$ denotes the set of  places of $F$ which do not belong to $S$ and  at which $M$ has bad reduction,  and $\F_v$ denotes the residue field of $v$.

\medskip

Remark \ref{Qvo}.1. There are variants of this question. We may add the assumption $[F:F_0]\leq d$ for some fixed $d\geq 1$. We may replace $\prod_{v \in S} \sharp(\F_v)$ by $H_{\heartsuit, S}(M)$ assuming that $S_0$ contains all non-archimedean places of $F_0$ at which $M$ has bad reduction.

The answer Yes to any one of these questions is an analogue of a conjecture of Vojta in \cite{Vo} (see Section 24 of \cite{Vo}).

\medskip

Remark \ref{Qvo}.2. If we follow the analogy with conjectures in \cite{Vo}, the reader may think that $\epsilon>0$ should be fixed and $H_{\La}(M)$ should be replaced by $H_{\La}(M)^{\epsilon}$. However, we have $H_{\La}(M)^{\epsilon}= H_{\La'}(M)$ where $\La'$ is obtained from $\La$ by replacing $c(i)$ with $\epsilon c(i)$ and replacing $t(w,d)$ with $\epsilon t(w,d)$ (then $\La'$ is also ample). 

\medskip

Remark \ref{Qvo}.3. In the case $G$ is reductive, by \ref{spadepos}, we may have a variant of Question 4 (1)  replacing $H_{\spadesuit}(M)H_{\La}(M)^{-1}$ by $H_{\spadesuit}(M)^{1-\epsilon}$ for a fixed $\epsilon>0$. 

\medskip

{\bf Question 4 (2).} Assume we are in Situation (II). 
 Let $\La$ be ample (\ref{ample}). Fix $B$.  Is the following statement  true? 

If $f\in \sM_{\text{hor},\text{gen}}(B, X(\C))$ and $c\in \R$, we have 
  $$N^{(1)}_{f,\heartsuit}(r)+N_{\text{Ram}(\Pi)}(r) \geq_{\text{exc}} T_{f,\spadesuit}(r)-T_{f,\La}(r)+ c$$ where $\geq_{\text{exc}}$ means that the inequality $\geq $ holds for any $r$ outside 
  some subset of $\R_{\geq 0}$ of finite Lebesgue measure. 
Here $N^{(1)}_{f,\heartsuit}(r)$ is as in \ref{N^1}. 

 \medskip
 
 Remark \ref{Qvo}.4. Question 4 (2)  treats the Hodge-Nevanlinna version of the conjecture of Griffiths in the usual Nevanlinna theory treated in \cite{Vo} Section 14, Section 26. 
 
 \medskip
 
 Remark \ref{Qvo}.5. There is a variant of Question 4 (2) in which we replace $N^{(1)}_{f, \heartsuit}(r)$ by $N_{f,\heartsuit}(r)$. 
 
 \medskip
 
 Remark \ref{Qvo}.6. If $G$ is reductive, by \ref{spadepos}, Question 4 (2) becomes equivalent to the question in which we replace $T_{f,\spadesuit}(r)-T_{f,\La}(r)$ in the above by $(1-\epsilon)T_{f,\spadesuit}(r)$ ($\epsilon>0$). 
\medskip

{\bf Question 4 (3).}  Assume we are in Situation (III). Let $\La$ be ample (\ref{ample}). Fix $C$. Is the following statement true? 

There is $c\in \R$ such that  
 $\sharp(\Sig(\sH)) \geq h_{\spadesuit}(\sH)+c$
for any $\sH\in \sM_{\text{hor},\text{gen}}(C, X(\C))$. 

Here $\Sig(\sH)$ denotes the set of points in $C$ at which $\sH$ degenerates.

\medskip

Remark \ref{Qvo}.7. In Situation (III), various results are known concerning the comparison of $\sharp(\Sig(\sH))$ and kinds of heights of $\sH$ for a variation of Hodge structure $\sH$ on $C$ with log degeneration. For example, see Proposition 2.1 of \cite{VZ}. 

\medskip

Remark \ref{Qvo}.8. If the answer to Question 4 (2) is Yes, the answer to Question 4 (3) is Yes. This follows from \ref{Tandh} and \ref{Tandh2}.

\medskip

Remark \ref{Qvo}.9.  If $X(\C)$ is of dimension $1$, Question 4 (3) has an affirmative answer. In this case, if $f$ is not constant, then for $R:= \Sigma(\sH)$, we have a canonical injective map $f^*(\Omega^1_{\bar X(\C)}(\log))\to \Omega^1_C(\log R)$ and hence the comparison of degrees gives $h_{\spadesuit}(\sH) \leq \text{deg}(\Omega^1_C(\log R))= 2g_C-2+\sharp(R)$.

Remark \ref{Qvo}.10. There is a variant of Question 4 (3) in which we replace $\sharp(\Sig(\sH))$ by $h_{\heartsuit}(\sH)$. 
\end{sbpara}

\begin{sbpara}\label{gentype} 
In this Question 5, we assume $G$ is reductive. 

\medskip

{\bf  Question 5 (1).}
Are the following (i) and (ii) equivalent?

\medskip

(i) $X(F)$ is finite for any finite extension $F$ of $F_0$ in $\C$.

(ii) Any $f\in \sM_{\text{hor}}(\C, X(\C))$ is constant.

\medskip

{\bf Question 5 (2).}
Are the following (i)--(iii) equivalent? 

\medskip

(i) $X_{\text{gen}}(F)$ is finite for any finite extension $F$ of $F_0$ in $\C$ in Situation (I). 

(ii)  Any $f\in \sM_{\text{hor},\text{gen}}(\C, X(\C))$ is constant.  

(iii) There is an ample $\La$ (\ref{ample}) such that  $h_{\spadesuit}(\sH) -h_{\heartsuit}(\sH) \geq h_{\La}(\sH)$ in Situation (III)  for  any $C$ and for any $\sH\in \sM_{\text{hor},\text{gen}}(C, X(\C))$.

\medskip

{\bf Question 5 (3).} 
Is the following true? 

\medskip
 
Let $F$ is a finite extension of  $F_0$ in $\C$ and let $S$ be a finite set of places of $F$ contining all archimedean places of $F$. Then there are only finitely many $M\in X(F)$ which are of good reduction outside $S$.

\medskip

Remark \ref{gentype}.1. The answer Yes to Question 5 (3)  is a generalization of Shafarevich conjecture for abelian varieties proved by Faltings (\cite{Fa0}) to $G$-motives.

Remark \ref{gentype}.2. The last section of Koshikawa \cite{Ko1} contains an affirmative result on the Shafarevich conjecture for motives assuming the finiteness of the number of motives with bounded heights and assuming the semi-simplicity of the category of pure motives. 

\medskip

Renark \ref{gentype}.3. Via the analogy in \ref{Dioconj} below, Proposition \ref{Sh0} supports the answer Yes to Question 5 (3).

\end{sbpara}

\begin{sbrem}\label{Dioconj} 

The following conjecture of Bombieri-Lang-Vojta is well known in Diophantine geometry.

Conjecture \ref{Dioconj}.1. Let $F_0\subset \C$ be a number field, let $V$ be a scheme of finite type over  $O_{F_0}$, and let $V'$ be a closed subscheme of $V$. 
Then the following conditions are equivalent.

(i) The image of any morphism $\C\to V(\C)$ of complex analytic spaces is contained in $V'(\C)$.

(ii) For any finite extension $F$ of $F_0$ and for any finite set $S$ of places of $F$ containing all archimedean places of $F$, $V(O_S) \smallsetminus V'(O_S)$ is finite.

Question 5 basis on the analogy with this conjecture.

In Question 5 (1), we consider an arithmetic version $\bar X$ of $\bar X(\C)$ (if it exists)  in place of $V$ and we consider $\bar X \smallsetminus X$ in place of $V'$. 

In Question 5 (2), we consider $\bar X$ in place of  $V$.

In Question 5 (3),  we consider $X$ in place of $V$ taking the empty set as $V'$.

\end{sbrem}

The following Question 6 presents inequalities which are similar to those in Question 4, but the difference is that the generic parts ($X_{\text{gen}}(F)$ etc.) are considered in Question 4 whereas the total space ($X(F)$ etc.) are considered in Question 6.

\begin{sbpara}\label{hyperb} {\bf Question 6.} 

\medskip

 {\bf (1)} Assume we are in Situation (I) and assume $G$ is reductive. Is there an ample $\La$ (\ref{ample}) satisfying the following?
Fix   a finite set of places $S_0$ of $F_0$ which contains all archimedean places of $F_0$. 
Then there exists $c\in \R_{>0}$ such that

\medskip
 $(\prod_{v\in\Sig_S(M)} \sharp(\F_v))\cdot D_{S_0}(F/F_0)\geq c^{[F:F_0]}\cdot  H_{\La}(M)$ 

\medskip
\noindent
 for any finite extension $F$ of $F_0$ in $\C$  and for any $M\in X(F)$.

Here $S$ and $\Sig_S(M)$ are as in Question 4 (1).

\medskip

{\bf (2)} Assume we are in Situation (II) and  assume $G$ is reductive. Fix $B$. Is there an ample $\La$ satisfying the following?
If $f\in \sM_{\text{hor}}(B, X(\C))$ and $c\in \R$, we have
  $$N^{(1)}_{f,\heartsuit}(r)+ N_{\text{Ram}(\Pi)}(r)\geq_{\text{exc}} T_{f,\La}(r)+c.$$

\medskip
{\bf (3)} Assume we are in Situation (III) and assume $G$ is reductive. Fix $C$. Are there an ample $\La$
and $c\in \R$ such that 
 $$\sharp(\Sig(\sH)) \geq h_{\La}(\sH)+c$$
for any $\sH\in \sM_{\text{hor}}(C, X(\C))$?
Here $\Sig(\sH)$ denotes the set of points in $C$ at which $\sH$ degenerates.

\medskip

{\bf (4)} Assume we are in Situation (I). Let $S_0$ be a finite set of places of $F_0$ containing all archimedean places of $F_0$ and all non-archimedean places at which $M_b$ has bad reduction. Then is there an ample $\La$ and $c>0$ such that $$H_{\La}(M)\geq c^{[F:F_0]} H_{\heartsuit,S}(M)$$ for any finite extension $F$ of $F_0$ in $\C$ and for any $M\in X(F)$? Here $S$ denotes the set of all places of $F$ lying over $S_0$.

\medskip

{\bf (5)}  Assume we are in Situation (II). Fix $B$. Is there an ample $\La$ such that  $$T_{f,\La}(r) \geq_{\text{exc}} N_{f,\heartsuit}(r)+c$$ for any $f\in \sM_{\hor}(B, X(\C))$ and $c\in \R$?

\medskip

{\bf (6)}  Assume we are in Situation (III). Fix $C$. Is there an ample $\La$ such that $$h_{\La}(\sH) \geq h_{\heartsuit}(\sH)$$ for any $\sH\in \sM_{\hor}(C, X(\C))$?

\medskip

Remark \ref{hyperb}.1.  If the answer to Question 6 (2) is Yes, then the answer to  Question 6 (3) is Yes. This follows from \ref{Tandh} and \ref{Tandh2}. 

\medskip

Remark \ref{hyperb}.2. If the answer to  Question 6 (1) is Yes and if (i) $\Rightarrow$  (iv)  in Question 2 is true, we have the answer Yes
 to Question 5 (3). 

\medskip
Remark \ref{hyperb}.3.  If the answer to Question 6 (4) (resp. (5), resp. (6)) is true, by the method of \cite{Vo0}, we can show the equivalence between Question 6 (1) (resp. (2), resp. (3)) and its variant using $H_{\heartsuit,S}$ (resp. $N_{f,\heartsuit}$, resp. $h_{\heartsuit}$) in Remark \ref{Qvo}.1. (We move $K$ to prove this equivalence.) But we do not discuss this in this paper.

\end{sbpara}

\subsection{Speculations on the number of motives of bounded height}\label{ss:Manin}

\begin{sbpara}\label{331} Fix an ample $\La$ (\ref{ample}). Consider Situation  (I). We fix a finite extension $F$ of $F_0$ in $\C$.

 For $t\in \R_{>0}$, define
 $$N(H_{\La},t) := \sharp\{M\in X(F)\;|\; H_{\La}(M)\leq t\}\quad N_{\text{gen}}(H_{\La},t) := \sharp\{M\in X_{\text{gen}}(F)\;|\;H_{\La}(M)\leq t\}.$$
 We expect that these numbers are finite (Question 3 in Section 4.1). 
\end{sbpara}
\begin{sbpara}\label{Q7}  {\bf Question 7.} Do we have 
$$N(H_{\La},t) = c t^a \log(t)^b(1+o(1)), \quad N_{\text{gen}}(H_{\La}, t)= c' t^{a'}\log(t)^{b'}(1+o(1))$$
for some constants $a, a'\in \R_{\geq 0}$, $b, b'\in \frac{1}{2}\Z$, $c, c'\in \R_{\geq 0}$? Do we have $a,a'\in \Q$ in the case $c(i), t(w,d)\in \Q$ in the definition of $\La$?

\end{sbpara}

\begin{sbpara}\label{alpha}
Define 
$\alpha\in \R_{\geq 0}\cup \{\infty\}$ to be the inf of all $s\in \R_{\geq 0}$ satisfying the following condition (i).

(i) For any connected projective smooth curve $C$ over $\C$ and for any $\sH\in \sM_{\text{hor}}(C, X(\C))$, we have $$sh_{\La}(\sH) + h_{\spadesuit}(\sH) - h_{\heartsuit}(\sH)\geq 0.$$
\medskip
($\alpha$ is defined to be $\infty$ if such $s$ does not exist).

\end{sbpara}

\begin{sbrem} The inequality $sh_{\La}(\sH) + h_{\spadesuit}(\sH) \geq h_{\heartsuit}(\sH)$ can be written also as 
$\text{deg}(f^*\Omega^1_{\bar X(\C)}) + s h_{\La}(\sH)\geq 0$ ($f=\sH\in \sM_{\hor}(C, X(\C))$) by using the partial toroidal compactification $\bar X(\C)$ of $X(\C)$. 
\end{sbrem}

The next question arises following the analogy with the conjecture \cite{BM} of Batyrev-Manin. 
\begin{sbpara}\label{Q8}

{\bf Question 8.} Assume that $X_{\text{gen}}(F)$ is not empty.
Do we have $$\lim_{t\to \infty} \frac{\log(N_{\text{gen}}(H_{\La},t))}{\log(t)} =  \alpha<\infty?$$ 

\end{sbpara}

\begin{sbpara} {\bf Example.} In Example \ref{MEx4}, let $F_0=\Q$, $M_1=\Z\oplus \Z(1)^n$ for some $n\geq 1$, let the polarization $p_0$ on $\Z$ be the standard one, and let the polarization $p_{-2}$ on $\Z(1)^n$  be the direct sum of the standard polarizations on $\Z(1)$.
 Let $\sQ$, $G$ and $\sG$ be as in \ref{MEx4}. Then 
$\sQ={\bf G}_m$, $\sG={\bf G}_a^n$,  and $G$ is the semi-direct product of $\sG$ and $\sQ$ in which the action of $\sQ$  on $\sG$ via inner automorphisms is the natural action of ${\bf G}_m$ on ${\bf G}_a^n$.

We have $X(\C)= \text{Ext}^1_{\Q\text{MHS}}(\Z,\Z(1)^n)= (\C^\times)^n$. $X(F)= \text{Ext}^1_{\text{MM}(F, \Z)}(\Z, \Z(1)^n)$, where $\text{Ext}^1_{\text{MM}(F,\Z)}$ is the category of mixed motives with $\Z$-coefficients over $F$. Let $X'(F)=(F^\times)^n$. We have a canonical map $X'(F) \to X(F)$ and a philosophy of mixed motives suggests that this map should be bijective. We will compute the number of motives of bounded height replacing $X(F)$ by $X'(F)$ and replacing $X_{\text{gen}}(F)$ by the inverse image $X'_{\text{gen}}(F)$ of $X_{\text{gen}}(F)$ in $X'(F)$. Then $(a_i)_{1\leq i\leq n}\in X'(F)=(F^{\times})^n$ belongs to $X'_{\text{gen}}(F)$ if and only if $a_1, \dots, a_n$ are linearly independent in the $\Q$-vector space $F^\times \otimes \Q$. For $a\in (F^\times)^n$, let $M_a\in X(F)$ be its image. 

We take $\La$ (\ref{331}) such that $t(0,2)=1$. 
We have 
$$H_{\La}(M_a)= \exp(\sum_v (\sum_{i=1}^n \log(|a_i|_v)^2)^{1/2}),$$ where $v$ ranges over all places of $F$. (This is proved as Part I, 1.7.8 which treats  the case $n=1$.)

We prove  $$\lim_{t\to \infty} \frac{\log(N_{\text{gen}}(H_{\La},t))}{\log(t)} = 1= \alpha.$$

First we prove $\alpha=1$. If  $\sH\in \sM_{\text{hor}}(C, X(\C))$ corresponds to the extension $0\to \Z(1)^n \to \sH \to \Z\to 0$ defined by a family $(f_i)_{1\leq i\leq n}$ of non-zero elements $f_i$ of the function field of $C$, then
$$h_{\La}(\sH)=\sum_{x\in C}   (\sum_{i=1}^n \text{ord}_x(f_i)^2)^{1/2},$$
$$h_{\spadesuit}(\sH)=0,$$
$$h_{\heartsuit}(\sH)= \sum_{x\in C}  \text{GCD}\{|\text{ord}_x(f_i)|\;|\; 1\leq i\leq n\},$$
where GCD is the greatest common divisor. 
Hence $sh_{\La}(\sH)+ h_{\spadesuit}(\sH) \geq h_{\heartsuit}(\sH)$ means $sh_{\La}(\sH)\geq h_{\heartsuit}(\sH)$. We have clearly $h_{\La}(\sH)\geq h_{\heartsuit}(\sH)$. On the other hand, if we take $f_i=1$ for $2\leq i\leq n$, we have 
$h_{\La}(\sH)=h_{\heartsuit}(\sH)= \sum_{x\in C} |\text{ord}_x(f_1)|$. These prove $\alpha=1$.

Next we prove $\lim_{t\to \infty} \log(N_{\text{gen}}(H_{\La},t))/\log(t) =1$.

 Let $V:={\bf G}_m^n$, $\bar V_1=({\bf P}^1)^n$, and let $D_1$ be the divisor $\bar V_1\smallsetminus V$ with simple normal crossings. Let 
$\bar V_2$ be the  blowing-up of $\bar V_1$ along all intersections of two irreducible components of $D_1$ (i.e. the bowing-up of $\bar V_1$ by the product ideal of $\sO_{\bar V_1}$ of the ideals which define intersections of two irreducible components of $D_1$). Let $D_2$ be the divisor $ \bar V_2\smallsetminus V$ on $\bar V_2$ with simple normal crossings.
Then for $i=1,2$, the height function $H_i$  associated to the divisor $D_i$ on $\bar V_i$
is described as follows. 
For $a=(a_i)_{1\leq i\leq n}\in (F^\times)^n$, 
$$H_{(1)}(a)= \exp(\sum_v \sum_{i=1}^n |\log(|a_i|_v)|),\quad H_{(2)}(a)= \exp(\sum_v \max\{|\log(|a_i|_v)|\;|\; 1\leq i \leq n\},$$
where $v$ ranges over all places of $F$. Since $D_i$ is rationally equivalent to $-K_i$ where $K_i$ is the canonical divisor of $\bar V_i$, we have the work \cite{BT} shows that
$\lim_{t\to \infty} \log(N(H_{(i)},t))/\log(t)=1$
for $i=1,2$ by \cite{BT}. 
Since $H_{(1)}(a) \geq H_{\La}(M_a) \geq H_{(2)}(a)$, this  shows 
 $\;\lim_{t\to\infty} \log(N(H_{\La}, t))/\log(t)=1$. We can obtain $\;\lim_{t\to \infty} \log(N_{\text{gen}}(H_{\La},t))/\log(t)=1$ from it easily.

\end{sbpara}

\begin{sbpara}\label{Q10.0} 

 The following Question 9 is a refined version of Question 8. In Question 9, we consider the type of monodromy.

Consider the quotient set $\Sig/\sim$ of $\Sig$ (\ref{Sig}) where $\sim$ is the following equivalence relation. For nilpotent operators $N, N'\in \Lie(\sG)$, $\R_{\geq 0}N\sim \R_{\geq 0}N'$ if and only if there are $(g,t)\in \sG(\C)$ and $c\in \C^\times$ such that $N'=c\text{Ad}(g)(N)$ in $\Lie(\sG)_\C$. Let $\Sig'$ be a subset of $\Sig/\sim$ which contains the class of the cone $\{0\}$.

\medskip

Let $S$ be a finite set of places of $F$ which contains all archimedean places of $F$ and all non-archimedean places of $F$ at which $M_b$ is of bad reduction. 
We define the height function $h_{\heartsuit, \Sig'}(\sH)$ and the counting function $N_{\text{gen},S, \Sig'}(H, t)$, which take the shapes of the monodromy operators into account.

For  a connected projective smooth curve $C$ and for $\sH\in \Mh(C, \bar X(\C))$, 
define $h_{\heartsuit, \Sig'}(\sH)= \sum_x  h_{\heartsuit, x}(\sH)\in \Z_{\geq 0}$ where $x$ ranges over all points of $C$ which satisfy the following condition (i).

\medskip

(i) The image of $x$ under the composite map $$C\overset{\sH}\to \bar X(\C) \to \sG(\Q)\bs (\Sig \times (\sG({\bf A}^f_\Q)/K))\to \sG(\Q)\bs \Sig \to \Sig/\sim$$ belongs to $\Sig'$. Here $\sG(\Q)$ acts on $\Sig$ by conjugation.

\medskip

Let $N_{\text{gen}, S, \Sig'}(H_{\La}, t)$ be the number of $M\in X_{\text{gen}}(F)$ which appear in the definition of $N_{\text{gen}}(H_{\La},t)$ and which satisfy  the following condition (ii) at each place $v\notin S$ of $F$.

(ii) There is an element $\text{class}(\R_{\geq 0}N)$ (with $N\in \Lie(G_1)$ nilpotent) 
of $\Sig'$ which satisfies the following:
Let $p=\text{char}(\F_v)$. For any prime number $\ell\neq p$, there are $t,c\in \bar \Q_{\ell}^\times$ and  an  isomorphism $\bar \Q_{\ell}\otimes_{\Q_{\ell}} M_{et,\Q_{\ell}}\cong H_{0,\bar \Q_{\ell}}$ of $\bar \Q_{\ell}$-vector spaces which preserves $W$ and via which $\langle\;,\;\rangle_w$ on $\bar \Q_{\ell} \otimes_{\Q_{\ell}} \gr^W_wM_{et,\Q_{\ell}}$ corresponds to $t^w\langle \; , \;\rangle_{0,w}$ on $\gr^W_wH_{0,\bar \Q_{\ell}}$ for any $w\in \Z$ and  
the local monodromy operator $N'_v$ corresponds to $cN$.   Furthermore, there are $t,c\in \bar F_v^\times$ and  an  isomorphism $\bar F_v \otimes_{F_{v,0,\text{ur}}} D_{\pst}(F_v, M_{et,\Q_p})\cong H_{0, \bar F_v}$ of $\bar F_v$-vector spaces 
which preserves $W$ and via which  $\langle\;,\;\rangle_w$ on $\bar F_v\otimes_{F_{v,0,\text{ur}}} \gr^W_wD_{\pst}(F_v, M_{et,\Q_p})$ corresponds to $t^w\langle \; , \;\rangle_{0,w}$ on $\gr^W_wH_{0,\bar F_v}$  for any $w\in \Z$ and 
the local monodromy operator $N'_v$ corresponds to $cN$.

\end{sbpara}

\begin{sbpara}　The above $h_{\heartsuit,\Sig'}$ can be expressed 
in the form 
$$h_{\heartsuit,\Sig'}(\sH)= \sum_{x\in C} e(x), \quad (f^*I_{\Sig'})_x=m_x^{e(x)}$$
where $I_{\Sig'}$ is the invertible ideal of $\sO_{\bar X(\C)}$ which coincides with $I_{\bar X(\C)}$ in \ref{barX4} at points of  $\bar X(\C)$ whose classes in $\Sig/\sim$ belong to $\Sig'$ and coincides with $\sO_{\bar X(\C)}$ at  other points of $\bar X(\C)$. 
 
 \end{sbpara}
 
\begin{sbpara}\label{Q10}  Let the notation be as in \ref{Q10.0}. Let $F$ be a number field. Let $S$ be a finite set of places of $F$ containing all archimedean places. 
Define the modified version $$ \alpha_{\Sig'} \;\;   \quad \text{of}\quad   \alpha$$ by replacing  
$$sh(\sH)+ h_{\spadesuit}(\sH) -h_{\heartsuit}(\sH)\geq 0\quad \text{by} \quad sh(\sH)+ h_{\spadesuit}(\sH) - h_{\heartsuit,\Sig'}(\sH)\geq 0.$$ 

\medskip
{\bf Question 9}. Assume $X_{\text{gen}}(F)\neq \emptyset$. Do we have 
$$\lim_{t\to \infty} \frac{\log(N_{\text{gen},S, \Sig'}(H_{\La},t))}{\log(t)} =  \alpha_{\Sig'}?$$

\end{sbpara}

\end{document}